\batchmode
 \documentclass[11pt]{amsart}
 \usepackage{amssymb}
 \usepackage{amsmath}
\usepackage{amscd} 
\begin{document}
    \newtheorem{Theorem}{Theorem}[section]
    \newtheorem{Proposition}[Theorem]{Proposition}
    \newtheorem{Lemma}[Theorem]{Lemma}
    \newtheorem{Corollary}[Theorem]{Corollary}
   \newcommand{\ra}{\rightarrow}   
   \newcommand{\ul}{\underline}
   \newcommand{\ve}{\varepsilon}
   \newcommand{\varp}{\varphi} 
   \newcommand{\al}{\alpha}
    \newcommand{\aut}{{\rm Aut}} 
     \newcommand{\mbf}{\mathbf t} 
     \newcommand{\Z}{\boldsymbol{Z}} 
     \newcommand{\C}{\boldsymbol{C}} 
     \newcommand{\R}{\boldsymbol{R}} 
     \newcommand{\B}{\boldsymbol{P}} 
     \newcommand{\Q}{\boldsymbol{Q}}

\title{Anti-self-dual bihermitian structures on Inoue surfaces} 

\author{A. Fujiki \and M. Pontecorvo}
\footnote{Research partially supported by Grant-in-Aid for Scientific Research (No.18340017) MEXT, 
and PRIN, Italy: Metriche riemanniane e variet\'a differenziabili} 

\begin{abstract} In this article we show that 
any hyperbolic Inoue surface (also called Inoue-Hirzebruch surface of even type) 
admits anti-self-dual bihermitian structures. 
The same result holds for any of its small deformations 
as far as its anti-canonical system is non-empty.  
Similar results are obtained for parabolic Inoue surfaces.  
Our method also yields a family of anti-self-dual {\em hermitian} metrics 
on any half Inoue surface. 
We use the twistor method of Donaldson-Friedman \cite{df} for the proof. 
\end{abstract}

\maketitle 

\section{Introduction} 

Let $M$  be a compact smooth oriented four dimensional manifold.  
A {\em bihermitian structure} on  $M$  is a triple $\{[g],J_1,J_2\}$ 
consisting of a conformal class $[g]$ of a Riemannian metric $g$ and 
two complex structures $J_i, i=1,2$, such that $([g],J_i)$ define a conformal 
hermitian structure on  $M$, and $J_i$ are compatible with the oritentation of $M$, and 
are inequivalent to each other in the sense that $J_1\neq \pm J_2$ 
when considered as integrable almost complex structures.  
It is called an {\em anti-self-dual bihermitian structure} 
if, further, $(M,[g])$ is an anti-self-dual structure in the sense of \cite{ahs}.  
Note that such a structure is always (twisted) generalized K\"ahler in the sense of \cite{gu04} 
if the hyperhermitian case is excluded as we do in this paper. 

The second-named author \cite{pnt6} with a supplement by Dloussky \cite{dl} 
has shown the following: 
Let  $M$  be a compact smooth oriented four-manifold admitting 
an anti-self-dual bihermitian structure $\{[g],J_1,J_2\}$ 
which is not hyperhermitian.  
Let  $S=(M,J_i)$ be the associated compact complex surface. 
Then the anti-canonical system $|-K_S|$ admits a disconnected member. 
In particular if  $S$  is minimal, $S$ is either a Hopf surface, a parabolic Inoue surface, 
or a hyperbolic Inoue surface.  In general,  $S$ is obtained 
from its minimal model $\bar{S}$  by blowing up a finite number of points 
on a fixed anti-canonical divisor.  
Note that all these surfaces have the underlying $C^{\infty}$ manifold  
$M=M[m]:=(S^1\times S^3)\#m\bar{\B}^2$ 
where  $\bar{\B}^2$ denotes the complex projective plane with reversed oritentation. 

In this paper we show the existence 
of anti-self-dual bihermitian structures on any hyperbolic Inoue surface and 
also on any of its small deformations 
preserving the unique anti-canonical divisor on it.  
This is thought of as a partial converse of the above result of \cite{pnt6}. 
One of our main results is more precisely as follows.    
Let $\{S,{}^{\mbf}S\}$ be any pair of hyperbolic Inoue surfaces $S$ and its 
transposition ${}^{\mbf}S$ with second Betti number $m$.  
Then there exists a family  
of anti-self-dual bihermitian structures $\{[g]_t,J_t,{}^{\mbf}J_t\}$ 
on  $M[m]$ with real smooth $m$-dimensional parameter  $t$ 
such that $(M[m],J_t) \cong S$  and $(M[m],{}^{\mbf}J_t) \cong {}^{\mbf}S$ 
independently of the parameter $t$ (Theorem \ref{ma}). 
The same result also holds 
for any properly blown-up hyperbolic Inoue surface (cf. \S 3) and 
for any of its small deformations with an effective (and disconnected) anti-canonical divisor 
(Theorem \ref{mad}). 
Moreover, we prove similar results also for parabolic Inoue surfaces (Theorem \ref{map}). 
Finally, the same method also yields the existence of 
a real $m$-dimensional family of anti-self-dual {\em hermitian} structures 
on any (properly blown-up) half Inoue surfaces (Theorem \ref{mah}), 
which never carry anti-self-dual bihermitian structures. 

For the proof we use a variation of the twistor method 
due to Donaldson-Friedman \cite{df} in the spirit of \cite{kp}. 
Namely instead of an anti-self-dual bihermitian triple $\{[g], J_1, J_2\}$ 
we construct the twistor space corresponding to $[g]$ and two pairs of 
mutually conjugate elementary, i.e., degree-1, surfaces  $\{S^+_i,S^-_i\}$ on it  
giving rise to $\pm J_i, i=1,2$.  
The twistor spaces  
associated to self-dual metrics on  $m\B^2$ constructed by Joyce \cite{jo}, 
studied in detail in \cite{fjj}, play a crucial role in our construction. 

All our examples yield anti-self-dual (bi)hermitian - 
also (twisted) generalized K\"ahler - and locally conformally K\"ahler structures, 
which are new except possibly for the parabolic Inoue case and its deformations, 
which should be compared with the examples of LeBrun [23]. 

We now give a brief description of each section: after 
some preliminaries on deformation theory of pairs of complex spaces in Section 2, 
we recall in  Section 3 
basic properties of hyperbolic, half and parabolic Inoue surfaces and 
study the Kuranishi family of deformations of associated pairs. 
All the Inoue surfaces are known to be obtained 
as a deformation of a singular toric surface $\hat{S}$ \cite{na84}.  
In Section 4 we formulate the result in terms of the Kuranishi family 
of deformations of a pair.  
This concludes the first half of the paper and will be used constantly 
in the following sections where we take up our construction of twistor spaces. 

First in Section 5 
we recall the basic properties of a Joyce twistor space $Z$ according to \cite{fjj} 
and explain how to obtain its singular model $\hat{Z}$ 
together with the natural anti-canonical divisor $\hat{S}$ by analytic modifications. 
Then in Section 6 we study the local structures of the singularities of the pair 
$(\hat{Z},\hat{S})$ and its automorphism group. 
We state our main results in Section 7 and their proof will be subsequently 
provided in Section 8.  
The main point is to show the vanishing of the obstructions for 
smoothing of the pairs $(\hat{Z},\hat{S})$. 
Technically, this section is the most delicate part of the paper. 
In Section 9 we prove that a twisted version of our previous construction yields 
anti-self-dual {\em hermitian} metrics on half Inoue surfaces. 
Finally, in Section 10 we summarize differential geometric implications 
of our results including the relations with generalized K\"ahler and locally conformally 
K\"ahler structures.

\section{Preliminaries} 
\vspace{5 mm}

\ul{Sheaves of logarithmic forms}

For a complex space $X$ we shall denote by $\Omega _{X}$ 
the sheaf of germs of holomorphic $1$-forms on  $X $ and by 
$\Theta_{X}$ the sheaf of germs of holomorphic vector fields on $X $. 
$\Theta_{X}$ is the dual of $\Omega _{X}$. 

Let now  $Y$  be a reduced Cartier divisor on  $X$.  
In our case of interest  $X$  is smooth or with 
at worst normal crossings singularities 
along a smooth connected hypersuface $D$ and  
$Y$  has only mild singulariteis.  

We define the {\em sheaf} $\Omega _{X}(\log Y)$ 
{\em of logarithmic $1$-forms} on  $X$  along $Y$ 
to be the sheaf of germs of meromorphic $1$-forms $\omega $ on  $X$  
such that for any local equation  $f=0$ of $Y$  in  $X$,
both  $f\omega $ and  $fd\omega $ are holomorphic $1$-forms on  $X$. 
On the other hand, we define the sheaf  $\Theta_{X}(-\log Y)$ 
of logarithmic vector fields on  $X$  along $Y$ 
to be the sheaf of germs of holomorphic vector fields  $v$ on  $X$  
such that  $v(f)/f$  is again holomorphic with  $f$ as above, namely 
germs of those holomorphic vector fields which are tangent to $Y$. 
It is easy to see that they are coherent analytic sheaves on  $X$.  
At a smooth point of $X$ at which  $Y$  has at worst normal crossings singularities 
both the sheaves $\Omega _{X}(\log Y)$ and 
$\Theta_{X}(-\log Y)$ are free (cf.\ \cite{del}). 

We also consider the subsheaf 
\begin{equation}\label{om'} 
\Omega'_{X}(\log Y) 
\end{equation}  
of $\Omega_{X}(\log Y)$ locally generated by  $\Omega_{X}$ and the element  $df/f$  
for a defining equation  $f=0$ of $Y$ in $X$.  

\vspace{2 mm} 
\ul{Deformation theory} 

Let  $X$  and  $Y$  be as above and suppose that they are compact. 
A {\em deformation} of the pair $(X,Y)$ is a triple 
\begin{equation}\label{kur} 
 f: ({\mathcal X},{\mathcal Y}) \ra  T,\ X_o=X,\ o\in T 
\end{equation}  
where  $f$  is a flat morphism ${\mathcal X} \ra  T$ of complex spaces 
which induces a flat morphism ${\mathcal Y} \ra  T$; 
especially ${\mathcal Y}$ is a Cartier divisor on  ${\mathcal X}$.  
A {\em log-deformation} of the pair $(X,Y)$ is a deformation of $(X,Y)$ 
such that any local Cartier irreducible component of $Y$  remains locally irreducible 
under deformations.  
Especially the number of irreducible components of $Y$ remains the same under deformations. 

As usual we have the notions of the Kuranishi family (semiuniversal family) 
and the universal family of such deformations.  
Suppose that (\ref{kur}) is a Kuranishi family 
of log-deformations (resp.\ deformations) of $(X,Y)$.  
In this case we call the base space $T$ the associated {\em Kuranishi space}.  
The cohomological description of the infinitesimal deformation spaces of $(X,Y)$, 
the obstrcution space for the deformations and the space of 
infinitesimal automorphisms, are given by the following: 

\begin{Proposition}\label{st} 
The Kuranishi space  $T$  is smooth if 
$Ext^2_{O_X}(\Omega_X(\log Y),O_X)=0$ 
{\em (}resp.\ 
$Ext^2_{O_X}(\Omega'_X(\log Y),O_X)=0${\em )}. 
$Ext^1_{O_X}(\Omega_X(\log Y),O_X)$ 
{\em (}resp.\ 
$Ext^1_{O_X}(\Omega'_X(\log Y),O_X)${\em )} 
is naturally identified with the tangent space of  $T$  
at the reference point.  Moreover, if  $X$  is weakly normal, i.e., 
the Riemann extention theorem holds on  $X$ we have 
\[ Ext^0_{O_X}(\Omega_X(\log Y),O_X)=
Ext^0_{O_X}(\Omega'_X(\log Y),O_X)=
H^0(X,\Theta_X(-\log Y)) \] 
and if these vector spaces vanish, 
the Kuranishi family (\ref{kur}) is universal in a small neighborhood of $o$. 
\end{Proposition} 

In order to determine these vector spaces 
we use the local to global spectral sequence for $Ext$ functors. 
For instance in the case of log-deformations this takes the following form 
\begin{equation}\label{spc} 
E^{p,q}_2:=H^p(X, {\mathcal Ext}^q_{O_X}(\Omega_{X}(\log Y),O_X)) 
\Longrightarrow Ext^{p+q}_{O_X}(\Omega_{X}(\log Y),O_X) 
\end{equation}  
giving rise to the five term exact sequence in case  $X$  is weakly normal:
\begin{align}\label{lctog} 
0 & \ra &   H^1(X,\Theta_{X}(-\log Y)) 
  & \ra &  Ext^1_{O_X}(\Omega_{X}(\log Y),O_X) 
  & \ra  H^0(X,{\mathcal Ext}^1_{O_X}(\Omega_{X}(\log Y),O_X)) \\
  & \ra  &  H^2(X,\Theta_{X}(-\log Y))  \nonumber 
  & \stackrel{c}{\ra} &  Ext^2_{O_X}(\Omega_{X}(\log Y),O_X) & 
\end{align}  

\vspace{2 mm} 
\ul{A general lemma on Ext}

Let  $X$  be a complex space and  $Y$  a Cartier divisor on  $X$.  
Let  $N=[Y]|Y$ be the normal bundle of  $Y$  in  $X$.  
Then for any coherent analytic sheaf  $F$  on  $Y$ 
we have the following comparison theorem of Ext. 

\begin{Lemma}\label{extc} 
The notations being as above we have the following natural isomorphisms: 
\begin{gather}\label{isom}
Ext^i_{O_X}(F\otimes N,O_X) \cong 
Ext^{i-1}_{O_Y}(F, O_Y),  \ \ i\geq 0 , \\ 
{\mathcal Ext}^i_{O_X}(F\otimes N,O_X) \cong 
{\mathcal Ext}^{i-1}_{O_Y}(F, O_Y),  \ \ i\geq 0 ,   \label{liso}
\end{gather} 
where the second isomorphisms are those of $O_X$-modules. 
\end{Lemma} 

{\em Proof}.  
Applying \cite[p.72,Prop.2.9]{ak} for  $ E=N$ and  $G=O_{X}$ in the notation there, 
we obtain a spectral sequence 
\begin{equation}\label{sz} 
E^{p,q}_2:= 
Ext ^{p}_{O_Y}
(F,{\mathcal Ext}^q_{O_X}(N, O_X)) 
\Longrightarrow   
Ext^{p+q}_{O_X}(F\otimes N,O_X).  
\end{equation}  
Since $Y$ is a Cartier divisor in  $X$, we have 
${\mathcal Ext}^q_{O_X}(N, O_X)=0, q\neq 1$, 
and for $q=1$ 
\begin{eqnarray*} 
{\mathcal Ext}^1_{O_X}(N, O_X)  \cong 
{\mathcal Ext}^1_{O_X}(O_X(Y)\otimes_{O_X}O_Y, 
O_X) \\ 
\cong 
 {\mathcal Ext}^1_{O_X} 
(O_Y, O_X)\otimes_{O_X}O_X(-Y)
\cong N\otimes_{O_X} O_X(-Y) \cong O_Y
\end{eqnarray*} 
(cf.\ \cite[p.74,Prop.3.4]{ak} for the third isomorphism).  
Thus (\ref{sz}) yields the desired results. \hfill q.e.d.  

\vspace{2 mm} 
\ul{Automorphism groups}

For a compact complex space  $X$  and its subspace $Y$ denote by Aut$(X,Y)$ 
the group of automorphisms  of  $X$  preserving  $Y$.  
This has a natural structure of a complex Lie group with 
its Lie algebra naturally identified with  $H^0(\Theta_X(-\log Y))$.  
We denote by Aut$_0(X,Y)$ the identity component of $Aut(X,Y)$.  

\vspace{2 mm} 
\ul{Cycle of rational curves} 

A {\em cycle of rational curves} on a smooth surface is a compact connected curve 
$C$  which is either an irreducible rational curve with a single node or 
is a reducible curve with $k$ nodes 
whose irreducible components are nonsingular rational curves 
$C_i, 1\leq i\leq k, k\geq 2$, such that $C_i$ and $C_{i+1}$ intersects at a single point 
and there exists no other intersections, where $C_{k+1}=C_1$ by convention. 
We write such a $C$ as $C=C_1+\cdots +C_k$.  
Then the sequence $(b_1,\ldots ,b_k)$ 
of opposite self-intersection numbers  $b_i=-(C_i)^2$ is 
called the {\em weight sequence} of  $C$, 
which are considered modulo cyclic permutations and reversing the order. 

\vspace{2 mm} 

\ul{Toric surfaces}

Let  $G=\C^{*2}$ be the algebraic two-torus. 
Let  $S$  be a projective toric surface with fixed $G$-action with open orbit  $U$.  
The complement  $C:= S-U$  forms a cycle of rational curves 
$C=C_1+\cdots +C_k$ and is an element of the anti-canonical system $|-K_S|$ of $S$. 
We shall call  $C$ the {\em anti-canonical cycle} of  $S$ and 
denote the toric surface also by the pair  $(S,C)$. 
In this case $\Theta_S(-\log C)$ is free, i.e., 
\[ \Theta_S(-\log C)\cong O_S^{2}.  \] 
The weight sequence of  $C$ is also called the weight sequence of  $S$. 

\vspace{2 mm} 
\ul{Notation}  

For a sheaf  $F$  on a complex space  $X$  
we set  $h^i(X,F)=\dim H^i(X,F)$ for any integer  $i$.   

\section{Inoue surfaces} \label{inou} 

Let  $S$  be a compact connected complex surface.  
It is called a surface {\em of class VII} 
if its first Betti number $b_1=1$ and its Kodaira dimension $\kappa =-\infty$. 
It is called {\em of class VII}$_0$ if it is further minimal, i.e., 
contains no $(-1)$-curves.  
(A {\em $(-1)$-curve} is a nonsingular rational curve with self-intersection number $-1$.) 

\ul{Hopf surfaces}

$S$ is called a {\em Hopf surface} if its universal covering is isomorphic to  $\C^2-0$. 
A Hopf surface is a surface of class VII$_0$ with second Betti number $b_2=0$. 
It is called a {\em diagonal Hopf surface} 
if it is isomorphic to the quotient of $\C^2-0$ 
by an infinite cyclic group generated by a transformation of the form 
$(z,w) \ra (\al z,\beta w), \ 0< |\al|, |\beta |<1 $.   
Such a surface is diffeomorphic to the product of spheres $S^1\times S^3$.   

In this case the images of $\{z=0\}$ and $\{w=0\}$ give 
two nonsingular elliptic curves $E_1$ and  $E_2$ on  $S$. 
We use the following characterization of a diagonal Hopf surface, 
which is due to Kato-Nakamura \cite[Theorem 5.2]{na84} 
when the algebraic dimension $a(S)=0$, and 
is due to Kodaira \cite[Theorems 28, 31]{kd} when $a(S)=1$, 
(Note the difference of the definition of class VII$_0$ here and 
the one originally given by Kodaira.)  

\begin{Lemma}\label{knk} 
Let  $S$  be a surface of class VII$_0$ with infinite cyclic fundamental group.  
If  $S$  contains two smooth elliptic curves, $S$  is a diagonal Hopf surface. 
\end{Lemma} 

\ul{Inoue surfaces} 

Denote by VII$^+_0$ the class of surfaces of class VII$_0$ 
with positive second Betti number.  
The first examples of surfaces of class VII$^+_0$ 
were discovered by Inoue \cite{ino}\cite{i}, 
which we shall divide into three classes according to 
Nakamura \cite{na84} as those of hyperbolic, half and parabolic Inoue surfaces.  
(The first two surfaces are also called {\em Inoue-Hirzebruch surfaces}.) 
For the purpose of this paper it is most convenient to use 
the characterization of these surfaces similar to Lemma \ref{knk} 
due to Nakamura \cite[(8.1)(7.1)(9.2)]{na84}) as their definitions.   

\vspace{3 mm} 
{\em Definition}. 
Let $S$ be a surface of class VII$^+_0$ with $m:=b_2>0$. 
$S$ is called a {\em hyperbolic} (resp.\ {\em parabolic}) {\em Inoue surface} 
if $S$ contains two cycles of rational curves 
(resp.\ one cycle of rational curves and a nonsingular elliptic curve $E$).  
In both cases we denote by  $C$ the union of all these curves. 
$S$ is called a {\em half Inoue surface} if it contains a cycle $C$ 
of rational curves with $C^2<0$ such that 
the number of its irreducible components equals $m$. 

\vspace{3 mm} 
In all cases there are no other curves 
on the surface and $C$ is the unique maximal (reduced) curve on  $S$.  
In the hyperbolic or half Inoue case 
the number of irreducible components in $C$ 
equals  $m$, 
and in the parabolic case 
the number of irreducible components of the unique cycle equals $m$. 

All the known examples of surfaces of class VII$^+_0$ 
including Inoue surfaces 
are oritentation-preservingly diffeomorphic to  
\[M[m]:=(S^1\times S^3)\#m\bar{\B}^2 \] 
where $\bar{\B}^2$ is the complex projective plane with oritentation reversed; 
in fact they are all obtained 
as complex-analytic deformations of a blown-up Hopf surface. 
In particular they all have infinite cyclic fundamental groups. 

Starting from any of the Inoue surfaces we obtain other surfaces of 
class VII 
by blowing up successively the nodes of the cycles of rational curves 
on them. 
These surfaces again contain 
two or one cycles of rational curves, or one cycle of rational curves and 
a nonsingular elliptic curve.  
We call such surfaces  {\em properly blown-up} 
hyperbolic, half, or parabolic Inoue surfaces.  
(We include the case where the blowing down is trivial; thus 
in this terminology a hyperbolic 
Inoue surface is also a properly blown-up hyperbolic Inoue surface.) 
The reason why we consider these surfaces as well is that they arise 
equally naturally in our construction in Section 5, while 
from the viewpoint of anti-self-dual or bihermitian structures 
the minimality of the surfaces should be irrelevant.

\ul{Anti-canonical curves}  

We call a member $C$  of the anti-canonical system $|-K_S|$ of a surface  $S$  
an {\em anti-canonical curve}.  In this case we write simply $-K_S=C$ or $K_S+C=0$.  
For a diagonal Hopf surface  $S$  we have by \cite[(96)]{kd}
\begin{equation}\label{kbh} 
	-K_S = E_1 + E_2 
\end{equation}  
in the previous notation and for a properly blown-up hyperbolic or parabolic 
Inoue surface 
\begin{equation}\label{atcc} 
 -K_S  = C 
\end{equation}  
for the unique maximal curve  $C$  on it (cf.\ \cite{na84} for the minimal case;  
the general case is immedeately deduced from this case.) 
In particular any diagonal Hopf surface or 
hyperbolic or parabolic Inoue surface admits   
a disconnected anti-canonical curve.  
We may speak of {\em the} anti-canonical curve on a hyperbolic or parabolic 
Inoue  surface.  
In the half Inoue case $|-K_S|$ is empty, but $C$  becomes the unique member 
of the system $|-(K_S+L)|$, where 
$L$ is the unique non-trivial holomorphic line bundle on  $S$ with $L^2(=2L)$ trivial.  
(We call $C$ {\em $L$-twisted} anti-canonical curve.)  

In the minimal case the converse as in the following lemma holds true. 
This lemma, originally due to Nakamura, is crucial for our whole investigation. 
(See Section 10 for a proof.)  

\begin{Lemma}\label{aatc}
Let  $S$  be a compact complex surface of class VII$_0$ with infinite cyclic fundamental group. 
Suppose that there exists a disconnected anti-canonical curve on $S$. 
Then  $S$  is either a hyperbolic or parabolic Inoue surface 
or a diagonal Hopf surface. 
\end{Lemma} 

Let $S$ be as in the lemma and  $C$  an anti-canonical curve on it.  
Let  $h: \tilde{S} \ra  S$  be the blow-up of 
a finite number of points on  $C$ with exceptional curve $B$. 
Then 
by the adjunction formula  
\begin{equation}\label{ajc} 
-K_{\tilde{S}}=h^*(-K_S)-B . 
\end{equation}  
$\tilde{S}$  also admits a disconnected anti-canonical curve $\tilde{C}$ 
which is mapped surjectively onto  $C$.  
From (\ref{ajc}) we also deduce the following lemma 
which extends Lemma \ref{aatc} to non-minimal case (cf.\ \cite[Cor.3.14]{pnt6}). 

\begin{Lemma}\label{atc} 
Let  $S$  be a compact complex surface of class VII 
with infinite cyclic fundamental group. 
Suppose that there exists a disconnected anti-canonical curve  $C$ 
on $S$. Then the minimal model  $\bar{S}$  of  $S$  is 
either a hyperbolic or parabolic Inoue surface 
or a diagonal Hopf surface, and  $S \ra \bar{S}$  is obtained by blowing up $\bar{S}$ 
at a finite number of points (possibly infinitely near) 
of the image $\bar{C}$ of $C$. 
Moreover, $\bar{C}$ is an anti-canonical curve  on  $\bar{S}$. 
\end{Lemma} 

{\em Proof}. 
Let $h: S \ra \bar{S}$ be the blowing down map to the minimal model. 
We show that $\bar{C}$ is 
a disconnected anti-canonical curve on  $\bar{S}$.  We have only to prove this 
in the case of one point blown-up; the general case then follows 
by induction.  By (\ref{ajc}) we get  $|-K_{\tilde{S}}|=|-K_S(-x)|$, i.e., 
an anti-canonical curve on  $\bar{S}$ is identified with 
a anti-canonical curve on  $S$  which passes through the blown-up point $x$. 
In particular  $\bar{C}$ is an anti-canonical curve and $x \in \bar{C}$ and $C=\bar{C}$ in 
the above correspondence. 
Suppose that $\bar{C}$ is connected. 
If $x$ is a smooth point of $\bar{C}, C$ is a proper transform of $\bar{C}$  and 
is connected.  If  $x$ is a singular point of  $\bar{C}$, $B$ is in the support of 
$C$ and intersect with the proper transform of every branch passing through $x$. 
Thus  $C$  is again connected.  This shows that under our assumption 
$\bar{C}$ must be disconnected.  
Then applying Lemma \ref{atc} 
we obtain the possible structures of the surface $\bar{S}$. 
\hfill q.e.d. 

\vspace{3 mm} 
We next state an analogue of the above result in the half Inoue case. 
Let  $S$  be a compact complex surface of class VII 
with infinite cyclic fundamental group.  Then there exists a unique 
non-trivial holomorphic line bundle $L$ with $L^2=1$. 
We denote this line bundle by $L=L_S$ in what follows. 
We have the associated unramified double covering $u: \tilde{S} \ra S$ such that 
$u^*L$ is trivial.  

\begin{Lemma}\label{atch} 
Suppose that $S$ contains an $L$-twisted connected anti-canonical curve $C$ 
such that $u^{-1}(C)$ is disconnected in  $\tilde{S}$. 
Then the minimal model  $\bar{S}$  of  $S$  is a half Inoue surface 
or a diagonal Hopf surface,  
and  $S$  is obtained by blowing up $\bar{S}$ 
at a finite number of points (possibly infinitely near) on the image $\bar{C}$ of $C$. 
\end{Lemma} 

The proof is easily obtained by passing to $\tilde{S}$ 
associated to  $L$  and then applying Lemma \ref{atc}.  But we refer 
the detailed proof with more precise structures of  $S$ in this special case 
to the short note \cite{fp1}. 

\vspace{2 mm} 
\ul{Transpositions of hyperbolic Inoue surfaces}

In \cite{za} Zaffran defined for any hyperbolic Inoue surface  $S$  
its transposition  ${}^{\mbf}S$, 
which is again a hyperbolic Inoue surface with 
 ${}^{\mbf}({}^{\mbf}S) = S$.  
We shall recall its definition and basic properties. 

Let  $S$  be a hyperbolic Inoue surface.   
Let  $C^\al, \al =1,2$, be the two cycles of rational curves on  $S$.  
Then there is a geometric way of choosing one of the two (cyclic) numberings of 
the irreducible components of each $C^\al$ up to cyclic permutations due to Dloussky \cite{dlu}, 
which we shall now explain.  We shall call this a {\em canonical} numbering for $C^\al$.    

 In general, a domain  $D$  in  $S$ is called a {\em spherecal shell} if 
it is isomorphic to a domain in $\C^2$  
bounded by two concentric spheres. 
It is called {\em global} if $S-D$ is connected.  $D$ has thus two boundaries  
$\partial _+D$  and  $\partial _-D$  
which are strictly pseudoconvex and pseudoconcave respectively. 

We use the characterization of 
the canonical numbering in the form of the following lemma (cf.\ \cite{dlu}).  
Let  $B_1+\cdots +B_{h_\al}$ be a cyclic numbering of the irreducible components of $C^\al$.  
We assume that $h_\al > 1$ since otherwise the numbering in question is unique.  

\begin{Lemma}\label{d}
Suppose that there exists 
a global spherical shell  $D$  in  $S$  which 
intersects with $C^\al$ in a domain $U$ in  $B_1$. 
Among the two connected components of $B_1-U$, let 
$V$  be the component which contains $B_1\cap B_2$.  Then the numbering above 
is canonical if $\partial V=\partial D_-\cap B_1$ (instead of $\partial D_+\cap B_1$). 
\end{Lemma} 

Accordingly we may also speak of the {\em canonical} weight sequence of each $C^\al$ 
up to cyclic permutations. 
We further recall the following facts: 

a)  Let  $S$  and $S'$ be hyperbolic Inoue surfaces with two cycles of 
rational curves  $C^\al$ and $C^{'\al }$ respectively, $\al =1,2$. 
Then 
$S$ and $S'$ are isomorphic if and only if 
the canonical weight sequences of $C^\al$ and $C^{'\al }$ 
coincide up to cyclic permutations for one (and then both) of $\al =1,2$,  
after interchanging $C^{1}$ and $C^{2}$ if necessary.  
(See Remark 1.1 of \cite{za}.) 

b) For any hyperbolic Inoue surface $S$ there exists up to isomorphisms  
a unique hyperbolic Inoue surface $S'$ such that 
the canonical weight sequences 
of $C^\al$ and $C^{'\al }$ for one (and then both) of $\al =1,2$ 
are reverse to each other up to cyclic permutations,   
after interchanging $C^{1}$ and $C^{2}$ if necessary.  
$S'$ is called the {\em transposition} of  $S$  and is denote by ${}^{\mbf}S$.  
(See \cite{za}.) 

The {\em transposition} of a half Inoue surface is defined similarly, 
considering only the unique cycle instead of two cycles.  
By reducing to the corresponding minimal model 
we can also speak of the notion of transpositions of properly blown-up 
hyperbolic or half Inoue surfaces.  

\vspace{2 mm} 
\ul{Weight sequences} 

Let  $S$  be a hyperbolic Inoue surface with second Betti number $m$, and  $C^\al, \al =1,2$, 
the two cycles of rational curves on  $S$.  
The weight sequences of $C^1$ and $C^2$ are of the following form 
up to cyclic permutations and the interchange of $C^\al$:  
\begin{gather}\label{ws1} 
 (k_1+2,[k_2-1],\ldots ,k_{2n-1}+2,[k_{2n}-1]) \\ 
 ([k_1-1],k_2+2,\ldots ,[k_{2n-1}-1],k_{2n}+2) \label{ws2}  
\end{gather} 
where $n$ and $k_i, 1\leq i\leq 2n$, are positive integers 
and for a positive integer $l$, 
$[l]$ stands for the sequence $(2,\ldots, 2)\ $($l$ times), while [0] denotes 
the empty sequence \cite[(6.8)Th.]{na84}(cf.\ \cite[(2)(3)]{za}).  
However, the case $n=1$ and $k_1$ (resp.\  $k_2$) $=1$ is 
exceptional;   
in this case we should replace 
$k_2+2$ (resp.\ $k_1+2$) by $k_2$ (resp.\ $k_1$) 
\cite[(1.4)]{na84}.  Note that $m=\sum_{1\leq i\leq 2n} k_i$.  
Conversely, given $n$ and $k_i$ arbitrarily as above,  
there exists a hyperbolic Inoue surface with $b_2=m$ and 
with the above weight sequences. 
Hyperbolic Inoue surfaces are determined by 
the pair of weight sequences as above up to at most 
two non-isomorphic surfaces which are transpositions of each other. 
(The last statements holds true also for the properly blown-up case.) 
In particular there are only a countable number of hyperbolic surfaces 
up to isomorphisms.  

\vspace{2 mm} 

\ul{Isomorphism classes of parabolic Inoue surfaces}

For a parabolic Inoue surface $S$  with second Betti number $m$,  
the weight sequence of its unique cycle is given by $[m]$ for $m>1$ and $0$ for $m=1$ 
in the above notation, while 
the elliptic curve  $E$  on it has the self-intersection number  $E^2=-m$. 

Parabolic Inoue surfaces with fixed second Betti number 
are parametrized by the punctured unit disc 
$D^*=\{|d|<1\}$ (cf.\ \cite[(1.1)]{na84}).  
So we may write $S= S_d$  for some $d\in D^*$.  
The parameter $d$ is geometrically interpreted as follows.  
Let  $u: U \ra  S$  be the universal covering of  $S$.  
Then for the unique elliptic curve  $E$  in  $S$, 
we get an infinite cyclic covering  $v: \tilde{E}:= u^{-1}(E) \ra  E$. 
Let  $\gamma $  be a fixed 
generator of the covering transformation group.  
Then there exists a unique 
complex number $\al$ with $0<\al <1$ such that 
with respect to an isomorphism  $w: \tilde{E} \stackrel{\sim}{\ra}  \C^*=\C^*(s)$,  
$\gamma $ takes the form $\gamma (s)=\al s$.  
This number $\al$  is independent of the choice of $\gamma $ and 
the isomorphism  $w$  and depends only on the isomorphism class of  $v$.  
In fact, if $S=S_d$ the construction \cite{na84} clearly shows that $d=\al $. 
In particular we get

\begin{Lemma}\label{prbi}
Let  $S$  be a parabolic Inoue surface with fixed second Betti number $m$, 
and  $v: \tilde{E} \ra  E$  as above.  
Then the isomorphism class of  $S$  is determined by the isomorphism class 
of the infinite cyclic covering $v$.  $D^*$ is thus the moduli space of 
parabolic Inoue surface. 
\end{Lemma} 

\ul{Real structure on Inoue surfaces}  

Let  $J$  be a complex structure on a smooth manifold  $M$.  
Let  $S:=(M,J)$ be the resulting complex manifold and 
$\bar{S}:=(M,-J)$ its complex conjugate.  Then $S$ and $\bar{S}$ are 
biholomorphic if and only if  $S$  admits an anti-holomorphic diffeomorphism; 
in particular a real structure, i.e., an anti-holomorphic 
involution.  In this context we note the following: 
\begin{Lemma}\label{rsi}
Any hyperbolic, half, or parabolic Inoue surface $S$ 
has a natural real structure.  The same is true for a proper blowing-up of 
any such surface.  
\end{Lemma} 
{\em Proof}.  
The universal covering  $U$  of  $S$  
is covered by coordinate neighborhoods, such that in the intersection of 
any two of them the two coordinates are related 
by Laurent {\em monomials} (cf.\ \cite[\S 1]{na84}).  
Hence the complex conjugations with respect to 
each such coordinates are compatible in the intersections and 
give a real structure  $\tilde{\mu}$  on  $U$.  
Moreover, a generating covering transformation 
of $U \ra S$ is also given by Laurent monomials (cf.\ \cite[\S 1]{na84}) 
and hence the real structure $\tilde{\mu}$ descends 
to a real structure $\mu$ on  $S$.  
Moreover, since the nodes of the anti-canonical divisors 
of these surfaces are fixed points of the real structure, $\mu$ lifts to its 
proper blowing-ups. 
\hfill q.e.d. 
 
\vspace{3 mm} 

\ul{Deformations of Inoue surfaces}  

Let  $S$  be a properly blown-up hyperbolic, half or parabolic Inoue surface. 
For brevity we refer to the hyperbolic (resp.\ half, resp.\ parabolic) case 
as {\em Case-H} (resp.\ {\em Case-H$'$}, resp.\ {\em Case-P}) in what follows. 
Let $C$ be the unique maximal curve on  $S$.  
In other words, $C$  is the unique anti-canonical curve in Case H and -P, while 
it is the unique $L$-twisted anti-canonical curve in Case H$'$ where $L=L_S$. 

Let $n: \tilde{C}:=\coprod _{1\leq d\leq b} \tilde{C}_d \ra C$ be 
the normalization of  $C$, where $\tilde{C}_d$ are the normalizations of 
the irreducible components of  $C$, and  $b=m$ in Case-H or -H$'$ and $=m+1$ in Case-P. 
It is known that  Aut$_0(S,C) = \{e\}$ in Case-H or H$'$ and $\C^*$ in Case-P (cf.\ \cite[Prop.2.5]{dl0} and \cite{pik} for Case-H and -H$'$ and \cite{hau} for Case-P). 
\begin{Lemma}\label{vai}
In Case-H or -H$'$ we have $h^i(\Theta_{S}(-\log C))=0$ for $i=0,1,2$, 
while in Case-P 
we have   
$h^i(\Theta_{S}(-\log C))=1$ for $i=0,1$ and $=0$ for $i=2$. 
\end{Lemma} 

{\em Proof}.  
First we consider Case-H and -P. 
We have the following exact sequence 
\begin{equation}\label{sho} 
\begin{CD}
0	@>>> \Theta_{S}(-C) @>>> \Theta_{S}(-\log C) 
@>>> \Theta_{C}	 @>>> 0 
\end{CD}
\end{equation}  
and the isomorphisms  
\[\Theta_{C}\cong \oplus_d n_*\Theta_{\tilde{C}_d}(-(0_d+\infty_d))) \cong 
\oplus_d n_*O_{C_d} \] 
where $0_d$ and $\infty_d$ are 
the inverse images of the nodes of $C$ in $C_d$.  
Since $C=-K_{S}$ by (\ref{atcc}), we have  $\Theta_{S}(-C)\cong \Omega_S$ so that 
$h^i(\Theta_{S}(-C))=h^{i}(\Omega_{S})=h^{1,i}$, where $h^{p,q}$ denote 
the Hodge numbers. Hence we get 
$h^i(\Theta_{S}(-C))=0$ for $i=0,2$ and $h^1(\Theta_{S}(-C))=h^{1,1}=b_2(S)=m$.  
Thus taking the long exact sequence associated to (\ref{sho}) we obtain 
$h^2(\Theta_{S}(-\log C))=0$ (cf.\ \cite[Th.1.3]{nato}) and the exact sequence 
\begin{equation}\label{hpb} 
\begin{CD}
 0
 @>>> H^0(\Theta_{S}(-\log C))   
 @>>>  \C^b  
@>>> \C^m 
 @>>> H^1(\Theta_{S}(-\log C))   \\ 
 @>>> \C^k    @>>> 0 
\end{CD}
\end{equation}  
where  $k=0$ in Case-H and $=1$ in Case-P.   
On the other hand, as we noted before the lemma 
 $h^0(\Theta_{S}(-\log C))=0$ in Case-H  and $=1$ in Case-P.  
Thus the lemma is proved in these cases.   

In Case-H$'$ take the canonical finite unramified double covering 
$(\tilde{S},\tilde{C})\ra (S,C)$  with covering involution $\iota $.  
Then $H^i(\Theta_{S}(-C))$ are naturally identified 
with the subspaces  $H^i(\Theta_{\tilde{S}}(-\tilde{C}))^\iota $ 
of $\iota $-fixed elements for all $i$.  
From this the results follow from those in Case-H. 
\hfill q.e.d.  

\vspace{3 mm} 
As follows from the above proof in Case-P the restriction map  $\Theta_{S}(-\log C) \ra \Theta_{E}$ induces 
a natural isomorphism 
\begin{equation}\label{eli} 
 H^1(\Theta_S(-\log C))\cong H^1(\Theta_{E}) 
\end{equation}  
where  $E$  is the elliptic component  of $C$.  

The two sheaves 
$\Omega'_{S}(\log C)$ (cf.\ (\ref{om'})) and $\Omega_{S}(\log C)$ 
coincide except at the nodes of  $C$.  Let $B$ be the set of nodes of $C$. 
Then more precisely we have the following: 
\begin{Lemma}\label{wlo} 
We have an exact sequence 
\begin{equation}\label{oO1} 
0 \ra  \Omega'_{S}(\log C) \ra  
\Omega_{S}(\log C) \ra  \oplus_{p\in B} \C_p \ra  0, 
\end{equation}  
where $\C_p $ is the skyscraper sheaf with support $p$ and with fiber $\C$. 
\end{Lemma} 
{\em Proof}.  
The fact that the quotient 
$\Omega_{S}(\log C)/\Omega'_{S}(\log C) $ is 
isomorphic to $\C_p$ at each $p\in B$ is seen by checking the image of 
$\Omega'_{S}(\log C)$ by the Poincare residue map 
$P:  \Omega_{S}(\log C) \ra  O_1\oplus O_2 $, 
where  $O_s, s=1,2$, are the structure sheaves of the two irreducible components 
of  $C$ at  $p$.  
In the local model  $(\C^2(x,y), xy=0)$ of $(S,C)$, $P$ takes the form 
$ a(x,y)dx/x + b(x,y) dy/y  \ra  (a(x,0), b(0,y)) \in O_1\oplus O_2$, 
while the elements $\Omega'_{S}(\log C) $
are of the form  $f(x,y)(dx/x+dy/y)$ so that their images are 
given by  $(f(x,0),f(0,y))$.  The assertion thus holds. \hfill q.e.d.  

\vspace{3 mm} 
Since  
${\mathcal Ext}^{i}_{O_{S}}(\C_p,O_{S}) = 0$ for $i\neq 2$ and 
$= \C$ for  $i=2$, by applying 
${\mathcal Ext}^{i}_{O_{S}}(-,O_{S})$  to the sequence (\ref{oO1}) 
we have the following:  
\begin{Corollary}\label{exiso1}
1) ${\mathcal Hom}_{O_{S}}(\Omega'_{S}(\log C),O_{S})\cong 
\Theta_{S}(-\log C)$. 

2) There is a natural isomorphism 
\begin{equation}\label{ilm1} 
{\mathcal Ext}^{1}_{O_{S}}(\Omega'_{S}(\log C),O_{S}) 
\cong  \oplus_{p\in B} \C_p . 
\end{equation}  
	
3) 
${\mathcal Ext}^{i}_{O_{S}}(\Omega'_{S}(\log C),O_{S})=0$  for   $i\geq 2$.
\end{Corollary} 

In view of Corollary \ref{exiso1} and Lemma \ref{vai} 
the local to global spectral sequence 
for  $Ext^{i}_{O_{S}}(\Omega'_{S}(\log C),O_{S})$   
yields the following:  
\begin{Lemma}\label{vext}
We have 
\begin{gather}\label{vextt}
  Ext^{2}_{O_S}(\Omega'_{S}(\log C),O_{S})=0, \\ 
 0 \ra H^{1}(\Theta_{S}(-\log C))	
 \ra Ext^{1}_{O_S}(\Omega'_{S}(\log C),O_{S})  	\label{u}
  \ra \oplus_{p\in B} \C_p \ra 0, \\
  Ext^{0}_{O_S}(\Omega'_{S}(\log C),O_{S})=H^0(\Theta_{S}(-\log C)), 
\end{gather} 
where the sequence {\em (}\ref{u}{\em )} is exact. 
\end{Lemma} 
Let 
\begin{equation}\label{kr1} 
 g: ({\mathcal S},{\mathcal C}) \ra T,\ (S_o,C_o)=(S,C),\ o\in T 
\end{equation}  
be the Kuranishi family of deformations of the pair $(S,C)$. 
Also for any  $p\in B$ we denote by  
\[ g_p : {\mathcal C}(p) \ra T_p,\ (C_o,o)\cong (C,p),\ o\in T_p \]  
the Kuranishi family of deformations of the isolated singularity $(C,p)$, 
where $T_p$ is smooth of dimension one.  
Any deformation of 
$(S,C)$  induces a deformation of $(C,p)$ 
and correspondingly we have a versal map  $\tau_p: T \ra  T_p$. 
The fiber  $T(p):= \tau_p^{-1}(o)$ is uniquely determined independently of 
the choice of $\tau_p$.  Thus a point  $t\in T$  is outside of  $T(p)$ precisely when the 
two irreducible components of $C$  passing through $p$  are merged together 
in $C_t$ to become one smooth curve locally.  

In Case-P we also consider 
the Kuranishi family 
\[  e: {\mathcal E} \ra T_E,\ E_o=E,\ o\in T_E \]  
of the elliptic curve $E$ in  $S$, 
where $T_E$ is smooth of dimension one. 
Since a deformation of $(S,C)$ induces a deformation of $E$, 
we have the (unique) versal map $\tau_E: T \ra T_E$.  

  In Case-H and -P 
we write  $C$  as the disjoint union $C=C^{1}\cup C^{2}$ of two curves, where 
$C^{\al } , \al = 1,2$, are cycles of rational curves in Case-H, and 
$C^{1}$ is a cycle of rational curves and $C^{2}=E$ in Case-P.  
The fiber  $C_t, t\in T$, is similarly a disjoint union 
$C_t=C^{1}_{t}\cup C^{2}_{t}$, 
where  $C^{\al }_{t}$ is either a cycle of rational curves or a smooth elliptic curve. 
Recall that $\#B = m$.  

\begin{Proposition}\label{wa} 
Let  $(S,C)$  be as above. 
Then the Kuranishi space $T$ is smooth of dimension $m$ {\em (}resp.\ $m+1${\em )} 
in Case-H or -H$'$ {\em (}resp.\ -P{\em )}. 
Moreover, in each case we have the following: 

a{\em )} Case-H or -H$'$:  The product map 
\[ \Pi_p \tau _p : T\ra  \Pi_{p\in B} T_p\]  
is isomorphic; in particular $T(p)$ is a smooth hypersuface in  $T$ for each $p\in B$.  
The family is universal at each point of $T$.  
Accordingly,   $\dim Ext^{1}_{O_{S_t}}(\Omega'_{S_t}(\log C_t),O_{S_t})=m$ 
independently of  $t \in T$. 

b) Case-P: 
The product map 
\[ \Pi_p \tau _p \times  \tau _E: T\ra  \Pi_p T_p \times  T_E\]  
is isomorphic. The family is not universal.  
In fact  $h^0(\Theta_{S_t}(-\log C_t)))=1$ {\em (}resp.\ $=0${\em )} and  
$\dim Ext^{1}_{O_{S_t}}(\Omega'_{S_t}(\log C_t),O_{S_t})=m+1$ {\em (}resp.\ $=m${\em )} for 
for  $t \in I$ {\em (}resp.\ $\notin I${\em )}, 
where $I:= (\Pi_p \tau_p)^{-1}(o), o\in T_E$, is a submanifold of dimension one. 
\end{Proposition} 

{\em Proof}.  
The smoothness of $T$ follows from (\ref{vextt}). 
We have $\dim T = \dim Ext^{1}_{O_S}(\Omega'_{S}(\log C),O_{S})$ and 
the latter is identified with the claimed value by Lemmas \ref{vai} and \ref{vext}. 
By Lemma \ref{vai} and the upper semicontinuity of cohomology, we have 
$h^0(\Theta_{S_t}(-\log C_t)) =0$  independently of  $t$ in Case-H or -H$'$. 
Hence this family is universal at each point of $T$. 
The third arrow of (\ref{u}) is identified 
with the differential of $\Pi_p \tau _p$ 
at the base point.  
The first assertion of a) follows from this. 

In Case-P  $\Pi_p \tau _p$ is a submersion and the inverse image 
$I$ of the reference point 
is identified with the local moduli space of  $S$  as a parabolic Inoue surface,  
whose tangent space is identified with $H^1(\Theta_{S}(-\log C))$.  
The differential of the restriction of  $\tau_E$ to $I$ is identified with 
the isomorphism (\ref{eli}).  From this we get the first assertion in b). 
  The rest follows from the fact that Aut$_0(S_t,C_t)=\{e\}$ for 
any $t \notin I$. \hfill q.e.d.  

\vspace{3 mm} 
In the next Proposition \ref{seis} we assume that $S$ is 
a properly blown-up hyperbolic or parabolic Inoue surface.  
We write the set $B$ of nodes as the disjoint union $B=B_1\cup  B_2$ in the obvious way, 
where $B_2=\emptyset$ in Case-P. 
Then on the structure of the surfaces  $S_t$ in the family 
we have the following proposition. 

\begin{Proposition}\label{seis}
Let $S$ be as above. 

1) For any $t\in T$,  $C_t$ is the unique anti-canonical curve on $S_t$. 

2) 
If $t\notin T(p)$ for some $p\in B$ (i.e., $t\neq o$ in Case-H), 
$S_t$ is not minimal and its minimal model $\bar{S}_t$ is either 
a hyperbolic or parabolic Inoue surface or a diagonal Hopf surface. 

3) $\bar{S}_t$ is a diagonal Hopf surface 
if and only if  $t\notin T(p)$ for any $p\in B$. 
$\bar{S}_t$ is a a parabolic Inoue surface if and only if 
$t\notin T(p)$ for any $p\in B_\al$ for one of $\al$ ($=1$ or $2$), 
but not for both in Case-H 
(resp.\ $t\in T(p)$ for some  $p\in B_1$ in Case-P). 
\end{Proposition} 

{\em Remark 3.1}. 
We may call $(S,C)$ an {\em anti-canonical pair} in the sense that 
$C$ is an anti-canonical curve on  $S$. 
The above lemma implies that 
the Kuranishi family (\ref{kr1}) of $(S,C)$  is actually 
a Kuranishi family of $(S,C)$ as an anti-canonical pair. 
Thus we can identify our Kuranishi family with 
the family constructed by Nakamura in Lemma 5.7 of \cite{na84}. 
(Indeed, we can show that the family (\ref{kr1}) is realized as a subfamily 
of the Kuranishi family of  $S$ itself.) 
However, in \cite{na84} neither the smoothness of  $T$ nor the precise structures 
of $T$ as above is clear. 
\vspace{3 mm} 

{\em Proof}.  
1)  
Consider the short exact sequence 
\[
\begin{CD}
0	@>>> O_{S_t}(-C_t) @>>> O_{S_t} @>>> O_{C_t}	 @>>> 0 
\end{CD}
\] 
and the associated long exact sequence 
\[
\begin{CD}\label{sotsp}
 @>>> H^1(O_{S_t})  
 @>>>H^ 1(O_{C_t})  
 @>>>	H^2(O_{S_t}(-C_t))   @>>> 
\end{CD}
\] 
Together with Serre duality and the upper semicontinuoity of cohomology 
this yields 
\[ 2=h^ 1(O_{C_t})  \leq h^1(O_{S_t})+h^0(K_t+C_t) \leq 1+h^0(K+C)=2.\] 
Since $h^1(O_{S_t})=1, $
we get that 
$h^0(K_t+C_t)=1$ for all $t$. 
Then any non-vanishing element $u_0$ of $H^0(K+C)$ extends locally to 
elements $u_t$ of $H^0(K_t+C_t)$, which is again non-vanishing since so is 
$u_0$.  Thus $K_t+C_t = 0$ as desired.  The uniqueness then follows from the 
inequality  $h^0(-K_t)\leq h^0(-K)=1$. 

2) is then a consequence of Lemma \ref{atc}. 

3)  
$C^\al_{t}, \al =1,2$,  is a smooth elliptic curve 
if and only if all the nodes of $C^\al$ are smoothed 
in the deformation $C^\al_{t}$, and this is precisely the condition that 
$t\notin T(p)$ for all $p\in B_\al$. 
From this the conclusion follows from Lemma \ref{knk}.  \hfill q.e.d.  

\vspace{3 mm} 
When  $S$  is a properly blown-up half Inoue surface, 
the statement analogous to Proposition \ref{seis} is given as follows: 

\begin{Proposition}\label{seis'}
Suppose that $S$ is a properly blown-up half Inoue surface.  Then: 

1) For any $t\in T$,  $C_t$ is the unique $L_t$-twsited anti-canonical curve 
on $S_t$, where  $L_t=L_{S_t}$. 

2) 
If $t\neq o$, 
$S_t$ is not minimal and its minimal model $\bar{S}_t$ is 
either a half Inoue surface or a diagonal Hopf surface. 

3) $\bar{S}_t$ is a diagonal Hopf surface 
if and only if  $t\notin T(p)$ for any $p\in B$. 
\end{Proposition} 

As in the case of Lemma \ref{atch} we refer the proof of this 
proposition to the short note \cite{fp1}. 

\vspace{5 mm} 
\section{Deformations of rational surface with a nodal curve} \label{rswn} 

Let  $\tilde{S}$  be a projective toric surface acted by $G:=\C^{*2}$,  
and  $\tilde{C}=\tilde{C}_1+\cdots +\tilde{C}_{k+2}$ 
the anti-canonical cycle on  $\tilde{S}$,  
where we assume for simplicity that $k>2$. 
We put $\infty_i=0_{i+1}=\tilde{C}_i\cap \tilde{C}_{i+1}$, 
where the subscripts are considered cyclically modulo $k+2$. 

Suppose that 
there exist disjoint irreducible components $H$ and $E$ of $\tilde{C}$ 
with $H^2=1$ and $E^2=-1$ respectively. (In this case we call the toric surface $(S,C)$ 
{\em admissible}.) 
Take an isomorphism 
\begin{equation}\label{idf} 
\varphi : (H,0_H,\infty_H) \ra (E,\infty_E,0_E) \ \ \mbox{or} \ \ 
(H,0_H,\infty_H) \ra (E,0_E, \infty_E)
\end{equation}  
where $0_H=0_i$ if $H=\tilde{C}_i$ etc.  
In the latter case we call $\varphi$ {\em of twisted type} 
and in the former case {\em of untwisted type}. 
Let  $\hat{S}$ be the non-normal surface 
obtained by identifying the points $x\in H$ with $\varphi (x)\in E$.  

Let $n: \tilde{S} \ra \hat{S}$ be the natural map and 
denote the singular locus of  $\hat{S}$ by $\hat{F}=n(H)=n(E)\cong \B$. 
Let  $\tilde{C}^\al , \al =1,2$, be 
the connected components of the union of irreducible components of 
$\tilde{C}$ other than $H$ and  $E$.  Denote by $\hat{C}^\al$ 
their images in $\hat{S}$, and  put $\hat{C}=\hat{C}^1\cup \hat{C}^2$.  
When $\varphi$ is of untwisted type,  
$\hat{C}^\al$ are disjoint and each forms a cycle of rational curves on $\hat{S}$, 
while when $\varphi$ is of twisted type, 
$\hat{C}$ is connected and 
forms a single cycle of rational curves on  $\hat{S}$. 
In both cases $\hat{C}$ is an anti-canonical curve on $\hat{S}$ 
(cf.\ Lemma \ref{lll} below).  

In this section we study the smoothing of the singular surface  $\hat{S}$ under 
deformations.  This subject was studied extensively by Nakamura in \cite{nkm1}, 
\cite{na84}, \cite{nato}.  However, since we treat it from a little different point of view, 
we will describe some details here.

\vspace{3 mm} 
\ul{Automorphism groups}

\vspace{3 mm} 
We compute the identity component of the automorphism group of $(\hat{S},\hat{C})$. 
\begin{Lemma}\label{aut} We have 
\[ \rm{Aut}_0(\hat{S},\hat{C}) \cong \C^*\] 
\end{Lemma} 
{\em Proof}. 
Fix an affine coordinate  $z_i$ on each  $\tilde{C}_i$ 
with  $0_i$ and $\infty_i$ corresponding to $0$ and $\infty$.  
With respect to this coordinate the action of $G$ on  $\tilde{C}_i$ is written as 
$z_i \ra \chi _i(g)z_i, g\in G$, for a unique character  $\chi _i$ of $G$, 
which is independent of the choice of $z_i$. 

In view of the natural inclusion of algebraic groups 
Aut$_0(\hat{S},\hat{C})\hookrightarrow$ Aut$_0(\tilde{S},\tilde{C})$,  
it suffices to show 
that there exists a unique one dimensional subgroup 
of Aut$_0(\tilde{S},\tilde{C})\cong G$ which descends to an automorphism group 
of $(\hat{S},\hat{C})$.  
A one-parameter subgroup  $\rho $ of $G$  induces a 
$\C^*$-action on  $\hat{S}$  
if and only if  $\varphi$ is $\rho$-equivariant with respect to 
the induced $\rho $-actions on  $H$ and $E$.  
When $\varphi$ is of untwisted (resp.\ twisted) type 
$\varphi$ is written as $z_E = \varphi (z_H) = a/z_H$ (resp.\ $= az_H$) 
for some $a\neq 0$, where $z_H=z_i$ etc.\ as before.  
Thus the condition becomes 
\begin{equation}\label{roh} 
\chi _H\rho =-\chi _E\rho \ \ \mbox{(resp.\ } \ \chi_H\rho = \chi _E\rho ) 
\end{equation}  
if $\varphi$ is of untwisted (resp.\ twisted) type, where $\chi _H=\chi_i$ if 
$H=C_i$ etc. 
Under our assumption that $k>2$, we easily see that $\chi _H \neq \pm\chi _E$, 
and hence (\ref{roh}) defines a unique one-parameter 
subgroup as desired. \hfill q.e.d.  

\vspace{3 mm} 
It is convenient to formulate the above result in a more formal way as follows.  
Let $ M:= \Z^2$ be the free abelian group consisting of 
all the characters $\chi : G \ra  \C^*$ of 
$G$ and $N:= \Z^2$ 
the free abelian group consisting of all the one-parameter subgroups  
$\rho :\C^* \ra  G$ of $G$, written additively. 
We have a natural perfect pairing  $\langle , \rangle : M\times N \ra  \Z$, 
where $\langle \rho ,\chi \rangle =l$ if $\chi \rho(t) =t^l, l\in \Z, t\in \C^*$. 
We can define an orientation of  $M$ by the condition that $-\chi_{i-1},\chi_i$  
form an oriented basis of $M$  for any $i$. 

For $\chi\in M$ 
let  $\chi ^{\perp} $ be the unique element of  $N$ such that 
it is orthogonal to $\chi $, has the same length as $\chi $ and 
that $\chi^{\perp} $ and  $\chi $ define the positive oritentation on $\Z^2$. 
We have $(\chi +\chi ')^{\perp} =\chi ^{\perp} + \chi^{'\perp} $. 
On the other hand, 
$\langle \chi_i, \rho_i\rangle =0$,  
$\langle \chi_{i-1}^{\perp},\rho_i\rangle > 0$ and hence 
$\chi_i^{\perp}=\rho_i$. 
Thus we obtain the following supplement to Lemma \ref{aut} giving the 
explicit description of the one-parameter group in question. 
\begin{Lemma}\label{ps} 
Aut$_0(\hat{S},\hat{C})\cong \C^*$, being induced by 
the one-parameter subgroup $\rho_H +\rho_E$ {\em (}resp.\ $\rho_H - \rho_E${\em )} 
of $G$ if $\varphi$ is of untwisted {\em (}resp.\ twisted{\em )} type. 
\end{Lemma} 

\vspace{3 mm} 

\newpage 
\ul{Computation of Ext groups} 

Since $\hat{C}$ is a curve with normal crossings, the following is well-known: 
\begin{gather}
Ext^{1}_{O_{\hat{C}}}(\Omega_{\hat{C}},O_{\hat{C}}) \cong \label{exc}
\oplus_{p\in \bar{B}} \C_p \\ 
  H^1(\Theta_{\hat{C}}) = H^2(\Theta_{\hat{C}}) = 
  Ext^{2}_{O_{\hat{C}}}(\Omega_{\hat{C}},O_{\hat{C}})=0 , 
\end{gather} 
where  $\bar{B}$ is the set of nodes of $\hat{C}$ with $\#\bar{B}=k$ and $\C_p$ 
is the skyscraper sheaf at $p$ with fiber $\C$.  
Similarly we know that 
\begin{gather}\label{ex1}
{\mathcal Ext}^{i}_{O_{\hat{S}}}(\Omega_{\hat{S}},O_{\hat{S}}) =0, \ i\geq 2 \\ 
{\mathcal Ext}^{1}_{O_{\hat{S}}}(\Omega_{\hat{S}},O_{\hat{S}}) \cong  O_{\hat{F}} 
\label{ex2}
\end{gather} 
(cf.\ \cite{fr}) and hence the local to global spectral sequence yields the exact sequence: 
\begin{equation}\label{lctop} 
0 \ra    H^{1}(\Theta_{\hat{S}})
\ra Ext^{1}_{O_{\hat{S}}}(\Omega_{\hat{S}},O_{\hat{S}})  \ra 
H^0(O_{\hat{F}})\ra 0.  
\end{equation}  

We shall next show the following: 
\begin{Lemma}\label{tev}
\begin{equation}\label{dame} 
h^{0}(\Theta_{\hat{S}}(-\log\hat{C}))=1  \ \mbox{and} \ \ 
h^{q}(\Theta_{\hat{S}}(-\log\hat{C}))=0, \ q\geq 1. 
\end{equation}  
\end{Lemma} 

{\em Proof}.  
The first one follows from Lemma \ref{aut}. 
Let  $\hat{D}$ be the set of the 
two intersection points, say $r_\al$, $\al = 1,2$, 
of $\hat{F}$ and $\hat{C}$. 
Take the normalization exact sequence for  $n: \tilde{S}\ra \hat{S}$: 
\begin{equation}\label{sott} 
\begin{CD}
0 @>>> \Theta_{\hat{S}}(-\log\hat{C}) 
 @>>> n_*(\Theta_{\tilde{{S}}}(-\log\tilde{C}))  
 @>>>\Theta_{\hat{F}}(-\log \hat{D}) 
 @>>>	0 . 
\end{CD}
\end{equation}  
Since  $\Theta_{\tilde{S}}(-\log\tilde{C})\cong O_{\tilde{{S}}}^{\oplus 2}$ 
and $\Theta_{\hat{F}}(-\log \hat{D}) \cong O_{\hat{F}}$,  
from the long exact sequence associated to (\ref{sott}) 
we get the vanishing of $H^{2}(\Theta_{\hat{S}}(-\log\hat{C}))$ and 
the exact sequence 
\[
\begin{CD}\label{sots}
0 @>>> \C  @>>> \C^2  @>>> \C 
 @>>>	H^1(\Theta_{\hat{S}}(-\log\hat{C}))   @>>> 0 
\end{CD}
\] 
The lemma follows from this. 
\hfill q.e.d.  

\vspace{3 mm} 
Next we prove:  

\begin{Lemma}\label{syou}
${\mathcal Hom}_{O_{\hat{S}}}(\Omega_{\hat{S}}(\log \hat{C}),O_{\hat{S}})\cong 
\Theta_{\hat{S}}(-\log\hat{C})$ and   
${\mathcal Ext}^{i}_{O_{\hat{S}}}(\Omega_{\hat{S}}(\log \hat{C}),O_{\hat{S}})
\cong 
{\mathcal Ext}^{i}_{O_{\hat{S}}}(\Omega_{\hat{S}},O_{\hat{S}})$ for $i\geq 1$.  
\end{Lemma} 

{\em Proof}.  
The first isomorphism is well-known (cf.\ Proposition \ref{lg} below). 
For the second isomorphism 
we observe the sheaf exact sequence 
\begin{equation}\label{po} 
0\ra \Omega_{\hat{S}} \ra \Omega_{\hat{S}}(\log \hat{C}) 
\stackrel{P}{\ra} \oplus_l O_{C'_l} \ra 0 
\end{equation}  
where the last direct sum is over Cartier irreducible components $C'_l$ of 
$\hat{C}$, and $P$ is the Poincare residue map. 
Here irreducible components $C_i$ of $\hat{C}$ with $C_i\cap \hat{F}=\emptyset$ 
are Cartier irreducible components and 
the unions of two irreducible components passing through 
$r_\al, \al=1,2$, are the remaining Cartier irreducible components. 
Since ${\mathcal Ext}^{i}_{O_{\hat{S}}}(O_{\hat{C'_l}},O_{\hat{S}})\cong N'_l$ 
for $i=1$ 
and vanish otherwise, 
applying ${\mathcal Ext}^{i}_{O_{\hat{S}}}(-,O_{\hat{S}})$ to the above sequence 
we obtain the exact sequence of sheaves 
\begin{equation}\label{oto} 
\begin{CD}
0 @>>> \Theta_{\hat{S}}(-\log\hat{C}) 
 @>>> \Theta_{\hat{S}}
 @>>> \oplus_l N'_l \\ 
 @>>> {\mathcal Ext}^{1}_{O_{\hat{S}}}(\Omega_{\hat{S}}(\log \hat{C}),O_{\hat{S}})
  @>>>	{\mathcal Ext}^{1}_{O_{\hat{S}}}(\Omega_{\hat{S}},O_{\hat{S}}) 
  @>>>	0 
\end{CD}
\end{equation}  
where $N'_l$ is the normal bundle of $C'_l$ in  $\hat{S}_l$.  
The lemma follows easily from this.   \hfill q.e.d. 
 
\vspace{3 mm} 
Together with 
(\ref{lctog}), (\ref{ex1}), (\ref{ex2}) and Lemma \ref{tev} 
the lemma implies the following: 
\begin{Lemma}\label{exp} We have 
\begin{gather}
Ext^{2}_{O_{\hat{S}}}(\Omega_{\hat{S}}(\log \hat{C}),O_{\hat{S}})=0 
\ \ \mbox{and} \ \ \label{elo1} \\ 
Ext^{1}_{O_{\hat{S}}}(\Omega_{\hat{S}}(\log \hat{C}),O_{\hat{S}})\cong 
H^0(O_{\hat{F}})\cong \C. \label{elo2}
\end{gather} 
\end{Lemma} 
Let  $\hat{B}:= \bar{B}-\hat{D}$ be 
the set of nodes of $\hat{C}$ outside $\hat{F}$ 
with  $\#\hat{B} = m:=k-2$. 
$\bar{B}-\hat{B}$ consists 
of the two  points  $r_\al$, $\al = 1,2$, as above. 

The two sheaves 
$\Omega'_{\hat{S}}(\log \hat{C})$ (cf.\ (\ref{om'})) and $\Omega_{\hat{S}}(\log \hat{C})$ 
coincide except at points in $\hat{B}$. Thus by Lemma \ref{wlo} we get 
a natural exact sequence: 
\begin{equation}\label{oO} 
0 \ra  
\Omega'_{\hat{S}}(\log \hat{C}) \ra  
\Omega_{\hat{S}}(\log \hat{C}) \ra  
\oplus_{p\in \hat{B}} \C_p \ra  0.  
\end{equation}

\vspace{3 mm} 
Then similarly to Corollary \ref{exiso1} we get 
\begin{Lemma}\label{exiso}
1) ${\mathcal Hom}_{O_{\hat{S}}}(\Omega'_{\hat{S}}(\log \hat{C}),O_{\hat{S}})\cong 
\Theta_{\hat{S}}(-\log \hat{C})$. 

2) There is a natural exact sequence 
\begin{equation}\label{ilm} 
 0 \ra O_{\hat{F}} \ra 
{\mathcal Ext}^{1}_{O_{\hat{S}}}(\Omega'_{\hat{S}}(\log \hat{C}),O_{\hat{S}}) 
\ra \oplus_{p\in \hat{B}} \C_p\ra 0. 
\end{equation}  
	
3) 
${\mathcal Ext}^{i}_{O_{\hat{S}}}(\Omega'_{\hat{S}}(\log \hat{C}),
O_{\hat{S}})=0$  for   $i\geq 2$.
\end{Lemma} 

In view of Lemmas \ref{tev} and \ref{exiso} 
the local to global spectral sequence 
for $Ext^{i}_{O_{\hat{S}}}(\Omega'_{\hat{S}}(\log \hat{C}),O_{\hat{S}})$   
yields the first assertion of the following:  
\begin{Proposition}\label{diff}
1) We have 
\begin{equation}\label{vn} 
 Ext^{2}_{O_{\hat{S}}}(\Omega'_{\hat{S}}(\log \hat{C}),O_{\hat{S}})=0. 
\end{equation}  

2) There exists a natural isomorphism: 
\begin{equation}\label{lctopp} 
c: Ext^{1}_{O_{\hat{S}}}(\Omega'_{\hat{S}}(\log \hat{C}),O_{\hat{S}}) 
\stackrel{\sim}{\ra} 
\C_{r} \oplus(\oplus_{p\in \hat{B}} \C_p). 
\end{equation}  
where $r$  is any one of  $r_\al, \al = 1, 2$. 
In particular $\dim Ext^{1}_{O_{\hat{S}}}(\Omega'_{\hat{S}}(\log \hat{C}),O_{\hat{S}})
= m+1 $. 
\end{Proposition} 

{\em Proof}.  We shall show 2). 
We first prove that the following sequence is exact: 
\begin{equation} \label{ooss} 
0 \ra \Omega'_{\hat{S}}(\log \hat{C}) \stackrel{a}{\ra}  
\Omega_{\hat{S}}(\hat{C}) \stackrel{b}{\ra} 
\Omega_{\hat{C}}\otimes \hat{N} \ra 0 , 
\end{equation}  
where $\hat{N}$ is the normal bundle of $\hat{C}$ in  $\hat{S}$. 
Indeed,  the restriction map $b$ is given locally 
by  $[dx/xy, dy/xy \ra dx|\hat{C}, dy|\hat{C}]$, 
where $|\hat{C}$ denotes the restriction to  $\hat{C}$.  
Note that the element $dx = -dy$, which generates  
$\Omega_{\hat{S}}(\log \hat{C})/\Omega'_{\hat{S}}(\log \hat{C})$, 
generates locally the torsion part 
$\tau$ of $\Omega_{\hat{C}}\cong \Omega_{\hat{C}}\otimes \hat{N} $, 
thus inducing the isomorphism 
$\Omega_{\hat{S}}(\log \hat{C})/\Omega'_{\hat{S}}(\log \hat{C})\cong \tau$.  
From this the assertion follows easily. 

We apply Lemma \ref{extc} to the pair $(\hat{S},\hat{C})$ and 
$F=\Omega_{\hat{C}}$,  
and obtain the isomorphism 
\begin{equation}\label{isomm} 
Ext^{i}_{O_{\hat{S}}}(\Omega_{\hat{C}}\otimes \hat{N},O_{\hat{S}}) \cong 
Ext^{i-1}_{O_{\hat{C}}}(\Omega_{\hat{C}}, O_{\hat{C}}). 
\end{equation}  
Substituting this isomorphism for $i=1$ to $Ext$ homomorphism 
obtained by applying 
$Ext^{1}_{O_{\hat{S}}}(-,O_{\hat{S}})$ to $b$ in (\ref{ooss})  
we get a map 
$Ext^{1}_{O_{\hat{S}}}(\Omega'_{\hat{S}}(\log \hat{C}),O_{\hat{S}}) 
\stackrel{\beta}{\ra} 
Ext^{1}_{O_{\hat{C}}}(\Omega_{\hat{C}},O_{\hat{C}})$, and 
its sheaf version  
${\mathcal Ext}^{1}_{O_{\hat{S}}}(\Omega'_{\hat{S}}(\log \hat{C}),O_{\hat{S}}) 
\stackrel{\beta '}{\ra} 
{\mathcal Ext} ^{1}_{O_{\hat{C}}}(\Omega_{\hat{C}},O_{\hat{C}})$, 
which is surjective since 
${\mathcal Ext}^{2}_{O_{\hat{S}}}(\Omega_{\hat{S}}(\hat{C}),O_{\hat{S}}) =0$. 
These fit into the following commutative diagram: 
\[
\begin{CD}
 Ext^{1}_{O_{\hat{S}}}(\Omega'_{\hat{S}}(\log \hat{C}),O_{\hat{S}}) @>\beta>> 	
Ext^{1}_{O_{\hat{C}}}(\Omega_{\hat{C}},O_{\hat{C}})			\\
 @VeVV @VdVV		\\
H^0({\mathcal Ext}^{1}_{O_{\hat{S}}}(\Omega'_{\hat{S}}(\log \hat{C}),O_{\hat{S}}))  @>v>> 
H^0({\mathcal Ext} ^{1}_{O_{\hat{C}}}(\Omega_{\hat{C}},O_{\hat{C}})) \cong 
\oplus_{p\in \bar{B}} \C_p	 	@>>> 	0		\\
 @.  @VpVV		\\
 	@. \C_{r}  \oplus_{p\in \hat{B}} \C_p 	 				
\end{CD}	
\] 
where $v=H^0(\beta ')$, $d$  is the isomorphism (\ref{exc}) and $p$ is the natural projection. 
Note that $e$ also is isomorphic by Lemma \ref{dame}.  
Then we put  $c:=pd\beta =pve$. 
In view of the exact sequence (\ref{ilm}),  
in order to prove that  $c$ is isomorphic it suffices to show that 
$pv:H^0({\mathcal Ext}^{1}_{O_{\hat{S}}}(\Omega'_{\hat{S}}(\log \hat{C}),O_{\hat{S}})) 
\ra \C_{r}  \oplus_{p\in \hat{B}} \C_p$ 
gives an isomorphism $u: H^0(O_{\hat{F}})\ra \C_{r}$ 
when restricted to the subspace $H^0(O_{\hat{F}})$ (cf.\ (\ref{ilm}). 
Indeed, along  $\hat{F}$, $\beta '$  becomes 
the natural sheaf surjection  $O_{\hat{F}}\ra \C_{r_1}\oplus \C_{r_2}$ and 
$u$ is just the associated map $H^0(O_{\hat{F}})\ra \C_{r}$, 
which is isomorphic.  \hfill q.e.d.  
 
\vspace{3 mm} 
Finally, we also record the following exact sequence deduced from (\ref{oO}):
\begin{equation}\label{eece} 
\begin{CD}
0 @>>>  Ext^{1}_{O_{\hat{S}}}(\Omega_{\hat{S}}(\log \hat{C}),O_{\hat{S}}) 	
@>>>  Ext^{1}_{O_{\hat{S}}}(\Omega'_{\hat{S}}(\log \hat{C}),O_{\hat{S}}) 	@>>> 			\oplus_{p\in \hat{B}} \C_p @>>> 0 . 
\end{CD}	
\end{equation}  

\vspace{3 mm} 
\ul{Properties of the Kuranishi families} 

Let 
\begin{equation}\label{krr} 
 \hat{g}: (\hat{\mathcal S},\hat{\mathcal C}) \ra \hat{T},\  
 (\hat{S}_o,\hat{C}_o)=(\hat{S},\hat{C}),\ o\in \hat{T} 
\end{equation}  
be the Kuranishi family of deformations of the pair $(\hat{S},\hat{C})$. 
For any  $p\in \bar{B}$ we denote by  
\[\hat{g}_p: \hat{{\mathcal C}}(p) \ra T_p,\ 
(\hat{C}_o,p_o)\cong (\hat{C},p),\ o\in \hat{T}_p,\]  
the Kuranishi family of deformations of the isolated singularity $(\hat{C},p)$. 
$\hat{T}_{p}$ is smooth of dimension one.  Any deformation of 
$(\hat{S},\hat{C})$ induces 
a deformation of $(\hat{C},p)$, 
and correspondingly we have a versal map  
$\hat{\tau}_p: \hat{T} \ra  \hat{T}_{p}$.  
The fiber  $\hat{T}(p):= \hat{\tau}_p^{-1}(o)$ is a hypersuface which is 
uniquely determined independently of the choice of $\hat{\tau}_p$. 

Recall that $m=k-2$ and  $\#\hat{B}=m$. 

\begin{Proposition}\label{arie}
Let the notations be as above. 
Then the Kuranishi space $\hat{T}$ is smooth of dimension $m+1$. 
The product map 
$\hat{\tau} 
:=\hat{\tau}_{r}\times \Pi_{p\in \hat{B}} 
\hat{\tau}_p : \hat{T}\ra \Pi_p \hat{T}_p$
is isomorphic, where $r=r_\al, \al = 1$ or  $2$. 
\end{Proposition} 

{\em Proof}. 
$\hat{\tau}$ is identified with the map $c$ in 
2) of Proposition \ref{diff}. Thus the proposition follows from 
that proposition.  
\hfill q.e.d.  

\vspace{3 mm} 
{\em Remark 4.1}. The proof shows that the singularities of  $\hat{C}$ 
at the two points $r_\al,\ \al =1,2$,  
and the singularities of $\hat{S}$ along $\hat{F}$ are simultaneously smoothed. 

\vspace{3 mm} 
Define $A:= \hat{T}(r)$.  
This is a smooth hypersurface by the above proposition, and 
is independent of the choice of $r=r_\al$ by the above remark.  
Also we consider the subspace $I:=\cap_{p\in \hat{B}}\hat{T}(p)$, 
which is a one dimensional smooth subspace of  $\hat{T}$ by the above proposition. 
Also we consider the subspace $I:=\cap_{p\in \hat{B}}\hat{T}(p)$, 
which is a one dimensional smooth subspace of  $\hat{T}$ by the above proposition. 

\begin{Lemma}\label{I}
The restriction of the family to $I$ is identified with the 
Kuranishi family of log-deformations of  $(\hat{S},\hat{C})$. 
\end{Lemma} 

{\em Proof}.  
In view of Lemma \ref{exp} and (\ref{eece}), for the Kuranishi family 
$(h: (\hat{{\mathcal S}'},\hat{{\mathcal C}}')\ra \hat{I}$,\  
$(S'_o,C'_o)= (\hat{S},\hat{C}),\ o\in \hat{I})$  of log-deformations 
of $(\hat{S},\hat{C})$,  
the base space  $\hat{I}$ is smooth of 
dimension one and is realized as a subspace of  $\hat{T}$. 
It is in fact a subspace of $I\subseteq \hat{T}$ 
since $I$ is the maximal subspace of $\hat{T}$ 
parametrizing log-deformations of  $(\hat{S},\hat{C})$. 
Then we must have $I=\hat{I}$ since both are smooth of dimension one.  
\hfill q.e.d. 

\vspace{3 mm} 
For later purpose (cf.\ Lemma \ref{defc}) we also 
give a more explicit construction of the above one-dimensional 
deformation by the method of Nakamura \cite[(4.2)]{nkm1}. 

First, we construct a local model of the deformations of the pair  
$(\hat{S},\hat{C})$ along  $\hat{F}$. 
For $m\in \Z$
let  $L_m$  be  the holomorphic line bundle of degree $m$  
on the complex projective line  $\B$, identified with its total space. 
Let  $f$  be the holomorphic function on  $V:=L_1\oplus L_{-1}$ given by 
the composition $V \ra  L_0  \ra  \C$, 
where the first arrow is given by the natural pairing and 
the second arrow is the natural projection 
from the product $L_0 = \B\times \C$.  
The fiber  $V_t:= f^{-1}(t),t\neq 0$, then gives a smoothing of the pair 
\[ (V_0,C_{0,0}\cup C_{0,\infty}))
:=(L_1\cup L_{-1},(L_{1,0}\cup L_{-1,0})\cup (L_{1,\infty }\cup L_{-1,\infty }))\] 
to  $(V_t,C_{t,0}\cup C_{t,\infty })$, 
where $C_{t,*}=p^{-1}(*)\cap V_t, *=0, \infty$, 
with  $p: V \ra  \B$ the natural projection.  
For any $t\neq 0$  
the projection  $q_t : V_t \ra  L_1-0$ is isomorphic and 
sends $C_{t,*}$ to $L_{1, *}-0$, where $0$ is the zero section.  

$(\hat{S},\hat{C})$ and  $(V_0,C_{0,0}\cup C_{0,\infty })$ are then isomorphic 
as germs along $\hat{F}$ and along the zero section respectively. 
(For this one uses the fact that $(\tilde{S},\tilde{C})$ is obtained from 
the toric projective plane with three fixed lines in general position, 
$(\B^2, l_0\cup l_1\cup l_2)$,  by blowing-up successively nodes 
on the anti-canonical cycle over the node $l_1\cap l_2$.) 

Since the deformation above is trivial off the zero section, 
the induced deformation of the germ  $(\hat{S},\hat{C})$ along $\hat{F}$ 
extends to a global log-deformation of $(\hat{S},\hat{C})$ 
which is trivial outside a neighborhood of $\hat{F}$.  Call a deformation 
obtained in this way a standard family of deformations of $(\hat{S},\hat{C})$. 

\begin{Lemma}\label{cag}
A standard family of deformations of $(\hat{S},\hat{C})$ is 
a Kuranishi family of log-deformations of $(\hat{S},\hat{C})$. 
\end{Lemma} 

{\em Proof}.  It suffices to show that 
the induced versal map $\tilde{\tau} : (\C,0) \ra (I,o)$ is isomorphic. 
This is true if the composite map $\hat{\tau}_r \tau : (\C,0) \ra (\hat{T}_r,o)$ 
is isomorphic.  But the latter is clear by the above construction. 
\hfill q.e.d.  

We call the original pair $(\tilde{S},\tilde{C})$ {\em minimal} 
if each irreducible component $D$  of  $\tilde{C}$  with $D^2=-1$ 
intersects either   $H$ or $E$. 
For $t \in \hat{T}$ 
we shall identify the fibers $(S_t,C_t)$ of $\hat{g}$ over $t$ in the next proposition. 
In order to state it 
we distinguish three cases: \\ 
Case-H: $\hat{C}^\al, \ \al =1,2$, are disjoint, 
and both of $\hat{C}^\al$ are reducible.  \\ 
Case-H$'$:  $\hat{C}$ is connected, i.e. 
the case where  $\varphi$ is of twisted type. \\ 
Case-P:  $\hat{C}^\al$ are disjoint, but 
one of $\hat{C}^\al$, say $\al =1$, is irreducible.  

Note that by our assumption  $k>2$ 
at most one of $\hat{C}^\al$ is irreducible. 

\begin{Proposition}\label{grlfb} 
1) Suppose that $t\notin A$.  Then 
$S_t$ is a smooth surface of class VII with 
second Betti number $m$.  
In Case-P and H (resp.\ Case-H$'$) 
$C_t, t\notin A$, is the unique anti-canonical curve 
(resp.\ $L_t$-twisted anti-canonical curve with $L_t=L_{S_t}$) on  $S_t$. 
The minimal model $\bar{S}_t$ of $S_t$ 
is isomorphic to one of the following surfaces in each case: 

Case-H:  a hyperbolic or parabolic Inoue surface or a diagonal Hopf surface, 

Case-H$'$: a half Inoue surface or a diagonal Hopf surface.  

Case-P: a parabolic Inoue surface or a diagonal Hopf surface. \\ 

2) The restriction of the family to $I$ is identified with the 
Kuranishi family of log-deformations of  $(\hat{S},\hat{C})$. 
Let $t\in I$ with  $t\neq o$.  Then 
the surface  $S_t$ is 
a properly blown-up hyperbolic (resp.\ half, resp.\ parabolic) Inoue surface 
in Case-H (resp.\ Case-H$'$, resp.\ Case-P).  
The isomorphism class of $S_t$ is independent of  $t$ in Case-H or -H$'$. 
If $(\tilde{S},\tilde{C})$ is minimal, then  $S_t$ also is minimal, namely 
$S_t$ is a hyperbolic (resp.\ half, resp.\ parabolic) Inoue surface. 
\end{Proposition} 

{\em Remark 4.2}. 
1) By Lemmas \ref{vai} and \ref{aut} we have 
$h^0(\Theta_{\hat{S}_t}(-\log \hat{C}_t))= 1$ for $t= o$,  
and $=0$ otherwise in Case-H and -H$'$, while $h^0(\Theta_{\hat{S}_t}(-\log \hat{C}_t))= 1$ 
for all $t\in I$ in Case-P. 
Correspondingly, the Kuranishi family 
of log-deformations over  $I$  above is universal in Case-P and not 
in Case-H or -H$'$.  

2) 
$\bar{S}_t$ is a diagonal Hopf surface if and only if 
$t\notin \hat{T}(p)$ for any $p\in \hat{B}$ by Lemma \ref{knk}.  
\vspace{3 mm} 

{\em Proof}.  
By the definition of $A$ 
it is clear that  $S_t$ is smooth if and only if 
$t\notin A$ (cf.\ Remark 4.1). 

First we consider the restriction of  $\hat{g}$  to  $I$, which 
may be identified with the Kuranishi family of log-deformations of   
$(\hat{S},\hat{C})$ by Lemma \ref{I} 
and also with the family constructed before. 
Lemma \ref{cag} by that lemma.  Note that $I\cap A=\{o\}$. 
In this case by \cite[Th.44]{kd3} 
the general fiber $S_t, t\in I-o$ is, topologically,  
obtained from $\hat{S}$ by a spherical modification (cf.\ \cite[(3.3)]{nkm1}).  
Then by \cite[(3.4)]{nkm1} $S_t$ has infinite cycle fundamental group and 
has $m=k-2$ as the second Betti number.  
Moreover, since  $\hat{C}$ is an anti-canonical 
curve on  $\hat{S}$ as we have already noted, $h0(-\hat{K})\geq 1$ and hence $h^0(\hat{K}^n)=0$ 
for all $n> 0$, where $\hat{K}=K_{\hat{S}}$ 
Thus by the upper semicontinuity 
we get that  $S_t,t\neq o$, all have Kodaira dimension zero.  Hence  $S_t$, are 
surfaces of class VII for all $t\neq o$ and hence for all $t\notin A$. 

Indeed, actually more precisely, when $(\hat{S},\hat{C})$ is minimal in the sense defined above,  
$S_t$  is a hyperbolic Inoue surface for any $t\in I-0$ by the more precise 
computation of the self-intersection numbers of the irreducible components of 
$C_t$ in  $S_t$  due to \cite[(5.13)(5.14)]{nato} (cf.\ also below).  
In fact, from Lemma \ref{cag} we may identify our family  $g'$ with 
that used by Nakamura in \cite{nato}. 
Using the unique extension theorem of $(-1)$-curves \cite{kds} 
the general case can easily be reduced to the minimal case 
by contracting all the $(-1)$-curves contained in the irreducible components of 
$C_t, t\in I$, successively and simultaneously to the points. This shows 2).  

Next we show 1). 
Suppose first that we are in Case-H or -P.  
We show that 
$K_t+C_t$ is trivial for all $t\notin A$, 
where  $K_t$ is the canonical bundle of $S_t$. 
Note first that the family is versal at any point of $\hat{T}$ 
by the openness of versality (cf.\ \cite{bi}). 
Thus by 1) of Proposition \ref{seis} 
for some open neighborhood $U$ of $I-\{o\}$ in  $\hat{T}$, 
this is true. 
Then by analytic continuation  $K_t+C_t$ are trivial for all $t \in \hat{T}-A$.  
Indeed, the component of the identity section 
$P:=\mbox{Pic}_0\hat{{\mathcal S}}/\hat{T}\ra \hat{T}$  
of the relative Picard variety associated to  $\hat{g}$ 
is a principal $\C^*$-bundle 
at least over $\hat{T}-A$ and 
$K_t+C_t$ defines a holomorphic section of   $P$  over  $\hat{T}-A$  
which is trivial over  $U$, and hence over the whole  $\hat{T}-A$.   
(In fact, $P$ itself can be shown to be seperated as the singular fiber  $\hat{S}$ is 
irreducible and reduced.) 

Once this is proved, 
noting that $C_t$ is disconnected 
we obtain the structure of the surface $\bar{S}_t$ 
as stated in the proposition by Lemma \ref{atc}.  
This finishes the proof of 1) in Case-H and -P. 
Case-H$'$ is treated similarly by using 1) of Proposition \ref{seis'} and 
Lemma \ref{atch}.  \hfill q.e.d.  

\vspace{3 mm} 
We now restrict to Case-H. 
For $t\in I-o$ 
we consider the log-deformations $(S_t,C_t)$ of $(\hat{S},\hat{C})$. 
Fixing $\al , \al =1,2$, 
write  $\hat{C}^\al= \hat{B}_0+\cdots +\hat{B}_{h_\al}$ cyclically 
such that $\tilde{B}_0\cap \hat{H}\neq \emptyset \neq \tilde{B}_{h_\al}\cap \hat{E}$, 
where $\tilde{B}_i = n^{-1}(\hat{B}_i)$.  
Then the cycle  $C^\al_{t}$ which is a deformation of $\hat{C}^\al$ 
is written as  $B_{0,t}+\cdots +B_{h_\al -1 ,t}$, 
where $B_{i,t}$ are deformations of $\hat{B}_i$ 
for  $i\neq  0$, and 
$B_{0,t}$ is one of $\hat{B}_{h_\al}+\hat{B}_0$.  
This also gives a natural numbering of the cycle $C_t^\al$, 
which we call tentatively a {\em toric} numbering. 
Clearly for $i \neq 0$ we have 
$(B_{i,t})^2=\tilde{B}_i^2$ for the self-intersection numbers.  
For $B_{0 , t}$ we note that 
$(B_{0, t})^2 = (\tilde{B}_{h_{\al}}+\tilde{B}_{0})^2$ if $h_\al > 0$ and 
$(B_{0, t})^2 = (\tilde{B}_{0})^2$ if $h_\al = 0$ 
since the intersection number is invariant under normalization and 
deformation. 
Thus we get the following: 
\begin{Lemma}\label{intt}
\begin{gather}\label{int}
(B_{0,t})^2 = \tilde{B}_{0}^2, \ \ \rm{if} \ \  h_\al =0, \\
(B_{0,t})^2 = \tilde{B}_{0}^2 +\tilde{B}_{1}^2+2, \ \ \rm{if}\ \  h_\al =1,  \\
(B_{0 , t})^2 = \tilde{B}_{0}^2 +\tilde{B}_{h_\al}^2, \ \ \rm{if}\ \  h_\al >1. 
\end{gather} 
\end{Lemma} 
This lemma is a prerequisite for proving 
the following result due to Nakamura \cite[(5.13)(5.14)]{nato}. 
\begin{Proposition}\label{cmp}
Any properly blown-up hyperbolic, half or parabolic Inoue surface $(S,C)$ is obtained by 
log-deformations of a rational surface $(\hat{S},\hat{C})$ 
with a nodal curve obtained 
from an admissible toric surfaces $(\tilde{S},\tilde{C})$  as above. 
\end{Proposition} 

Although in \cite{nato} only the minimal case has been treated, the 
generalization to the blown-up case is immediate by the simultaneous 
blowning-up of the family of deformations in the minimal case. 

\vspace{3 mm} 
A toric surface which gives rise to the given 
hyperbolic or half Inoue surface is not unique.  In fact, we note 
in the following result 
that Nakamura's construction acturally give rise to $m$ such toric surfaces 
according as which pair of irreducible components  $C^\al_t, \al =1,2$, is 
to be bent and broken. 

\begin{Proposition}\label{m} 
Let  $S$  be any properly blown-up hyperbolic Inoue surface with second Betti number $m$.  
Then there exist in general  $\bar{m}$  admissible toric surfaces $(\tilde{S},\tilde{C})$ 
which give rise to $S$  as a deformation $S_t$ of $\hat{S}$  as above, where 
$\bar{m}$ is the second Betti number of the minimal model of  $S$.  
\end{Proposition} 

{\em Proof}.  
Suppose first that $S$ is minimal and 
that the canonical weight sequences of  $S$  are given by (\ref{ws1}) and (\ref{ws2}) 
with  $n$ and $k_i$ with $n>0$, $k_i\geq 1$ and $1\leq i\leq 2n$  fixed. 
Let  $(\B^2, l_0\cup l_1\cup l_2)$  be 
a toric projective plane with three fixed lines in general position. 
Then our toric surface  $(\tilde{S},\tilde{C})$  are given by 
a finite succession of blowning-ups such that the center of each blowing up 
is mapped to the point  $l_1\cap l_2$ and such that 
the proper transform of  $l_2$ corresponds to the weight $0$  in (\ref{ws3}) below, 
where the last condition comes from the minimality assumption. 
Consider the following $k_{2n}$ 
canonical weight sequences of $\tilde{C}^1,  \tilde{C}^2$ 
considered modulo interchanging $\tilde{C}^1$ and $\tilde{C}^2$ 
with the same notational convention as in (\ref{ws1}) (\ref{ws2}): 
\begin{gather}\label{ws4}
k_1, \ \ \   [k_2-1],k_3+2, \ldots ,k_{2n-1}+2,[k_{2n}] \\
0, [k_1-1],k_2+2,\ldots ,[k_{2n-1}-1],k_{2n}+2 \label{ws3}
\end{gather} 
and 
\begin{gather}\label{ws5}
0, [\kappa -1],k_1+2,[k_2-1],\ldots ,k_{2n-1}+2,[k_{2n}-\kappa ] \\
\kappa , \ \ \    [k_1-1],\ldots , [k_{2n-1}-1],k_{2n}-\kappa +2\label{ws6}
\end{gather} 
where $0<\kappa <k_{2n}$.  
Then for each of the above pairs of weight sequences 
we can find a unique minimal admissible toric surface $(\tilde{S},\tilde{C})$ 
having this sequence as its canonical weight sequences.  
Since the admissible toric surfaces are determined completely by the 
canonical weight sequences, the uniqueness is clear, but 
in fact, one can show that there exists a unique way to obtain  $(\tilde{S},\tilde{C})$ 
from $(\B^2, l_0\cup l_1\cup l_2)$ by a finite succession 
of the blowing-ups as above. (The details is omitted.) 

By Lemma \ref{intt} it is immediate to see that any of the above toric surfaces 
gives rise to the hyperbolic Inoue surface with canonical weight sequences 
(\ref{ws1}) and (\ref{ws2}). 

Finally the canonical weight sequences of hyperbolic Inoue surfaces of the above form 
are determined up to cyclic permutations. 
After taking into account all the sequences after such permutations 
and recalling that $\sum_{1\leq i \leq 2n} k_i = m$, 
in all we obtain $m$ minimal admissible toric surfaces. 
From these we conclude the proof of the proposition in the minimal case. 

Suppose next that  $S$ is not minimal and is obtained as 
a log-deformation $S=S_t$ of 
an admissible toric surface  $\hat{S}$.   
As follows easily from Lemma \ref{intt} 
no $(-1)$-curve in $C$ is never 
an amalgamated deformation of the union of two irreducible components 
of $\hat{C}$  intersecting along  $\hat{F}$.  Thus any $(-1)$-curve in  $C$ 
is a deformation of a $(-1)$-curve in  $\hat{C}$  contained in the smooth locus of $\hat{S}$. 
Then by blowing down these $(-1)$-curves simultaneously we get in a canonical 
way the minimal model $\bar{S}$ of $S$  
as a log-deformation of the corresponding `minimal model' of $\hat{S}$.  
This reduces the number in question for $S$ to the corresponding number 
for its minimal model. 
\hfill q.e.d.

\vspace{3 mm} 
{\em Remark 4.3}. 
The $m$ admissible toric surfaces obtained in the proposition would all be 
distinct if the weight sequences are general enough.  
On the other hand, 
for a parabolic Inoue surface  $S$ with second Betti number  $m$ 
the admissible toric surface which gives rise to  $S$  
by the above process is uniquely given by the following 
weight sequences : 
\begin{gather}\label{wpb} 
	m \\
 0 \ [m]  . 
\end{gather} 

\vspace{3 mm} 
Let  $S$  be a properly blown-up hyperbolic Inoue surface.  
By the previous results $(S,C)$ is 
obtained as a log-deformation  $(S_t,C_t)$  
of some admissible toric surface  $(\hat{S},\hat{C})$.  
In this case for each cycle  $C^\al$  we have the toric numbering 
for its irreducible components as defined before Lemma \ref{intt}. 
We show that this toric numbering indeed coincides with the canonical numbering 
as characterized by Lemma \ref{d}. This implies also that the toric 
numbering is independent of the initial datum  $\hat{S}$. 

\begin{Lemma}\label{defc} 
The canonical numbering and toric numbering coincide in the sense explained above.
\end{Lemma} 

{\em Proof}.  
By Lemma \ref{cag} we may consider the standard family of deformations of 
$(\hat{S},\hat{C})$.  Then by using the notations there we may assume that 
$C^1_t = C_{t,0}$ and $C^2_t=C_{t,\infty}$ in a neighborhood $W$ of $\hat{F}$ 
(in the total space of the family).  For instance we consider $C^1_t$. 
$S_t$ contains the global spherical shell $U_t$ (identifying $S_t$ with $V_t$ in  $W$), 
and we have $U_t\cap C^1_t\subseteq B_{0,t}$. Then from the fact that 
tubular neighborhoods of the zero sections in $L_1$ and $L_{-1}$ is strongly pseudoconcave and 
and strongly pseudoconvex respectively and by the definition of toric numbering 
we immediately see by Lemma \ref{d} that the two numberings coincide. 
\hfill q.e.d.

\section{Construction of a singular twistor space} \label{cstr}

\vspace{3 mm} 
\ul{Twistor spaces}

Let  $M$  be an oriented compact $C^{\infty}$ 4-manifold and $[g]$ a self-dual structure, 
i.e., the conformal class of a self-dual metric, on  $M$. 
Denote by  $Z$  the twistor space associated to the self-dual manifold 
$(M,[g])$ with the twistor fibration  $t: Z \ra  M$, 
which is a $\B$-bundle where $\B$ is the complex projective line (\cite{ahs}).  

Any fiber of $t$  is called a {\em twistor line}. 
There exists an anti-holomorphic involution $\sigma $ of  $Z$,  
called the {\em real structure} of $Z$,   
which is fixed point free and preserves each twistor line. 
A complex surface  $S$  on  $Z$  is called {\em elementary} 
if the intersection number  $LS=1$ for any twistor line  $L$ on  $Z$.  
If  $S$ is an elementary surface, 
then its {\em conjugate} $\bar{S}:=\sigma (S)$ is again an elementary surface.  

There are two cases to consider: 
\begin{quote}  
Case 1: $S$  contains a twistor line, \\
Case 2: $S$  contains no twistor lines. 
\end{quote}    
In Case 1 the structure of  $S$  is as follows (cf.\ \cite[Lemmas 1.9, 1.10]{p92}): 
\begin{Lemma}\label{elm} 
Let  $S$  be an elementary surface in Case 1.  Then  $S$  contains 
precisely one twistor line, say  $L$, $S$ and $\bar{S}$ intersects transversally 
along  $L$,  $S$  is obtained from a complex projective plane  $\B^2$ by 
a succession of blowing-ups such that the total blowing-down   $S \ra  \B^2$ 
maps a neighborhood of  $L$  isomorphically onto a neighborhood of a line on $\B^2$.  
In particular $L^2 = 1$ in  $S$.  
$M$ is diffeomorphic to $m\B^2$ with $m=b_2(M)$, 
and the restriction of the twistor fibration $t$  to  $S$ 
is nothing but the smooth contraction of  $L$  to a point $t(L)$ of $M=m\B^2$. 
\end{Lemma}   

On the other hand, in Case 2 we easily check the following: 
\begin{Lemma}\label{c2} 
In Case 2,  $S\cap \bar{S} = \emptyset$ and 
they are mapped diffeomorphically onto $M$. 
If $J$ and $\bar{J}$  are the complex structures  on  $M$  
induced from  $S$  and $\bar{S}$  via this diffeomorphism respectively, 
then $\bar{J}=-J$.  Moreover, the self-dual structure  $[g]$ is compatible 
with $\pm J$  and gives the anti-self-dual hermitian surface 
$(M,\pm J,[g])$ which are complex conjugate to each other. 
\end{Lemma} 

\vspace{3 mm} 
\ul{Joyce twistor space $Z$}

Let  $m$  be a positive integer and  $m\B^2$ the connected sum of 
$m$ copies of complex projective plane  $\B^2$.   
Fix any such  $m$  and write  $M=m\B^2$. 
Denote by  $K:=S^1\times S^1$ the real $2$-torus.  
Then $M$ admits a finite number $\psi (m)$ of smooth effective $K$-actions on  $M$ 
up to diffeomorphisms.  (For instance  $\psi (m)=1$ for $m=1,2$ and $\psi (3)=3$.) 
For each such smooth $K$-action on  $M$ 
Joyce \cite{jo} constructed a smooth connected family 
of $K$-invariant self-dual conformal structures on  $M$, 
depending on $(m-1)$ real smooth parameter. 

  Fix any such $K$-invariant self-dual structure on  $M$ and denote it 
by  $[g]$.   Let  $Z$  be the associated twistor space 
with natural projection  $t: Z \ra  M$, making  $Z$  a smooth $\B$-bundle over  $M$. 
Denote the real structure of $Z$ by $\sigma$ as before. 
In general we call such a $Z$ a {\em Joyce twistor space}. 

The $K$-action on $M$ naturally lifts to a holomorphic $K$-action on  $Z$. 
This action on $Z$ then extends to a holomorphic action 
of the complexfication  $G= \C^*\times \C^*$ of $K$, which is 
an algebraic torus of dimension two. 

The structure of  $Z$  with this $G$-action has been studied in detail in \cite[\S 4, \S 6]{fjj}. 
We explain some of the structures of  $Z$  
which are important for us in what follows. 
We set  $k=m+2\geq 2$. 
Then  $Z$  admits exactly $k$ pairs of $G$-invariant elementary surfaces 
$\{(S^+_i,S_i^-)\}, 1\leq i\leq k$, with $\sigma(S_i^\pm )=S_i^\mp$.  
$S^\pm_i$ are projective smooth toric surfaces 
with respect to the induced $G$-action 
and $S^+_i$ and $S_i^-$ intersect transversally 
along a $G$-invariant twistor line  $L_i$ (cf.\ \cite[Prop.6.12]{fjj}).  
The self-intersection numbers of $L_i$ in  $S_i^\pm $ both equal one; 
\[ L_i^2=1 .\] 
The point  $p_i:=t(L_i)$ is a fixed point of  $K$  on  $M$ and 
this sets up bijective 
correspondences among the set of pairs of $G$-invariant elementary surfaces 
on  $Z$, the set of $G$-invariant twistor lines, and 
the set of $K$-fixed points on  $M$.  

The union  $S_i:=S_i^+\cup S_i^-, 1\leq i\leq k$, all belong 
to the fundamental system $|-\frac{1}{2}K|$, where 
$\frac{1}{2}K$ is the canonical square root of the canonical bundle  $K$  of  $Z$.  
The subspace $H^0(Z,-\frac{1}{2}K)^G$ of $H^0(Z,-\frac{1}{2}K)$ 
of $G$-invariant elements 
is two dimensional and the associated pencil $|-\frac{1}{2}K|^G$ is important 
for the study of the structure of $Z$.  
The base locus  $C$  of $|-\frac{1}{2}K|^G$  is 
a cycle of rational curves which are both 
$G$- and $\sigma$-invariant and is of the form 
\begin{equation}\label{clc} 
 	C = C_1^++\cdots +C_k^++C_1^-+\cdots +C_k^- , \ \ C^\pm_i\cong \B,  
\end{equation}  
with  $\sigma(C_i^\pm )=C_i^\mp$.  The $G$-action is free outside 
$C\cup (\cup_i L_i)$ (cf.\ \cite[Prop.4.4]{fjj}).  

A general member $S_0$ of $|-\frac{1}{2}K|^G$ is a smooth toric surface 
with respect to the induced $G$-action with anti-canonical cycle  $C$. 
In fact, any smooth member are all isomorphic to each other with the same 
weight sequence of the form  
\begin{equation}\label{0w} 
(a_1,\ldots ,a_k,a_{1},\ldots ,a_k).  
\end{equation}  
There exist precisely $k$ singular members of the pencil $|-\frac{1}{2} K|^G$. 
They are the surfaces $S_i$ above with two irreducible components $S^\pm_i$. 

We put 
$p_i^\pm =C_i^\pm \cap C_{i+1}^\pm, \ 1\leq i\leq k$,   
with the convention that $C_{k+1}^\pm =C_1^\mp $.  
Thus  $\sigma(p_i^\pm )=p_i^\mp$. 
$S^\pm_i$ contains exactly half of the cycle $C$, i.e., 
$C\subseteq S_i$ and 
the intersection $C^\pm_{(i)}:=S^\pm_i\cap C$ is a chain of rational curves given by 
\[ C^\pm_{(i)}=C_{i+1}^\mp +\cdots +C_{k}^\mp + C_{1}^\pm + \cdots +C_i^\pm .  \] 
Then the anti-canonical cycle  $B^\pm_i$ of the toric surface $S_i^\pm $ 
is written as 
\[ B^\pm_i=L_i+C^\pm_{(i)}. \] 
Moreover, for each $j\neq i$, $L_i$ intersects with  $S^\pm_j$ transversally 
at the unique points $p_i^\pm $ and $L_i\cap C=\{p_i^\pm\}$. 
The weight sequence of $S^\pm_i$ is then given by: 
\begin{equation}\label{wsws} 
 (1,a_{i+1}+1,a_{i+2},\ldots,a_k,a_1,\ldots ,a_i+1)  
\end{equation}  
independently of  $\pm$ (cf.\ \cite[(13)]{fjj}). 
We call the pair $(i,j)$ {\em minimal}, if for any $d$, $a_d=1$ implies that $d=l$ or $l+1$ 
where $l=i,j$.  
The next lemma is used in proving that our construction covers all the 
hyperbolic Inoue surfaces.  

\begin{Lemma}\label{tj} 
For any projective toric surface  $S$  
with a $(+1)$-curve $H$ in its anti-canonical cycle 
there exist a Joyce twistor space $Z$ as above and an index  $i, 1\leq i\leq k$,  
such that $S^\pm_i$ are both isomorphic to  $S$.  
\end{Lemma} 

{\em Proof}.  Consider the induced $K$-action on  $S$.  Then we may 
$K$-equivariantly contract the curve  $H$  to a point $x$  
of a smooth $K$-manifold $M$,  
which is necessarily  diffeomorphic to  $m\B^2$, 
where  $m+1$ is the second Betti number of  $S$.   
The point $x$ is one of the $k$ fixed points of the $K$-action on  $M$.  
Take any $K$-invariant self-dual structure on  $M$  of Joyce and 
in the associated twistor space 
take the pair $\{S^\pm_i\}$ of elementary surfaces corresponding 
to  $x$  in the sense mentioned above.  By Lemma \ref{elm} and 
the above description we see that  $S$  and $S^\pm_i$ are $K$-diffeomorphic 
with respect to the induced $K$-action.  Since the induced $K$-action 
determines the toric surface as a complex surface, we are done.   \hfill q.e.d.

 \vspace{2 mm} 
\ul{Blown-up twistor space $\tilde{Z}$} 
	
Fix $i, j$ with $1\leq i<j\leq k$ and write  $l$  for $i$ and/or $j$. 
(Note however that since we consider 
$i$ and $j$ cyclically modulo $k$, 
the roles of  $i$  and  $j$ are symmetric.) 

Let  $\mu : \tilde{Z} \ra  Z$  be the blowing-up 
with center the disjoint union $L_i\cup L_j$ and 
with exceptional divisors  $Q_l:=\mu ^{-1}(L_l)$, $l=i,j$.  
$Q_l$ are isomorphic to the product  $\B\times \B$  
with $\mu |Q_l : Q_l \ra  L_l$  identified with 
the projection $\B\times \B \ra  \B$, say to the first factor. 
Then the normal bundle $N_{Q_l/ \tilde{Z}}$ of $Q_l$ in  $\tilde{Z}$ is isomorphic 
to the line bundle $O(1,-1)$  of bidegree $(1,-1)$ on $Q_l\cong \B\times \B$.  

Let  $\tilde{S}^\pm_l$  be 
the proper transforms of  $S^\pm_l$ in $\tilde{Z}$.  
From the construction we see that 
$\tilde{S}^\pm_l$ are disjoint, 
but the intersection of any other pairs 
from the four surfaces  $\tilde{S}^\pm_l, l=i,j$, 
is non-empty and 
consists of a chain of rational curves, which is a connected component of 
the proper transform in $\tilde{Z}$ of the cycle  $C$. 
(See the formulae for its image in $\hat{Z}$ in (\ref{co1})(\ref{co2}) below.) 
Now write $\tilde{S}_l$ for the disjoint union $\tilde{S}^+_l\cup \tilde{S}^-_l$ and 
set   $\tilde{S}:=\tilde{S}_i\cup \tilde{S}_j$.  The latter is 
a connected surface with four irreducible components.  

The actions of  $G$  and $\sigma$ naturally lift to  $\tilde{Z}$ with  
$\sigma (\tilde{S}^\pm_l)=\tilde{S}^\mp_l$ and $\sigma (Q_l)=Q_l$. 
$\tilde{S}^\pm_l$ and $Q_l$ are $G$-invariant and will be considered 
as toric surfaces with respect to the induced $G$-action.  

We put  $H^\pm_l:=\tilde{S}^\pm_l\cap Q_l,  E^\pm_l:=\tilde{S}^\pm_l\cap Q_{l'}$, 
where $\{l,l'\}=\{i,j\}$. 
$H^\pm_l$ is mapped isomorphically onto  $L_l$ by $\mu$, 
and $\tilde{S}^\pm_l \ra  S^\pm_l$  is the blowing up of  $p^\mp_j$  
if $l=i$ and of $p^\pm_i$ if $l=j$ 
with exceptional curve $E^\pm_l$. 
Thus we get  
\[(H^\pm_l)^2=1\ \ \mbox{and} \ \ (E^\pm_l)^2=-1 \ \ \mbox{in}\ \  
\tilde{S}^\pm_l \] 
while 
$(H^\pm_l)^2=(E^\pm_l)^2=0$  in  $Q_l$.  Further 
the anti-canonical cycles $\tilde{B}^\pm_l$ of $\tilde{S}^\pm_l$ is 
given e.g. when $l=i$ by 
\[ \tilde{B}^\pm_i= H_i^\pm +\tilde{C}^\mp_{i+1}+\cdots +\tilde{C}^\mp_j+E^\pm_i+\tilde{C}^\mp_{j+1}+\cdots +\tilde{C}^\pm_i \]  
with the same weight sequence (independently of $\pm$) 
\begin{equation}\label{wsl} 
 (1,a_{i+1}+1,a_{i+2},\ldots,a_j-1, -1, a_{j+1}-1,\ldots ,a_k,a_1,\ldots ,a_i+1) 
\end{equation}  
where $\tilde{C}^\pm_d$  is 
the proper transform of  $C^\pm_d$ in  $\tilde{Z}, 1\leq d\leq k $.   
Similarly, the anti-canonical cycles 
$F_l$ of $Q_l$ is given by 
\[ F_l = H^+_l + E^+_{l'} + H^-_l + E^-_{l'}. \]  

\vspace{2 mm} 
\ul{Singular twistor space  $\hat{Z}$} 

Now 
we choose and fix  
an isomorphism of the pairs 
\begin{equation}\label{fai} 
 \varphi : (Q_i,F_i) \ra  (Q_j,F_j)   
\end{equation}  
which maps  $H_i^\pm $ (resp.\  $E^\pm_j$) to  $E^\pm_{i}$ (resp.\ $H^\pm_j$), 
thus interchanging the horizontal and vertical directions.   
Let  $\hat{Z} $  be 
the complex space obtained by identifying in  $\tilde{Z}$ 
the subspaces $Q_i$ and $Q_j$ via  $\varphi $.  
Let  $\nu : \tilde{Z} \ra  \hat{Z} $ be the quotient map, which is considered 
as the normalization map of  $\hat{Z} $.    
Let  $\hat{Q}:=\nu (Q_i)=\nu (Q_j)$ be the singular locus of  $\hat{Z} $ and  
$\hat{S}^\pm_l:=\nu (\tilde{S}^\pm_l)$  
the image of  $\tilde{S}^\pm_l$  in  $\hat{Z} $.  
Then $\hat{S}^\pm_l$ is a non-normal surface with singular locus 
$\hat{F}^\pm_l:=\nu (H^\pm_l)=\nu (E^\pm_{l})$  ($= \hat{Q}\cap \hat{S}^\pm_l $). 
The image  
\[ \hat{F}:= \hat{F}^+_i+\hat{F}^+_j+\hat{F}^-_i+\hat{F}^-_j \] 
of  $\nu(F_i)=\nu(F_j)$  in  $\hat{Q}$ 
belongs to the anti-canonical system on  $\hat{Q}$ and 
shall be called the anti-canonical cycle of  $\hat{Q}$.  
(Note that since the $G$-action is not $\varphi$-equivariant,  
$\hat{Q}$  has no natural structure of a toric surface in general 
(cf.\ Prop.6.3 below).)

Let  $\hat{S}_l=\hat{S}^+_l\cup \hat{S}^-_l$.  
Then  $\hat{S}:= \hat{S}_i\cup \hat{S}_j=\nu (\tilde{S})$  is 
a surface in  $\hat{Z} $ consisiting of four irreducible components $\hat{S}^\pm_l, 
l=i,j$.  
By our construction $\varphi $ maps the intersection points  
$\tilde{C}_i^\pm \cap  H_i^\pm$  and $\tilde{C}^\mp_{i+1}\cap H^\pm_i$ 
to $\tilde{C}^\mp_{j+1}\cap E^\pm_i$ and $\tilde{C}_j^\mp\cap E_i^\pm$ respectively. 
This implies that if we set $\hat{C}^\pm_d=\nu(\tilde{C}^\pm_d)$, 
the curves 
\[ \hat{C}^\mp_{j+1}+\cdots +\hat{C}^\pm_i \ \ \mbox{and} \ \ 
 \hat{C}^\mp_{i+1}+\cdots +\hat{C}^\mp_j \] 
form (four disjoint) cycles 
of rational curves on  $\hat{Z}$.  Moreover, we have 
\begin{align} 
 \hat{S}^{\pm}_i\cap \hat{S}^\pm_j & =  \label{co1}
\hat{C}^\mp_{j+1}+\cdots +\hat{C}^\pm_i \\ 
\hat{S}^\pm_i \cap \hat{S}^\mp_j & \label{co2} 
= \hat{C}^\mp_{i+1}+\cdots +\hat{C}^\mp_j. 
\end{align} 
In this way we see that 
each of the four surfaces $\hat{S}^\pm_l$ contains a pair of disjoint 
cycles of rational curves.  (Thus our choice of $\varp$  in (\ref{fai}) amounts to assuming that 
the restrictions of $\varp$  to all surfaces $\tilde{S}^\pm_l$ (cf.\ (\ref{idf}) are 
of untwisted type in the sense defined there.) 

We denote by  $\hat{C}^\pm_{(l)}$ 
the union of these cycles on $\hat{S}^\pm_l$.  
Note also that 
$\hat{S}^{+}_l$ and $\hat{S}^{-}_l$ are disjoint and 
no three of $\hat{S}^\pm_l$ have common points.  

In what follows it is convenient to distinguish the following two cases: 

Case-P (parabolic case)  $j=i+1$, or $(i,j)=(1,k)$  

Case-H (hyperbolic case)  otherwise. 

In other words, when we consider  $i, j$  cyclically modulo $k$, 
Case-P is precisely the case where  $i$  and $j$ are adjacent. In fact, 
precisely in this case one of the 
intersections $\hat{S}^\pm_i\cap \hat{S}^\pm_j=\hat{C}^\mp_j$ or  
$\hat{S}^\mp_j\cap \hat{S}^\pm_i=\hat{C}^\pm_i$ becomes irreducible and 
is a single rational curve with a node.  

Actually, the constructions in Section 6 and the arguments in Section 8 below 
for Case-H all apply to the twisted cases where 
where (some of) the restrictions become of twisted type. 
However, the results in terms of the bihermitian structures are somewhat 
different.  For this reason and also for the simplicity of exposition we treat 
the twisted cases separately in Section 9. 

\section{Structure of a singular twistor space} 

$\hat{Z}$ and $\hat{S}^\pm_l$ are complex spaces with normal crossing singularities. 
Therefore we may speak of 
the canonical bundles $\hat{K}$ and $K_{\hat{S}^\pm_l}$ of the respective spaces, 
corresponding to the dualizing sheaves. 

We also note that $\hat{S}$ and $\hat{S}^\pm$ are both Cartier divisors on  $\hat{Z}$, 
and similarly $C^\pm_{(l)}$ is a Cartier divisor on $\hat{S}^\pm_l$ 
(cf.\ b) below). 

\ul{Anti-canonical system}

We first identify the anti-canonical divisors of  $\hat{Z}$ and $\hat{S}^\pm_l$. 

\begin{Lemma}\label{lll} Let $\hat{K}$ and $K_{\hat{S}^\pm_l}$ be the canonical 
bundle on  $\hat{Z}$ and $\hat{S}^\pm_l$ respectively. 
Then we have 
\[ -\hat{K} = \hat{S}\ \ \mbox{and} \ \ -K_{\hat{S}^\pm_l}=\hat{C}^\pm_{(l)}.    \] 
\end{Lemma} 

{\em Proof}.  
Let  $S:=S^+_i+S^-_i+S^+_j+S^-_j$.  
$S$  is a member of the anti-canonical system  $|-K|$  of  $Z$. 
Then by the adjunction formula for the blowing up $\mu$ we get 
\begin{equation}\label{ad} 
  \tilde{K}=\mu^*K+Q_i+Q_j = -\mu^*(S)+Q_i+Q_j = -\tilde{S}-(Q_i+Q_j).  
\end{equation}  
where  $\tilde{K}$  is the canonical bundle of  $\tilde{Z}$. 
Then $\hat{K}$ is obtained 
by identifying  $(\tilde{K}+ Q_i)|Q_i \cong K_{Q_i}$ 
and $(\tilde{K}+ Q_j)|Q_j \cong K_{Q_j}$ along $\hat{Q}$ 
by the defining isomorphism $\varphi$ of $\hat{Z}$ (cf.\  \cite[(2.11)]{fr}). 
On the other hand, by (\ref{ad}) and the adjunction formula we have 
$(\tilde{K}+ Q_l)|Q_l = -\tilde{S}\cap Q_l = F_l$ 
and  $\varphi$ induces an isomorphism  $F_i \stackrel{\sim}\ra F_j$.  
Thus  $\hat{S}=\nu(\tilde{S})$ is a member of $-\hat{K}$, giving the first 
equality.  The proof of the second equality is similar.  \hfill q.e.d.

\vspace{3 mm} 
\ul{Local structure of ($\hat{Z},\hat{S}$)}

\vspace{2 mm} 
a) {\em Tangential points} 

We put $\hat{p}^\pm_n = \hat{C}^\pm_n\cap \hat{C}^\pm_{n+1}$ 
for $1\leq n\leq k$ with the convention $\hat{C}^\pm_{k+1}=\hat{C}^\mp_1$. 
Outside $\hat{Q}$ 
the intersections $\hat{S}^\pm_i\cap \hat{S}^\pm_j$ 
(resp.\ $\hat{S}^\mp_j\cap \hat{S}^\pm_i$) are 
transversal except at the points 
$\hat{p}^\pm_n$ with $j+1\leq n\leq k$ 
and $\hat{p}^\mp_n$ with $1\leq n\leq i-1$
(resp.\ $\hat{p}^\mp_n$ with $i<n<j$), 
where the corresponding two components of $\hat{S}$ 
intersect and are tangent to each other; in fact, locally, with respect to 
suitable local coordinates $x,y,z$ at such a point of  $\hat{Z}, \hat{S}$ has 
a local equation  
\begin{equation}\label{lt} 
 z(z-xy)=0.
\end{equation}  
Note that there exist in all $2m$ ($=(2k-4)$) such tangential points.  

\vspace{2 mm} 
b) {\em Points of $\hat{F}$}

Let $r$  be any one of the four singular points of  $\hat{F}$, i.e., the 
intersection points of the irreducible components  $\hat{F}^\pm_l$.  

The local structure of the pair  $(\hat{Z} ,\hat{S})$ at $r$  is described as follows.  
Let  $X = \C^2(u,v)$ and $Y=\C^2(x,y)$.  
Let  $A$  and  $D$  be the curves in  $X$  and $Y$  defined 
by  $uv=0$ and $xy=0$ respectively. 
Also denote by  $D_1$ and $D_2$ the irreducible components of  $D$  defined by 
$x=0$ and $y=0$ respectively.   
If we identify all these spaces with the germs 
at the origin they define, 
we have an isomorphism 
\begin{equation}\label{lcm} 
 (\hat{Z} ,\hat{S})\cong (A \times Y,A \times D) 
\end{equation}  
with $A \times D_s, s=1,2$,  
corresponding to the germs 
of the two (global) irreducible components of  $\hat{S}$ at $r$. 
The structure of $(\hat{Z} ,\hat{S})$ at a smooth point of  $\hat{F}$  is given 
by the germ at any point outside the origin.  
We conclude that 
$\hat{S}$  is a Cartier divisor in  $\hat{Z} $ 
since $xy$ is not a zero divisor on  $A \times Y$. 
We can then consider the 
logarithmic $1$-forms on $\hat{Z}$ along $\hat{S}$ (cf.\ Section 7).

\ul{Automorphism group}  

We determine the identity component of the automorphism group of $(\hat{Z},\hat{S})$. 
We first recall the following: 

\begin{Lemma}\label{autt}
Let  $Z$  be a Joyce twistor space associated to a $K$-invariant self-dual structures 
on $m\B^2$ with $m\geq 1$ and  $S:=S^+_i+S^-_i+S^+_j+S^-_j$  as before.  
Then $Aut_0(Z,S)\cong  G:=\C^{*2}$. 
For the blowing-up $\tilde{Z}$ of $Z$, we have a natural isomorphism 
$Aut_0(\tilde{Z},\tilde{S}\cup \tilde{Q}) \cong Aut_0(Z,S)$, 
where $\tilde{Q}=Q_i\cup Q_j$. 
\end{Lemma} 

{\em Proof}.  
See e.g.\ \cite[Proposition 5.5]{fjc} for the first assertion when $m>1$. 
A direct computation yields also the result when $m=1$, the details being omitted. 
Since the center of the blowing up $\mu$ is $G$-invariant, 
the second assertion is obvious.   \hfill q.e.d.  

\vspace{3 mm} 
Using the above lemma we shall 
show the corresponding result for the pair $(\hat{Z},\hat{S})$: 
\begin{Proposition}\label{au}
We have 
$Aut_0(\hat{Z},\hat{S})=\{e\}$ in Case-H and $\cong \C^*$ in Case-P. 
\end{Proposition} 

First we recall that 
$\varphi : Q_i \ra  Q_j$   induces isomorphisms 
$\varphi ^\pm _i: H^{\pm }_i  \ra  E^{\pm }_i$ and 
$\varphi ^\pm _j: E^{\pm }_j  \ra  H^{\pm }_j$.  Then, with respect to to the 
natural $G$-equivariant isomorphisms 
$Q_l \cong  H^+_l\times E^{+}_{l'}$, $l=i,j$,  
we may write $\varphi =\varphi ^+_i\times \varphi ^+_j$. 

Hence, given one-parameter subgroup $\rho : \C^* \ra  G$ 
with the induced $\C^*$-actions on  $(Q_l,H^+_l, E^+_l), l=i,j$, 
the following two conditions are equivalent: 

1) $\varphi:Q_i\ra Q_j$  is $\rho$-equivariant, and 

2) $\varphi^+_l, l=i,j$, are both $\rho$-equivariant. 
 
In this case the one-parameter subgroup 
corresponding to the curve  $H^+_l$ in $\tilde{S}^+_l$ 
is the one-parameter subgroup $\mu ^+_l:=-\rho _l+\rho _{l+1} $ 
corresponding to $L_l$, where  
$\rho_d$ is the one-parameter subgroup corresponding to  $\C_d$ 
(cf.\ \cite[Prop.6.12 and p.241(10)]{fjj}).  
On the other hand, 
the one-parameter subgroup $\nu ^+_l$ corresponding to $E^+_l$ is given by 
 $\nu ^+_l=\mp(-\rho _{l}-\rho _{l+1}) $, 
where we take  $-$-sign (resp.\ $+$-sign) for 
$l= i$ (resp.\ $j$) (cf.\ \cite[p.235(5)]{fjj}).  
Hence by Lemma \ref{ps} 
the one-parameter subgroup 
which makes  
$\varphi ^+_{i}$ (resp.\ $\varphi ^+_{j}$ ) 
equivariant is 
\[ -\rho _i+\rho _{i+1}-\rho _j-\rho _{j+1} \ \ \mbox{(resp.}   
-\rho _i+\rho _{i+1}+\rho _j+\rho _{j+1}) \] 
up to signs.  
(The assumption $k>2$ made in Lemma \ref{ps} corresponds to the condition $m>0$ here.) 
Hence the equivariancy of $\varphi $ is given by the coincidence of these 
two subgroups.   Namely  $\rho _j=\rho _{i+1}$ or $\rho _i=-\rho _{j+1}$.  
Since $i<j$, this implies that 
$j=i+1$ or $i=k+j+1$ and the latter holds only when  $j=k$ and $i=1$. 
Namely, in the cyclic sense we have  $j=i+1, 1\leq i\leq k$.  

From this we get the following: 

\begin{Lemma}\label{ij}
Let  $G_1$  be the maximal connected subgroup of  $G$  
such that  $\varphi $  is $G_1$-equivariant. 
In Case-H, $G_1$ reduces to the identity, and in Case-P, 
$G_1 \cong  \C^*$.  
\end{Lemma}

\vspace{3 mm} 
{\em Proof of Proposition \ref{au}}.  
Since  $\nu$ is the normalization, we have the natural inclusion 
$Aut_0(\hat{Z},\hat{S})\subseteq Aut_0(\tilde{Z},\tilde{S}\cup \tilde{Q})$, and the latter is isomorphic to $G$ by Lemma \ref{autt}. 
On the other hand, with respect to this inclusion 
an element  $g \in G$  is contained in  $Aut_0(\hat{Z},\hat{S})$ 
if and only if  $g$  commutes with  $\varphi$.  
Thus the proposition follows from Lemma \ref{ij}. 
\hfill q.e.d.  

\vspace{3 mm} 
Proposition \ref{au} has the following implication. 

\begin{Proposition}\label{isoa}
In Case-H  the isomorphism class of $(\hat{Z},\hat{S})$ 
is independent of the choice of  $\varphi$. 
In Case-P there exists a one-parameter family of 
isomorphisms $\varphi_t: (Q_1,F_1)\ra (Q_2,F_2), t\in \C^*$,   
such that the corresponding pairs $(\hat{Z}_t,\hat{S}_t)$ form a 
non-trivial family of log-deformations of $(\hat{Z},\hat{S})$ 
and exhausts all non-isomorphic pairs obtained from some $\varphi$. 
\end{Proposition} 

{\em Proof}.  
Let $I=Isom((Q_1,F_1),(Q_2,F_2))$  be the space of 
isomorphisms of  $(Q_1,F_1)$ to $(Q_2,F_2)$.  $I$  has a natural structure 
of an algebraic 
principal homogeneous space of $G$. 
Then with respect to the algebraic action of  $G$  on $I$ 
defined by  $\psi  \ra g\psi g^{-1}, g\in G$, 
the identity component of the stabilizer group at $\varphi $ 
is precisely identified with the algebraic subgroup $G_1$ of $G$ in Lemma \ref{ij}.  
Thus 
by Proposition \ref{au} the first assertion is 
immediate from this since 
the $G$-action on $I$ above is then transitive. 

Similarly, in Case-P the $G$-action on  $I$ 
is not transitive and admits a one dimensional quotient  
isomorphic to $\C^*$.  Then any closed one dimensional subspace of  $I$  
which is mapped surjectively to this quotient parametrizes the 
family of isomorphisms  $\varphi$ with the properties of the 
proposition.    \hfill q.e.d.  

\section{Statement of main theorems}  \label{ldt}
 
\vspace{3 mm} 

We retain the notations of the previous sections.  
Moreover, we denote by 
$U$  the smooth locus  $\hat{Z}_{reg}=\hat{Z}-\hat{Q}$ of  $\hat{Z}$. 
Let  $\hat{B}$  be the set of tangential points  $\hat{p}^\pm_n$ 
and set  $V=U-\hat{B}$.  
Thus  $\hat{S}$ is a divisor with normal crossings on  $V$.   

The structure of the sheaf $\Omega_{\hat{Z}}(\log \hat{S})$ (cf.\ \S 2) 
along the singular locus  $\hat{Q}$  of  $\hat{Z}$ can be read from its structure 
at any of the singular points  $r$  of $\hat{F}$.  
In the notations of (\ref{lcm}) 
let  $p: A	\times Y \ra  A$  and $q: A	\times Y \ra  Y$  be 
the natural projections. Then we have 
\begin{equation}\label{w} 
\Omega _{\hat{Z}} \cong p^*\Omega_A\oplus q^*\Omega_Y \ \mbox{and}\  
\Omega _{\hat{Z}}(\log S) \cong p^*\Omega_A\oplus q^*\Omega_Y(\log D). 
\end{equation}  
such that the natural inclusion 
$\iota_{\hat{Z}} : \Omega _{\hat{Z}}\ra \Omega _{\hat{Z}}(\log \hat{S})$ 
is given locally by  
\begin{equation}\label{lp} 
id_A \oplus \iota_Y : 
p^*\Omega_A\oplus q^*\Omega_Y \ra 
p^*\Omega_A\oplus q^*\Omega_Y(\log D), 
\end{equation}  
where  $id_A$  denotes the identity of $A$  
and $\iota_Y$ the natural inclusion on  $Y$. 
\begin{Proposition}\label{lg} 
1)  $\Omega_{\hat{Z}}(\log \hat{S})$ and $\Theta_{\hat{Z}}(-\log \hat{S})$ are 
locally free on  $V$  and 
reflexive on $U$.  In particular both sheaves have homological codimension $\geq 2$ 
on $U$. 

2) $\Theta_{\hat{Z}}(-\log \hat{S})$ is isomorphic to the dual module 
of $\Omega _{\hat{Z}}(\log \hat{S}) $ on the whole  $\hat{Z}$.  

3)  There exists an exact sequence of  $O_{\hat{Z}}$-modules 
\begin{equation} \label{p} 
0 \ra \Omega_{\hat{Z}} \ra \Omega_{\hat{Z}}(\log \hat{S}) \stackrel{b}{\ra} 
\oplus_{\{\pm ,l=i,j\}}O_{\hat{S}^\pm_l} \ra 0 
\end{equation}  
where  $b$ is the (Poincare) residue homomorphism.  
\end{Proposition} 

{\em Proof}.  {\em On}  $U$ 
the assertions  are all 
special cases of the results due to K.\ Saito \cite{sai} 
(cf.\ (1.6), (1.7) and (2.9) of \cite{sai}). 
In particular the assertion 1) is true (cf.\ \cite[(1.21) Corollary]{st}). 
Moreover, he showed that 
there exists a natural perfect pairing on $U$ 
\[\al _U : \Omega_{U}(\log \hat{S}) 
\times \Theta_{U}(-\log \hat{S})  \ra O_U \] 
making  $\Omega_{U}(\log \hat{S}) $  and  
$\Theta_{U}(-\log \hat{S}) $  the dual $O_U$-modules of each other. 

We shall show that the pairing $\al _U $ and the exact sequence (\ref{p}) both 
extend to the whole  $\hat{Z}$.  
First we note the following two properties of 
$\Theta_{\hat{Z}}(-\log \hat{S}) $ on the whole $\hat{Z}$: 

a) There exist no local sections of  $\Theta_{\hat{Z}}(-\log \hat{S}) $  whose support is dimension $\leq 2$. 

b) Any local section of  $\Theta_{\hat{Z}}(-\log \hat{S}) $  defined outside an analytic subset, say $J$,  
of codimension $\geq 2$ extends holomorphically across  $J$.  

In fact, e.g. b) follows from the fact that 
in the notation of the above definition  
the quotient $v(f)/f$, which is holomorphic outside $J$, extends 
to a holomorphic function across $J$.   

Now as for the extension of $\al _U$, 
$\al _U$  extends trivially 
to $\al _W$ on $W:= \hat{Z}-\hat{F} = U\cup (\hat{Q}-\hat{F})$ 
since at the points of  $\hat{Q}-\hat{F}$,  
$\Omega_{\hat{Z}}(\log \hat{S}) =\Omega_{\hat{Z}}$, 
$\Theta_{\hat{Z}}(-\log \hat{S}) =\Theta_{\hat{Z}}$ and 
$\Theta_{\hat{Z}}$ is the dual of $\Omega_{\hat{Z}}$.  
Since  $\hat{F}$ is of codimension $\geq 2$  in $\hat{Z}$, 
the weak normality of  $\hat{Z}$ implies 
that $\al _W$ further extends to yield a natural pairing 
\[\al _{\hat{Z}}  : \Omega_{\hat{Z}}(\log \hat{S}) 
\times \Theta_{\hat{Z}}(-\log \hat{S})  \ra O_{\hat{Z}} \] 
on the whole  $\hat{Z}$.  
Moreover, by the above properties of  $\Theta_{\hat{Z}}(-\log \hat{S}) $  
we see that the induced map  
$\Theta_{\hat{Z}}(-\log \hat{S}) \ra \Omega_{\hat{Z}}(\log \hat{S})^*$  
is isomorphic since it is already isomorphic on $W$.  
In particular the assertion 2) is proved. 

Finally, from the local description of the inclusion $\iota_{\hat{Z}} $ 
(\ref{lp}) and 
from the standard Poincare residue exact sequence  
\begin{equation}\label{pc} 
0 \ra \Omega_Y \ra \Omega_Y(\log D) \stackrel{b_Y}{\ra} \oplus_{s=1,2} O_{D_s} \ra 0 
\end{equation}                                  
for $(Y,D)$ 
we readily obtain the exact sequence (\ref{p}) extending the one obtained outside 
$\hat{Q}$  by \cite{sai}.  \hfill q.e.d.  

\vspace{3 mm} 
	
\ul{Log-deformations of $(\hat{Z},\hat{S})$} 

We consider the log-deformations of the pair $(\hat{Z},\hat{S})$.  
This amounts to considering deformations 
of the pair $(\hat{Z},\hat{S})$  
which induce deformations of each irreducible components 
$\hat{S}^\pm_l$ of  $\hat{S}$.  
Let  
\begin{equation}\label{kr'} 
 g: ({\mathcal Z} ,{\mathcal S} ) \ra T,\ (Z_o,S_o)=(\hat{Z},\hat{S}),\ o\in T,
\end{equation}  
be the Kuranishi family of log-deformations of the pair  $(\hat{Z},\hat{S})$.  
For any  $t\in T$,  $Z_t$ and  $S_t$ shall denote respectively the fibers 
over  $t$  of the projections ${\mathcal Z}  \ra  T$  and ${\mathcal S} \ra  T$. 
$S_t$ consists of four irreducible components $S^\pm_{l,t}$ which are 
deformations of  $\hat{S}^\pm_l$  respectively.  
For a fixed $l$, $S^\pm_{l,t}$ are mutually 
disjoint since this is true at $t=o$. 
Similarly consider the Kuranishi family of deformations 
of the pair  $(\hat{Z},\hat{S})$ which are locally trivial at each point 
of  $\hat{Z}$ (resp.\ at each point of  $\hat{Q}$, resp.\ 
at each tangential point $p=\hat{p}^\pm_n$ of $U=\hat{Z}_{\rm{reg}}$.) 
These are subfamilies $h$ (resp.\ $h_{\hat{Q}}$, resp.\ $g_p$) 
of  $g$  for a unique subspace  $A$  ((resp.\ $A(\hat{Q})$, resp.\ $T(p)$) of $T$.  
Clearly we have  $A=A(\hat{Q})\cap (\cap_{p\in \hat{B} }T(p))$.  

First we note the following: 
\begin{Theorem}\label{ii}
Let $g: ({\mathcal Z} ,{\mathcal S} ) \ra T$ be the Kuranishi family of 
$(\hat{Z},\hat{S})$ as above.  Then the following hold: 

1) $T$  is smooth of dimension $3m$ in Case-H 
{\em (}resp.\ $3m+1$ in Case-P{\em )}.  

2) $A(\hat{Q})$ and $T(p)$ are smooth hypersufaces of $T$  passing through $o$ such that 
$D:=A(\hat{Q})\cup (\cup_{p\in \hat{B} }T(p))$ is a divisor with normal crossings in  $T$.  
In particular $I:= \cap_{p\in \hat{B} }T(p)$ is a smooth subspace 
of $T$ of dimension $m$ in Case-H (resp.\ $m+1$ in Case-P) 
and $A$ is a smooth hypersuface of $I$. 
\end{Theorem}

We set $C^\pm_{l,t}:= S^\pm_{l,t}\cap (S^+_{l',t}\cup S^-_{l',t})$ with $\{l,l'\} = \{i,j\}$. 
As for the structure of the surfaces  $S^\pm_{l,t},\ l=i,j$, for $t\in T-A(\hat{Q})$ 
we shall show the following: 

\begin{Theorem}\label{iit}
1) Assume that  $t \in T-A(\hat{Q})$.  
Then the fibers $Z_t$ and $S^\pm_{l,t}$ are all smooth, and 
$S^\pm_{l,t}$ are surfaces of class VII.  In Case-H {\em (}resp.\ Case-P{\em )} their 
minimal models $\bar{S}^\pm_{l,t}$ are 
either a hyperbolic or parabolic {\em (}resp.\ a parabolic{\em )} Inoue surface or 
a diagonal Hopf surface.  Each 
$S^\pm_{l,t}$ is obtained from $\bar{S}^\pm_{l,t}$ by blowing-up, 
as described in Lemma \ref{atc},  
a finite number of {\em (}possibly infinitely near{\em )} points 
on the image $\bar{C}^\pm_{l,t}$ of  $C^\pm_{l,t}$ in $\bar{S}^\pm_{l,t}$. 

2) Assume that $t \in I-A(\hat{Q})$.  
Then in Case-H {\em (}resp.\ Case-P{\em )}, $S^\pm_{l,t}$ are all 
properly blown-up hyperbolic {\em (}resp.\ parabolic{\em )} Inoue surfaces. 
In Case-H 
the isomorphism class of $S^\pm_{l,t}$ is independent of  $t$, 
$S^+_{l,t}$ and $S^-_{l,t}$ are isomorphic to each other, and 
$S^\pm_{i,t}$ and  $S^\pm_{j,t}$ are transpositions of each other.  
If $(i,j)$ is minimal, 
they are hyperbolic {\em (}resp.\ parabolic{\em )} Inoue surfaces.  
If $t \in T-D$, $S^\pm_{l,t}$ are blown-up diagonal Hopf surfaces.  

3) In Case-H  the Kuranishi family $g$ is universal.  
\end{Theorem} 

The proofs of Theorem \ref{ii} and Theorem \ref{iit} will be 
given in Section \ref{coh}. 

\vspace{3 mm} 
\ul{Real deformations and twistor spaces}

In the construction above, 
suppose that we have taken $\varphi : Q_i \ra Q_j$  
to be $\sigma $-equivariant, which is always possible. 
Then 
$(\hat{Z},\hat{S})$ has the induced real structure 
(denoted by the same letter $\sigma$) 
which interchanges $\hat{S}^+_l$ and $\hat{S}^-_l, l=i,j$. 
 
First we consider Case-H, namely we assume that $|j-i|>1$.  
Then by Theorem \ref{iit} 
$g$ is universal and the real structure $\sigma $ on $\hat{Z}=Z_o$ 
extends canonically to the family  $g: ({\mathcal Z} ,{\mathcal S} ) \ra  T$.  
Denote by $T^{\sigma} $ the set of fixed points of $\sigma $, 
which is a real submanifold of  $T$ of real dimension  $3m$.  
It is not contained in any proper analytic subset of  $T$.  
For any point $t\in T^{\sigma} $, the fiber $(Z_t,S_t)$ has 
the induced real structure $\sigma _t$.  Recall that we set 
$M[m]=(S^1\times S^3)\# m\bar{\B}^2$.  

\begin{Theorem}\label{tw}
For any $t\in T^{\sigma} -A(\hat{Q})$, 
the fiber  $Z_t$,  together with the induced real structure 
$\sigma_t$,  is 
a twistor space of an anti-self-dual bihermitian structure $([g]_t,J_{1,t},J_{2,t})$ 
on the smooth oriented manifold  $M[m]$ such that 
$(M[m],\pm J_{1,t})\cong S^\pm_{i,t}$ and 
$(M[m],\pm J_{2,t})\cong S^\pm_{j,t}$.  
The structure of the surfaces $S^\pm_{l,t}$ are given by Theorem \ref{iit}, 1).  
Moreover, this family gives a universal family of 
anti-self-dual bihermitian structures on  $M[m]$ at each point of 
$t\in T^{\sigma} -A(\hat{Q})$ with $3m$ real parameters. 
\end{Theorem} 

{\em Proof}.  The fact that  $Z_t$  is a twistor space associated to 
a self-dual structure on the $C^{\infty}$ $4$-manifold 
$(S^1\times S^3)\#m\B^2$, or equivalently, to an anti-self-dual structure on 
$M=M[m]$, is similar to \cite[4.2]{df} 
except that in this case 
$M$ is a ``self-connected sum'' of $ m\B^2$  
as we have identified $Q_i$ and $Q_j$ in a single manifold  $\tilde{Z}$.   
Here the degrees of the surfaces  $S^\pm_{l,t}$ are 
all equal to one as well as the original surfaces $S^\pm_l$ in  $Z$.   
Thus 
$\{S^+_{l,t}, S^-_{l,t}\}, l=i,j$, are two pairs of elementary surfaces in  $Z_t$.   
By Lemmas \ref{elm} and \ref{c2}, 
they are in Case 2, and hence give an anti-self-dual bihermitian structure on 
$M[m]$.   The rest follows immediately from Theorem \ref{iit}. 
\hfill q.e.d.  

\vspace{3 mm} 
There exist nice subfamilies of this universal family.  
Most typically, 
when  $t\in (T^{\sigma}\cap I) -A$, we know by Theorem \ref{iit} 
that $S^\pm_{i,t}$ and $S^\pm_{j,t}$ are both 
properly blown-up hyperbolic Inoue surfaces whose isomorphism class is 
independent of  $t$.  
More precisely, we shall show the following: 

\begin{Theorem}\label{ma}
Let  $S$ be an arbitrary properly blown-up hyperbolic Inoue surface. 
Let $m$ be the second Betti number of $S$ and $\bar{m}$ that of its minimal model.   
Then there exist $\bar{m}$ families 
of anti-self-dual bihermitian structures $([g]_t; J_{1,t},J_{2,t})$ 
on  $M[m]$ with real smooth $m$-dimensional parameters $t$ 
such that  $(M[m],\pm J_{1,t})$ and  $(M[m],\pm J_{2,t})$ 
are biholomorphic respectively to $S$ and to its transposition  ${}^{\mbf}S$, 
independently of $t$.   
\end{Theorem}  

{\em Proof of Theorem \ref{ma}}.  By Proposition \ref{m} there exists an admissible toric surface 
$\tilde{S}$ 
such that the given hyperbolic Inoue surface  $S$  is obtained by smoothing 
the rational surface $\hat{S}$ with nodal curve obtained from $\tilde{S}$ 
via the procedure of Section \ref{rswn}.  
On the other hand, let $\bar{S}$ be the surface obtained from 
$\tilde{S}$ by blowing down its $(-1)$-curve $E$.  Then by 
Lemma \ref{tj} there exists a Joyce twistor space $Z$ which contains $\bar{S}$ as 
one of its $G$-invariant elementary surfaces.  We may take  $\bar{S}=S^+_i$.  Then 
we get a unique number  $j$  such that the twistor line $L_j$ passes thruough the point 
$p \in \bar{S}=S^+_i$ which is the image of $E$.  (Changing the numbering cyclically 
we can assume that $i<j$.) 

Now starting from this Joyce twistor space $Z$  and the pair of indices  $(i,j)$ 
we perform the construction of Section \ref{cstr} and 
consider the universal family in Theorem \ref{tw}. We restrict the family 
to  $I\cap T^{\sigma} -A(\hat{Q})$ and obtain a real $m$-dimensional family 
of bihermitian structures on  $M[m]$. By Theorem \ref{iit} 
the corresponding pairs of elementary surfaces are as described in 
the theorem.  Finally, for the given  $S$, according to Proposition \ref{m} we have actually 
$\bar{m}$ choices of admissible toric surfaces $\tilde{S}$ and correspondingly we get 
$\bar{m}$ such families. 
\hfill q.e.d.  

\vspace{3 mm} 

{\em Remark 7.1}. 
The parameter  $t$  in the above theorem belongs to a complement of 
a real hyperplane in  $\R^m$ in a neighborhood of the origin.  In this sense 
the parameter space has actually two connected components.  However, 
our construction depends on the initial Joyce self-dual metrics.  
Once the $K$-action is fixed, they form a 
connected $(m-1)$-dimensional family parametrized by 
$m+2$ points on the real projective line  $\R\B^1$ up to the 
action of $PSL(2,\R)$.  It should still be checked if the global parameter 
space is connected or not. 
On the other hand, the choice of $K$-action and of the index $(i,j)$ gives 
discrete invariants for our construction. 
The $\bar{m}$ families in the theorem refers to $\bar{m}$ families 
with different discrete invariants 
but giving one and the same properly blown-up hyperbolic Inoue surface.  
Basically,  similar remarks apply also to the other theorems. 
(See \cite{fjo}.)    

\vspace{3 mm} 
Next,  
for a surface which is obtained from a properly blown-up hyperbolic Inoue surface 
by a small deformation,   
we can show a similar but weaker resut:  

\begin{Theorem}\label{mad}
Let  $S$ be an arbitrary properly blown-up hyperbolic Inoue surface 
and $C$ the unique anti-canonical curve on it.  
Let $m$ be the second Betti number of  $S$ and $\bar{m}$ that of 
its minimal model.  
Let $(S',C')$ be any fixed sufficiently small 
deformation of  $(S,C)$ in the Kuranishi family (\ref{kr1}).   
Then $S'$ admits $\bar{m}$ $m$-dimensional families of 
anti-self-dual bihermitian structures.  Namely  
there exist $\bar{m}$ families 
of anti-self-dual bihermitian structures $([g]_t; J_{1,t},J_{2,t})$ 
on  $M[m]$ with real and smooth $m$-dimensional parameters $t$ 
such that  $(M[m],J_{1,t})$ is 
biholomorphic to $S'$ independently of $t$.  
\end{Theorem} 

{\em Remark 7.2}.  
As noted in Remark 3.1 
$(S',C')$ is a deformation of $(S,C)$ as an anti-canonical pair 
with $C'$ disconnected as well as $C$.  By \cite[Th.4.1]{pnt6} 
the existence of a disconnected anti-canonical curve 
is a necessary condition for the existence of an anti-self-dual bihermitian structure. 
The above theorem is our strongest result toward the sufficiency of this condition, 
although the converse is in general not true as the diagonal Hopf surface case 
already shows. 

\vspace{3 mm} 

We shall give proofs of this and the next theorem in the next section. 
In Case-P our result is less complete. We state the result only in the minimal case 
for simplicity.  
\begin{Theorem}\label{map} 
For any $m>0$ 
there exists a real one-parameter family 
of parabolic Inoue surfaces with second Betti number $m$ such that 
for any member  $S$  of this family 
we have a family of anti-self-dual bihermitian structures $\{([g]_t; J_{1,t},J_{2,t})\}$ 
on  $M[m]$ with real and smooth $m$-dimensional parameter $t$ 
such that  $(M[m],J_{1,t})$ and $(M[m],J_{2,t})$ are both 
biholomorphic to $S$.  
\end{Theorem}  

\vspace{3 mm} 
{\em Remark 7.3}. 
1) It remains open to identify the parabolic Inoue surfaces 
which corresponds to the points of $I\cap T^{\sigma} -A(\hat{Q})$.  

2) For each fixed $m$ a Joyce twistor space  $Z$  
which produces Case-P in the minimal case is 
uniquely characterized, up to deformations, 
by the three equivalent conditions of Proposition 6.14 of \cite{fjj}.  
(We call these {\em LeBrun-Joyce twistor spaces}.)  
In this case the weight sequence (\ref{0w}) of  $S_0$  is given by 
$(a_1,\ldots ,a_k)=(1,m,1,2,\ldots , 2)$ and then  $(i,j)=(1,2)$ is 
the unique choice of the indices.  

\vspace{3 mm} 
{\em Remark 7.3}. 
One interesting probelm is to compare the above anti-self-dual bihermitian structures 
on parabolic Inoue surfaces with those constructed by LeBrun \cite{lb}.  
For instance we can ask if both coincide at least for some parameters.

\section{Proof of theorems} \label{coh} 

In this section we prove Theorems \ref{ii}, \ref{iit},  
\ref{mad} and \ref{map}. 
The main part of the proof consists 
in showing the following theorem, which is the corresponding 
cohomological computations of relevant $Ext$ and cohomology groups. 

\begin{Theorem}\label{pt} 
 Both $H^2(\hat{Z},\Theta_{\hat{Z}}(-\log \hat{S}))$ and 
$Ext^2_{O_{\hat{Z}}}(\Omega_{\hat{Z}}(\log \hat{S}),O_{\hat{Z}})$ vanish. 
We have a natural short exact sequence 
\begin{equation} \label{shoo} 
0 \ra  H^1(\hat{Z},\Theta_{\hat{Z}}(-\log \hat{S})) 
\ra  Ext^1_{O_{\hat{Z}}}(\Omega_{\hat{Z}}(\log \hat{S}),O_{\hat{Z}}) 
\stackrel{\hat{c}}{\ra} H^0(O_{\hat{Q}})\oplus(\oplus_{p\in \hat{B}} \C_p)\ra 0. 
\end{equation}  
In Case-H and Case-P we have respectively 
\begin{gather}
\dim H^1(\hat{Z},\Theta_{\hat{Z}}(-\log \hat{S})) = m-1, 
\ \ \mbox{and} \ \ \label{sh1}
\dim Ext^1_{O_{\hat{Z}}}(\Omega_{\hat{Z}}(\log \hat{S}),O_{\hat{Z}}) = 3m \\
\dim H^1(\hat{Z},\Theta_{\hat{Z}}(-\log \hat{S})) = m, 
\ \ \mbox{and} \ \  \label{sh2}
\dim Ext^1_{O_{\hat{Z}}}(\Omega_{\hat{Z}}(\log \hat{S}),O_{\hat{Z}}) = 3m+1.  
\end{gather}
Finally, 
\begin{equation}\label{ex0} 
\dim Ext^0_{O_{\hat{Z}}}(\Omega_{\hat{Z}}(\log \hat{S}),O_{\hat{Z}}) = 0 \ \mbox{in Case-H and} \  
=1 \ \mbox{in Case-P}. 
\end{equation}  
\end{Theorem} 

Recall here that $\hat{B}$ is the set of tangential points $p_n^\pm$ on  $\hat{Z}$. 
We start by determining the structure 
of ${\mathcal Ext}^i_{O_{\hat{Z}}}(\Omega_{\hat{Z}}(\log \hat{S}),O_{\hat{Z}})$. 

\begin{Lemma}\label{ext}
1) ${\mathcal Ext}^0_{O_{\hat{Z}}}(\Omega_{\hat{Z}}(\log \hat{S}),O_{\hat{Z}})\cong\Theta_{\hat{Z}}(-\log \hat{S})$ and hence $Ext^0_{O_{\hat{Z}}}(\Omega_{\hat{Z}}(\log \hat{S}),O_{\hat{Z}})\\ \cong 
H^0(\hat{Z},\Theta_{\hat{Z}}(-\log \hat{S}))$.  

2) ${\mathcal Ext} ^1_{O_{\hat{Z}}}(\Omega_{\hat{Z}}(\log \hat{S}),O_{\hat{Z}})\cong 
{\mathcal E}^1_{\hat{Q}} \oplus {\mathcal E} ^1_{\hat{B}}$, where ${\mathcal E}^1_Y$ has 
support in $Y \ (Y=\hat{Q}, \hat{B})$.  Moreover ${\mathcal E} ^1_{\hat{Q}} \cong O_{\hat{Q}}$ and 
${\mathcal E} ^1_{\hat{B}}\cong \oplus_{p \in \hat{B}} \C_p$, 
where $\C_p$ is the skyscraper sheaf at $p$.  

3) ${\mathcal Ext}^2_{O_{\hat{Z}}}(\Omega_{\hat{Z}}(\log \hat{S}),O_{\hat{Z}})=0$. 
\end{Lemma} 

{\em Proof}.  
Since ${\mathcal Ext}^0\cong {\mathcal Hom}^0$, the assertion 1) follows from 
Proposition \ref{lg}, 2).  
Applying ${\mathcal Ext}(-,O_{\hat{Z}})$ to the exact sequence (\ref{p}) 
we obtain a long sheaf exact sequence 
\begin{gather}\label{lon} 
0 \ra  \Theta_{\hat{Z}}(-\log \hat{S}) \ra  \Theta_{\hat{Z}} 
\stackrel{b}{\ra} \oplus_{\{\pm ,l=i,j\}}{\mathcal Ext}^1_{O_{\hat{Z}}}(O_{\hat{S}^\pm_l},O_{\hat{Z}}) \\
\ra  {\mathcal Ext}^1_{O_{\hat{Z}}}(\Omega_{\hat{Z}}(\log \hat{S}),O_{\hat{Z}}) \nonumber 
\stackrel{a}{\ra}  {\mathcal Ext}^1_{O_{\hat{Z}}}(\Omega_{\hat{Z}},O_{\hat{Z}})
\end{gather}  
On the other hand, since $\Omega_{\hat{Z}}(\log \hat{S})$ 
is locally free outside $\hat{Q}\cup \hat{B}$, the support of 
${\mathcal Ext}^1_{O_{\hat{Z}}}(\Omega_{\hat{Z}}(\log \hat{S}),O_{\hat{Z}})$ is contained in 
$\hat{Q}\cup \hat{B}$.  
Locally along  $\hat{Q}$  the map $a$ is induced from the map in (\ref{lp}). 
Since $q^*\Omega_Y(\log D)$ and $q^*\Omega_Y$ are locally free 
and hence their ${\mathcal Ext} ^1$'s 
vanish, $a$ is locally isomorphic to the identity 
${\mathcal Ext}^1_{O_{\hat{Z}}}(p^*\Omega_{A},O_{\hat{Z}}) \ra 
{\mathcal Ext}^1_{O_{\hat{Z}}}(p^*\Omega_A,O_{\hat{Z}})$.  
In particular  $a$  is isomorphic along  $\hat{Q}$. 
On the other hand, by Friedman \cite[Corollary 2.4]{fr} 
${\mathcal Ext}^1_{O_{\hat{Z}}}(\Omega_{\hat{Z}},O_{\hat{Z}}) \cong O_{\hat{Q}}$ along $\hat{Q}$. 
Also we see from (\ref{lp}) that $\Omega_{\hat{Z}}(\log \hat{S})$ and 
$\Omega_{\hat{Z}}$ are locally isomorphic.  Hence 
${\mathcal Ext}^2_{O_{\hat{Z}}}(\Omega_{\hat{Z}}(\log \hat{S}),O_{\hat{Z}})$ vanishes 
since so does ${\mathcal Ext}^2_{O_{\hat{Z}}}(\Omega_{\hat{Z}},O_{\hat{Z}})$ along $\hat{Q}$ 
(cf.\ \cite[Sect.2]{fr}). 

Thus 2) and 3) are shown along  $\hat{Q}$.  It remains to check these assertions 
at each point $p$ of   $\hat{B}$.  First of all, since the homological codimension of 
$\Omega_{\hat{Z}}(\log \hat{S})$ is two at $p$ by Proposition \ref{lg},  
the assertion 3) is true there.  We shall compute  
${\mathcal Ext}^1_{O_{\hat{Z}}}(\Omega_{\hat{Z}}(\log \hat{S}),O_{\hat{Z}})$ 
at $p$. Since  
${\mathcal Ext}^1_{O_{\hat{Z}}}(\Omega_{\hat{Z}},O_{\hat{Z}})=0$ at $p$, 
it is the cokernel of  $b$ in (\ref{lon}).  
 Since $\hat{S}^\pm_l$ are Cartier divisors, we have 
 ${\mathcal Ext}^1_{O_{\hat{Z}}}(O_{\hat{S}^\pm_l},O_{\hat{Z}}))\cong N^\pm_l$, 
 where 
 $N^\pm_l=[\hat{S}^\pm_l]|\hat{S}^\pm_l$ is the normal bundle of 
 $\hat{S}^\pm_l$ in  $\hat{Z}$.  

Now we work in the local model (\ref{lt}) so that we may put 
$\hat{Z} = \C^3(x,y,z)$.  
Let  $D_m, m=1,2$, be the irreducible components of  $\hat{S}$ at $p$  defined by 
the local equations $f_1:=z=0$ and $f_2:=z-xy=0$ respectively.  
Let  $N_m=Hom(I_m,O_{D_m})$ be the normal sheaves of  $D_m$  in  $\hat{Z}$, 
where $I_m$  is the ideal sheaf of $D_m$. 
Then  $b$  is given 
by $b(\theta )=(\theta (f_1)|D_1,\theta (f_2)|D_2)$ in terms of the above 
identifications.  
Then for $\theta =\partial /\partial x, \partial /\partial y, \partial /\partial z$ 
we obtain 
$(0,-y), (0,-x), (1,1)$ restricted to $(D_1,D_2)$.  We conclude immediately 
that the cokernel of  $b$, which has support in $p$, is in fact one dimensional.  
This shows 2) at $p$  and the lemma.   \hfill q.e.d.  

\vspace{3 mm} 
{\em Remark 8.1}.  For $p\in \hat{B}$ 
${\mathcal Ext}^1_{O_{\hat{Z}}}(\Omega_{\hat{Z}}(\log \hat{S}),O_{\hat{Z}})_p (\cong \C$) is considered to be the tangent space of the local versal log-deformation 
of the pair $(\hat{Z},\hat{S})$ considered as a germ at  $p$.  
The versal deformation $(\hat{Z}_t,\hat{S}_t) =(\C^3,\hat{S}_t)$ is given explicitly 
by the defining equation $z(z-xy+t) = 0$ of $\hat{S}_t$.  For $t\neq 0$, the intersection of 
the two irreducible components $D_{1,t}$ and $D_{2,t}$ of $\hat{S}_t$ is 
now smooth, along which $D_{m,t}$ are transversal.  
We shall denote this versal family by  
$g(p): ({\mathcal Z(p)} ,{\mathcal S(p)} ) \ra T(p),\ 
(Z(p)_o,S(p)_o)=(\hat{Z}(p),\hat{S}(p)),\ o\in T(p)$.   
The original Kuranishi family induces 
a log-deformation of the germ $(\hat{Z}(p),\hat{S}(p))$ and 
we have a versal map $\tau(p): T \ra T(p)$.  
 
\vspace{3 mm} 
We next compute the cohomology groups  $H^i(\Theta_{\hat{Z}}(-\log \hat{S}))$ 
by relating them to the corresponding cohomology groups of the blown-up twistor space 
$(\tilde{Z},\tilde{S})$ and of the original Joyce twistor space  $(Z,S)$, where 
$\tilde{S}=\tilde{S}_i\cup \tilde{S}_j$ and $S=S_i\cup S_j$.  
We first record the infinitesimal form 
of the results of Lemma \ref{autt} and Proposition \ref{au}. 

\begin{Proposition}\label{auti} We have 
$h^0(\Theta_{Z}(-\log S))=2,\ 
h^0(\Theta_{\tilde{Z}}(-\log (\tilde{S}+\tilde{Q})))=2$, \ and 
$h^0(\Theta_{\hat{Z}}(-\log \hat{S}))=0$ in Case-H and $=1$ in Case-P. 
\end{Proposition} 
 
Now we compare the cohomology groups of $(\hat{Z},\hat{S})$ with those of 
$(\tilde{Z},\tilde{S})$ via the normalization exact sequence 
\[ 0 \ra \Theta_{\hat{Z}}(-\log \hat{S}) \ra \Theta_{\tilde{Z}}(-\log (\tilde{S}+\tilde{Q})) \ra \Theta_{\hat{Q}}(-\log \hat{F}) \ra 0\] 
with associated long exact sequence 
\begin{gather}
0 \ra H^0(\Theta_{\hat{Z}}(-\log \hat{S}))\ra H^0(\Theta_{\tilde{Z}}(-\log (\tilde{S}+\tilde{Q}))) \stackrel{a}{\ra} H^0(\Theta_{\hat{Q}}(-\log \hat{F})) \ra \nonumber \\
\ra H^1(\Theta_{\hat{Z}}(-\log \hat{S})) \ra H^1(\Theta_{\tilde{Z}}(-\log (\tilde{S}+\tilde{Q}))) \ra H^1(\Theta_{\hat{Q}}(-\log \hat{F})) \ra  \\
\ra H^2(\Theta_{\hat{Z}}(-\log \hat{S})) \ra H^2(\Theta_{\tilde{Z}}(-\log (\tilde{S}+\tilde{Q}))) \ra H^2(\Theta_{\hat{Q}}(-\log \hat{F})) \ra  \nonumber 
\end{gather}
Here $\hat{F}$ is the anti-canonical cycle of the toric surface  $\hat{Q}$  and 
hence we have $\Theta_{\hat{Q}}(-\log \hat{F}) = O_{\hat{Q}}^2$; 
thus  $H^i(\Theta_{\hat{Q}}(-\log \hat{F}))=0$ for $i=1,2$ 
and $H^0(\Theta_{\hat{Q}}(-\log \hat{F}))=\C^2$. 
In view of Proposition \ref{auti} this implies that 
$a$  is isomorphic in Case-H 
and has one dimensional image in Case-P. 
Thus the above exact sequence reduces in Case-H to: 
\begin{gather}
 H^0(\Theta_{\tilde{Z}}(-\log (\tilde{S}+\tilde{Q}))) \cong H^0(\Theta_{\hat{Q}}(-\log \hat{F})) \cong \C^2 \label{h1} \\ 
 H^1(\Theta_{\hat{Z}}(-\log \hat{S})) \cong H^1(\Theta_{\tilde{Z}}(-\log (\tilde{S}+\tilde{Q}))) 
\end{gather}
and in Case-P to the two short exact sequences 
\begin{gather}
0  \ra \C \ra H^0(\Theta_{\tilde{Z}}(-\log (\tilde{S}+\tilde{Q}))) \stackrel{a}{\ra} \C \ra 0 \\ 
0  \ra \C \ra H^1(\Theta_{\hat{Z}}(-\log \hat{S})) \ra H^1(\Theta_{\tilde{Z}}(-\log (\tilde{S}+\tilde{Q}))) \ra 0 \label{p2} 
\end{gather}
In both cases we have 
\begin{equation}\label{iso^} 
H^2(\Theta_{\hat{Z}}(-\log \hat{S})) \cong H^2(\Theta_{\tilde{Z}}(-\log (\tilde{S}+\tilde{Q})))  
\end{equation}  

Next we compare the cohomology groups 
$H^i(\Theta_{\tilde{Z}}(-\log (\tilde{S}+\tilde{Q})))$ with those of $(Z,S)$. 
Namely we haveF
\begin{Lemma}\label{sita}
We get natural isomorphisms 
\begin{equation}\label{iso} 
 H^q(\Theta_{\tilde{Z}}(-\log (\tilde{S}+\tilde{Q}))) \cong H^q(Z, \Theta_Z(-\log S)), \ q\geq 0 .
 \end{equation}  
\end{Lemma} 
{\em Proof}.  
There exists a natural sheaf isomorphism 
 \[ \Theta_{\tilde{Z}}(-\log (\tilde{S}+\tilde{Q})) \cong \mu ^*\Theta_Z(-\log S) \] 
induced by  $\mu_*$.  Together with the Leray spectral sequence for 
\[ E^{p,q}_2:=H^p(R^q\mu_*\mu ^*\Theta_Z(-\log S)) \Rightarrow  H^{p+q}(\mu^*\Theta_Z(-\log S))\]  
and the fact that $R^q\mu_*O=0$ for $q>0$, we get the desired isomorphisms. 
(Note that in a neighborhood of $L_i$ and $L_j$ the sheaf  
$ \Omega_{Z}(\log S)$ is locally free so that the projection formula 
\[ R^q\mu_*\mu ^*\Theta_Z(-\log S)\cong R^q\mu_*O_{\tilde{Z}}\otimes \Theta_Z(-\log S) \] holds.) \hfill q.e.d.  

\vspace{3 mm} 
We have thus reduced our computation to that of the cohomology groups on  $Z$.  
Consider now the short exact sequence of $O_Z$-modules 
\[ 0 \ra \Theta_Z(- S) \ra \Theta_Z(-\log S)\ra  \Theta_S \ra 0\] 
and the associated cohomology exact sequence  
\begin{gather}
0 \ra H^0(\Theta_{Z}(-S)) \ra H^0(\Theta_Z(-\log S)) \ra H^0(\Theta_S) \ra \nonumber \\
  \ra H^1(\Theta_{Z}(-S)) \ra H^1(\Theta_Z(-\log S)) \ra H^1(\Theta_S) \ra \label{s} \\
  \ra H^2(\Theta_{Z}(-S)) \ra H^2(\Theta_Z(-\log S)) \ra H^2(\Theta_S) \ra \nonumber 
\end{gather}
We compute the dimensions of the spaces as follows.  
\begin{Lemma}\label{sll} We have 
\begin{gather}\label{s1}
  h^0(\Theta_{Z}(-S))=0,\ h^1(\Theta_{Z}(-S))=0,\ h^2(\Theta_{Z}(-S))=m+1,\  
 h^3(\Theta_{Z}(-S))=0.  \\
 h^0(\Theta_Z(-\log S))=2,\  h^3(\Theta_Z(-\log S))=0. \label{s2} \\
 h^0(\Theta_S)=2,\ 
 h^2(\Theta_S)=0,\ h^3(\Theta _S)=0. 
\label{s3}
\end{gather}
\end{Lemma} 

{\em Proof}.  
Since $K=-S$, 
we have  $h^i(\Theta_{Z}(-S))
=h^{3-i}(\Omega_{Z})=h^{1,3-i}(Z)$, 
where $h^{p,q}$ denotes the Hodge numbers. 
We know that $h^{k,0}(Z)=0$ for any twistor space and hence also $h^{0,k}(Z)=0$ 
since $Z$  is Moishezon. 
By the same reason we have $b_k(Z)=\sum_{p+q=k} h^{p,q}(Z)$ 
for the Betti numbers $b_k(Z)=b_k(M)+b_{k-2}(M)$, where $M=m\B^2$. 
In particular for odd $k$, we have $b_k(Z)=0$.  Thus 
$h^i(\Theta_{Z}(-S))=0$ if $i$  is odd.  Also from  $h^{1,1}=b_2=m+1$, we have 
$h^2(\Theta_{Z}(-S))=m+1$. 

Further, we show that $h^0(\Theta_{Z}(-S))=0$. 
Take any twistor line  $L$  and consider the standard short exact sequence 
\[     0 \ra  \Theta _L \ra  \Theta _Z|L \ra  N \ra  0  \]
where  $N\cong O(1)\oplus O(1)$  is the normal bundle of  $L$  in  $Z$. 
Since $K|L$ is of degree $-4$, we get that $h^0(\Theta _L(-S))= h^0(N(-S))=0$.  
Hence the above exact sequence tensored with $K$  yields 
$h^0(L, \Theta _Z(-S)|L)=0$, 
from which follows the desired vanishing since  $L$  is arbitrary. Thus (\ref{s1}) 
is proved. 

 The identity component of the automorphism group of  $Z$  is $\C^{*2}$ for $m\geq 2$ 
and it preserves  $S$.   Also in case $m=1, \C^{*2}$ is 
the maximal connected automorphism group of $(Z,S)$. 
Thus $h^0(\Theta_Z(-\log S))=2$.
We get $h^3(\Theta _S)=0$ since $\dim S =2$. Then 
together with (\ref{s1}) 
we deduce $h^3(\Theta_Z(-\log S))=0$ from the exact sequence (\ref{s}).  
This shows (\ref{s2}). 

 We show that 
$h^2(\Theta_S)=0$.  Let $\omega _S:=K[S]|S=O_S$ be the dualizing sheaf of  $S$. 
By Serre duality $H^2(\Theta_S)$ is dual to 
$Ext^0(\Theta_S,\omega _S)\cong Ext^0(\Theta_S,O_S)\cong H^0(\Omega_S^{**})$, 
where $\Omega_S^{**}$ denotes 
the double dual of $\Omega _S $. (See \cite[Lemma (2.9)]{fr} for the structure 
of  $\Omega_S^{**}$ outside tangential points.) 
The last space injects into 
$\oplus_{\pm,l} H^0(\Omega _{S^\pm_l})$ 
which vanishes.  
Hence $h^2(\Theta_S)=0$.  From the exact sequence (\ref{s}) 
together with (\ref{s1}) and (\ref{s2}) we get  $h^0(\Theta_S)=2$.  
(\ref{s3}) is proved.  \hfill q.e.d.  

\vspace{3 mm} 
The main part of the exact sequence (\ref{s}) now reduces to 
\begin{equation}\label{s'} 
0   \ra H^1(\Theta_Z(-\log S)) \ra H^1(\Theta_S) 
  \stackrel{\delta}\ra H^2(\Theta_{Z}(-S)) \ra H^2(\Theta_Z(-\log S)) \ra 0. 
\end{equation}  

\begin{Lemma}\label{ts}
$h^1(\Theta_S)=2m$. 
\end{Lemma} 

{\em Proof}.  
We consider the normalization exact sequence 
\begin{equation}\label{tg} 
0 \ra  \Theta_S \ra \oplus_{l,\pm }\Theta_{S^\pm _l}(-\log B^\pm _l) 
\stackrel{a}{\ra}  \oplus_\al  \Theta_{B_\al}((0+\infty)) \ra  0 
\end{equation}  
where $B^\pm _l$ is the anti-canonical cycle of $S^\pm _l$, 
and 
$B_\al$ are the irreducible components of the curve  $B:=C\cup L_i\cup L_j$;  
$0=0_\al$ and $\infty=\infty_\al$ are the two points of intersection 
of $B_\al$ and  the other irreducible components of $B$ (cf.\ (\ref{clc})). 
(The surjectivity at the tangential points of $a$ may be shown by using 
two sections 
$z\partial /\partial z-x\partial /\partial x, 
z\partial /\partial z-y\partial /\partial y$ 
of $\Theta_Z(-\log S)$ in the local notations of (\ref{lt}).) 

Note that  $B$  has $2k+2$ irreducible components. 
Let 
\begin{gather}
0 \ra H^0(\Theta_S) \ra \oplus_{l,\pm }H^0(\Theta_{S^\pm _l}(-\log B^\pm _l)) 
\ra \oplus_\al  H^0(\Theta_{B_\al}(0+\infty)) \ra \nonumber \\
  \ra H^1(\Theta_S) \ra \oplus_{l,\pm }H^1(\Theta_{S^\pm _l}(-\log B^\pm _l)) 
\oplus_\al  H^1(\Theta_{B_\al}(0+\infty)) \\ 
  \ra H^2(\Theta_S) \ra \oplus_{l,\pm }H^2(\Theta_{S^\pm _l}(-\log B^\pm _l)) 
\oplus_\al  H^2(\Theta_{B_\al}(0+\infty)) \nonumber \ra 
\end{gather}
be the associated long exact sequence. 

Since  $\Theta_{S^\pm _l}(-\log B^\pm _l))\cong O_{S^\pm _l}^2$, 
we have 
\[ h^0(\Theta_{S^\pm _l}(-\log B^\pm _l))=2\ \  \mbox{and}\ \  
h^i(\Theta_{S^\pm _l}(-\log B^\pm _l))=0, i>0. \] 
Similarly, we have 
$\Theta_{B_\al}((0_+\infty))\cong O_{B_\al}$ and hence 
\[ h^0(\Theta_{B_\al}(0+\infty))=1 \ \  \mbox{and}\ \  
h^i(\Theta_{B_\al}(0+\infty)), i>0. \] 
Thus we get $H^2(\Theta_S)=0 $ (deduced above by a different method), and 
the exact sequence 
\begin{gather}
0 \ra H^0(\Theta_S) \ra \oplus_{l,\pm }H^0(\Theta_{S^\pm _l}(-\log B^\pm _l)) 
\ra \oplus_\al  H^0(\Theta_{B_\al}(0+\infty)) 
  \ra H^1(\Theta_S) \ra 0
\end{gather}
with $\oplus_{l,\pm }H^0(\Theta_{S^\pm _l}(-\log B^\pm _l))\cong \C^8$ and 
$\oplus_\al  H^0(\Theta_{B_\al}(0+\infty)) \cong \C^{2k+2}$. 
Together with (\ref{s3}) we get  $h^1(\Theta_S)=2m$.  \hfill q.e.d.  

\vspace{3 mm} 
We next prove a lemma which will be used in the proof of Proposition \ref{vah} 
below. 
\begin{Lemma}\label{vae}
${\mathcal Ext}^1_{O_S}(\Theta_S,O_S)=0$.  
\end{Lemma} 

{\em Proof}.  
First we prove this at a non-tangential point, i.e., at a point $p$ 
where  $S$  has only normal crossings singularities.  
We apply  
${\mathcal Ext}^1_{O_S}(-,O_S)$ to the sequence (\ref{tg}) and obtain: 
\[
\ra \oplus_{l,\pm }{\mathcal Ext}^1_{O_S}(\Theta_{S^\pm _l}(-\log B^\pm _l), O_S)  
\ra {\mathcal Ext}^1_{O_S}(\Theta_S, O_S)  
\ra \oplus_\al  {\mathcal Ext}^2_{O_S}(\Theta_{B_\al}(0+\infty), O_S) 
\] 
It suffices to show that 
${\mathcal Ext}^1_{O_S}(\Theta_{S^\pm _l}(-\log B^\pm _l), O_S)=0$  and 
${\mathcal Ext}^2_{O_S}(\Theta_{B_\al}(0+\infty), O_S)=0$ at $p$. 
We prove this when $p$  is a general point, leaving similar arguments 
to the reader at four triple points.  
Let  $S_\al , \al =1,2$, be the irreducible components of $S$ 
passing through $p$ with structure sheaves $O_\al$ and put $D=S_1\cap S_2$, 
the singular locus of  $S$  at $p$. 
Locally at $p$, $\Theta_{S^\pm _l}(-\log B^\pm _l) \cong  O_1\oplus O_2\oplus O_S$  and 
$\Theta_{B_\al}((0+\infty))\cong \Theta_{B_\al}\cong O_D, D=S_1\cap S_2$. 
Therefore what we have to check is that 
${\mathcal Ext}^1_{O_S}(O_{S_\al}, O_S)=0$  and 
${\mathcal Ext}^2_{O_S}(O_D, O_S)=0$.  
For this we consider the short exact sequences 
\begin{gather}\label{shho}
 0 \ra I_\al  \ra O_S \ra  O_\al  \ra 0 \\ 
 0 \ra I_D \ra O_S \ra  O_D \ra 0 
\end{gather}
and the associated long ${\mathcal Ext}$-exact sequences. 
The desired assertion then follows from the following facts: 

1) 
$ {\mathcal Hom}_{O_S}(O_S ,O_S)  \ra  {\mathcal Hom}_{O_S}(O_\al , O_S) $ 
is surjective, being isomorphic to the quotient map  $O_S\ra O_\al ',  
\{\al , \al ',\} = \{1,2\}$ 

2) ${\mathcal Ext}^i_{O_S}(I_D, O_S)=0, i\geq 1$ (cf.\ \cite[Lemma 2.8]{fr}).  
  
\vspace{2 mm} 
Thus 
${\mathcal Ext}^1_{O_S}(\Theta_S,O_S)$ has support in the tangential points. 
But at any of these points $p$ we can find an exact sequence 
\[ 0 \ra O_S^2 \ra \Theta_S \ra  {\mathcal Q}  \ra 0\] 
where  ${\mathcal Q}$ has support in  $p$.  
Then it is easily seen that 
${\mathcal Ext}^1_{O_S}({\mathcal Q},O_S)=0$.  
(Take an exact sequence 
\[ 0 \ra F' \ra F \ra  {\mathcal Q}  \ra 0\] 
with some coherent $O_S$-modules with  $F$  free.  Then 
${\mathcal Hom}_{O_S}(F,O_S) \ra  {\mathcal Hom}_{O_S}(F',O_S)$ is surjective 
since ${\mathcal Q}$ has support in $p$ 
and $S$ is weakly normal.) 
By applying ${\mathcal Ext}^1(-,O_S)$ to the above exact sequence we get 
the desired vanishing of ${\mathcal Ext}^1_{O_S}(\Theta_S,O_S)$.  
\hfill q.e.d.

\begin{Proposition}\label{vah}
The map  $\delta$ in (\ref{s'}) is surjective, and we have  $H^2(\Theta_Z(-\log S))=0$. 
\end{Proposition} 

{\em Proof}.  By Serre dualty it suffices to show that the dual map 
$\gamma : H^1(\Omega_{Z}) \ra  Ext^1(\Theta _S,O_S)$ of  $\delta$ is injective, 
where we have 
used the isomorphism $\omega_S\cong O_S$.  
Consider the exact sequence 
\[ 
0  \ra    H^1(\Omega_S^{**}) 
   \ra   Ext^1_{O_S}(\Theta_S,O_S) 
   \ra  H^0({\mathcal Ext}^1_{O_S}(\Theta_S,O_S)) \] 
similar to (\ref{lctog}).  
By Lemma \ref{vae} 
we may identify $Ext^1_{O_S}(\Theta_S,O_S)$ with $H^1(\Omega_S^{**})$ 
and $\gamma $  with the natural map $\gamma ': H^1(\Omega_{Z}) \ra  H^1(\Omega_S^{**})$.  
Taking any of the irreducible components of  $S$, say $S^+_i$, 
we obtain a natural map  
$H^1(\Omega_S^{**}) \ra  H^1(\Omega^{**}_{S^+_i})\cong H^1(\Omega_{S^+_i})$.  
Composing $\gamma '$ with this 
we obtain the natural map 
$H^1(\Omega_{Z}) \ra H^1(\Omega_{S^+_i})$, which in turn is identified with 
the restriction map of corresponding complex cohomology groups 
$r: H^2(Z,\C) \ra H^2(S^+_i,\C)$.  

In fact we can prove the injectivity of  $r$  precisely as in the proof of 
\cite[Lemma 5.4]{fjc}, where we showed the injectivity of  $H^2(Z,\C) \ra H^2(S,\C)$ for 
a smooth member $S$ of  $|K^{-\frac{1}{2}}|$.  
Since in our case $S=S^+_i$ is an elementary surface,  
we have only to note the following: $S^+_i$ is obtained as an $m$-times blown-up of $\B^2$ 
so that $H^2(S^+_i,\C)$ is spanned by the exceptional curves and by the first Chern class. 
Thus the proposition is proved. 
(Note that the argument above is in principle similar to that for the vanishing  of 
$H^2(Z,\Theta_Z(-\log D))$ for a Joyce twistor space and a smooth element $D$ of 
$K^{-\frac{1}{2}}$ (cf.\ \cite[Th.5.1]{fjc})). \hfill q.e.d.  

\vspace{3 mm} 
By (\ref{s}) and Lemmas \ref{sll} and \ref{ts} we get: 

\begin{Corollary}\label{jg}
$h^1(\Theta_Z(-\log S)) = m-1$. 
\end{Corollary} 

Now we are in a position to prove Theorem \ref{pt}. 

\vspace{3 mm} 
{\em Proof of Theorem \ref{pt}}.  By Lemma \ref{ext} 3) 
	$H^0({\mathcal Ext}^2_{O_Z}(\Omega_{\hat{Z}}(\log \hat{S}),O_Z))=0$. 
	By 2) of the same lemma we get 
	$H^1({\mathcal Ext}^1_{O_Z}(\Omega_{\hat{Z}}(\log \hat{S}),O_Z))
	\cong H^1(O_{\hat{Q}})=0$. 
	Thus in view of (\ref{spc}) for $(X,Y)=(\hat{Z},\hat{S})$ 
the last arrow $c$ in (\ref{lctog}) is surjective. 
Since $H^2(\Theta_Z(-\log S))=0$ by Proposition \ref{vah}, 
we have the first two vanishings
by using the isomorphisms (\ref{iso}) and (\ref{iso^}). 
The sequence (\ref{lctog}) reduces to (\ref{shoo}) above, again by using 
Lemma \ref{ext}. 
Finally, from (\ref{h1})--(\ref{p2}) we get the dimensional counts 
of (\ref{sh1}) and (\ref{sh2}), and the final assertion comes from Proposition \ref{auti}. \hfill q.e.d.  

\vspace{3 mm} 
We still have to compare the deformations of the pair $(\hat{Z},\hat{S})$ 
with those of the subspaces $(\hat{S_l},\hat{C_l})$.  We start from 
the following:  

\begin{Lemma}\label{om}
There exists a natural exact sequence of $O_{\hat{Z}}$-modules: 
\begin{equation} \label{ool} 
0 \ra \Omega_{\hat{Z}}(\log \hat{S}) \stackrel{a}{\ra}  
\Omega_{\hat{Z}}(\log \hat{S}_{l'})(\hat{S}_l) \stackrel{b}{\ra} 
\Omega'_{\hat{S}_l}(\log \hat{C}_l)\otimes \hat{N_l} \ra 0 , 
\end{equation}  
where 
$\hat{N}_l:=N_{\hat{S}_l/\hat{Z}}$ is the normal bundle of $\hat{S}_l$ 
in $\hat{Z}$. 
\end{Lemma} 

{\em Proof}.  
The map $a$ is the natural inclusion.  
The map  $b$  is given by the tensor product of the natural restriction maps  
$\Omega_{\hat{Z}}(\log \hat{S}_{l'})\ra 
\Omega'_{\hat{S}_l}(\log \hat{C}_l)$ 
and  $[\hat{S}_l]\ra \hat{N_l} $.  Note that a defining equation of 
$\hat{S}_{l'}$ in $\hat{Z}$ restricts 
one of $\hat{C}_l(=\hat{S}_l\cap \hat{S}_{l'})$  
in  $\hat{S}_l$ so that the first restriction makes sense.  
This remark also implies that $b$ is surjective.  
Now what we have to show is that  Ker $b =$ Im $a$. 

Indeed, locally at points which are regular for both 
$\hat{Z}$ and  $\hat{C}_l$, the map $b$ takes the form 
\[ dx/xy,\ dy/y,\ dz/y  \ra  dx/x|\hat{S}_l,\ dy|\hat{S}_l=0,\ dz|\hat{S}_l  \] 
while the image of $a$  is generated by $dx/x,\ dy/y,\ dz$, 
where $x=0$ (resp.\ $y=0$) is 
the local equation of  $\hat{S}_{l'}$ (resp.\ $\hat{S}_l$). 
Here we consider sections of $\Omega_{\hat{Z}}(\log \hat{S}_{l'})(\hat{S}_l)$ 
as meromorphic 1-forms on  $\hat{Z}$ and identify those of 
$\Omega'_{\hat{S}_l}(\log \hat{C}_l)\otimes \hat{N_l}$ as sections of 
$\Omega'_{\hat{S}_l}(\log \hat{C}_l)$ regarding  $(1/y)|\hat{S}_l$ 
as giving the trivialization of $N_l$. The disired assertion is now 
obvious.  
By applying the product principle (\ref{w}) and (\ref{lp}) 
the same consideration applies also at singular points of $\hat{Z}$.  

It only remains to consider the tangential points.  If a local section of 
$\Omega_{\hat{Z}}(\log \hat{S}_{l'})(\hat{S}_l)$ at such a point $p$ is mapped to 
zero by  $b$, by what we have proved above it is in the image of a section $s$ of  
$\Omega_{\hat{Z}}(\log \hat{S})$ outside $p$.  
But by Proposition 7.1 the latter sheaf is reflexive at  $p$  and 
hence  $s$  extends across  $p$  as a section of the same sheaf.  The assertion is 
thus proved on the whole $\hat{Z}$.  \hfill q.e.d. 

\vspace{3 mm} 

In the same way as we get the isomorphism (\ref{isomm}) 
by using Lemma \ref{extc} 
we obtain the isomorphisms 
\begin{equation}\label{isoms} 
Ext^{i}_{O_{\hat{Z}}}(\Omega'_{\hat{S}_l}(\log \hat{C}_l)\otimes \hat{N}_l,O_{\hat{Z}}) 
\cong 
Ext^{i-1}_{O_{\hat{S}_l}}(\Omega'_{\hat{S}_l}(\log \hat{C}_l), O_{\hat{S}_l}). 
\end{equation}  

Comparing with the $Ext$ sequences associated with (\ref{ool}) 
we obtain a natural map 
\begin{equation}  \label{pp}
Ext^{1}_{O_{\hat{Z}}}(\Omega_{\hat{Z}}(\log \hat{S}),O_{\hat{Z}}) \stackrel{\al}{\ra} 
Ext^{1}_{O_{\hat{S}_l}}(\Omega'_{\hat{S}_l}(\log \hat{C}_l),O_{\hat{S}_l}). 
\end{equation}  
Together with part of the local to global sequences 
this fits into the following commutative diagram:  
\begin{equation}\label{tlo} 
\begin{array}{ccccc} 
H^1(\Theta_{\hat{Z}}(-\log \hat{S})) & 
\hookrightarrow & Ext^1_{O_{\hat{Z}}}(\Omega_{\hat{Z}}(\log \hat{S}),O_{\hat{Z}})  	
& \stackrel{\hat{c}}{\twoheadrightarrow} & 
H^0(O_{\hat{Q}})\oplus(\oplus_{p\in B} \C_p)  \\
& &  \al \downarrow   & &  u \downarrow      	 	\\
&  &  
Ext^1_{O_{\hat{S}_l}}(\Omega'_{\hat{S}_l}(\log \hat{C}_l),O_{\hat{S}_l})  	
&  \stackrel{c}{\cong}   & H^0(O_{\hat{F_l}})\oplus(\oplus_{p\in \hat{B}_l} \C_p)\\		
& &  \beta \downarrow   & &  v ||     	 	\\
&  &  
Ext^1_{O_{\hat{C}_l}}(\Omega_{\hat{C}_l},O_{\hat{C}_l})  	
&  \stackrel{d}{\cong}   & (\oplus_{p\in \bar{B}_l -\hat{B}_l} \C_p)
\oplus(\oplus_{p\in \hat{B}_l} \C_p)
\end{array}  
\end{equation}  
recalling the isomorphisms $c$ in (\ref{lctopp}) and $d$ in (\ref{exc}), 
where the top sequence is nothing but (\ref{shoo}) 
and 
$\bar{B}_l=\bar{B}^+_l\cup \bar{B}^-_l$ and 
$\hat{B}_l=\hat{B}^+_l\cup \hat{B}^-_l$ 
with $\bar{B}^\pm_l$ and $\hat{B}^\pm_l$ 
defined 
for each $\hat{S}^\pm_l$ as $\bar{B}$ and $\hat{B}$ are defined for $\hat{C}$ in Section 4. 
Note also that $c$ and $d$ are the direct sum of two isomorphisms 
\begin{equation}\label{ic} 
 c^\pm_l: Ext^1_{O_{\hat{S}^\pm_l}}(\Omega'_{\hat{S}^\pm_l}(\log \hat{C}^\pm_l),O_{\hat{S}^\pm_l})  	\cong H^0(O_{\hat{F^\pm_l}})\oplus(\oplus_{p\in \hat{B}^\pm_l} \C_p)
\end{equation}  
and 
\begin{equation}\label{id} 
d^\pm_l: Ext^1_{O_{\hat{C}^\pm_l}}(\Omega_{\hat{C}^\pm_l},O_{\hat{C}^\pm_l}) \cong 
(\oplus_{p\in \bar{B}^\pm_l -\hat{B}^\pm_l} \C_p)\oplus(\oplus_{p\in \hat{B}^\pm_l} \C_p).
\end{equation}  
By construction we have the natural identification of sets: $\hat{B}_l=B$.  
In particular the right vertical maps $u$ and $v$ give isomorphisms 
(identifications) of the second factors, while on the first factors these maps 
are isomorphic to the diagonal embeddings 
\begin{equation}\label{dgl} 
\C\ra \C^2 \ra \C^4.  
\end{equation}  

\vspace{3 mm} 
{\em Proof of Theorem \ref{ii}}.  
1) is an immediate consequence of Theorem \ref{pt} in view of 
Proposition \ref{st}.  
For 2), first note that $\al$ in (\ref{pp}) is identified with the differential of the 
versal map from the Kuranishi space  $T$  of deformations of $(\hat{Z},\hat{S})$ 
to the Kuranishi space $\hat{T}_l$ of deformations of the pair $(\hat{S_l},\hat{C_l})$, 
which is the disjoint union of $(\hat{S^\pm_l},\hat{C^\pm_l})$.  
Let  $v_{\hat{Q}}$ and $v_p, p\in B$, be the natural projections from 
$H^0(O_{\hat{Q}})\oplus(\oplus_{p\in B} \C_p)$ to $H^0(O_{\hat{Q}})\cong \C$ and 
to $ \C_p$ respectively.  Then the Zariski tangent spaces of  
$A(\hat{Q})$ and $T(p)$ are naturally identified with the kernels of 
$v_{\hat{Q}}\hat{c}$ and $v_p\hat{c}$ respectively, which are of codimension one. 

For $p\in \bar{B}_l$
let $T'_p$ be the Kuranishi space of deformations of the isolated singularity 
$(C_l,p)$, which is smooth of dimension one.  
 Then we have a versal map  $\tau '_p: (T,o) \ra (T'_p,o)$, whose differential 
is identified with  $d\beta \al$ 
composed with the projection to $\C_p$, which is a surjection.  
The inverse image  $\tau_p^{'-1}(o)$, which is 
independent of the choice of  $\tau '_p$, is easily identified with  $T(p)$ 
when $p\in \hat{B}=\hat{B}_l$ and with $A(\hat{Q})$ when $p\in \bar{B}_l-\hat{B}_l$ 
(independently of $p$). 
Thus by the properties of the diagram (\ref{tlo}) we see easily that 
 $T(p)$ and $A(\hat{Q})$ are both smooth of codimension one and 
$D$ has the properties stated in 2).  The rest of the assertions  
in 2) is obvious. \hfill q.e.d.  

\vspace{3 mm} 
{\em Proof of Theorem \ref{iit}}.  
First of all, 3) is a consequence of general theory \cite{pal} in view of 
(\ref{ex0}) in Theorem \ref{pt}.  
The Kuranishi family of $(\hat{Z},\hat{S})$ induces a deformation of 
each of the pairs $(\hat{S}^\pm_l,\hat{C}^\pm_l), l=i,j$, and 
the induced versal maps  $\tau^\pm_l  : T \ra  \hat{T}^\pm_l$  
between the corresponding Kuranishi spaces $T$ and   $\hat{T}^\pm_l$  
are submersions with 
$(\tau^\pm_l)^{-1}(A^\pm_l)=A(\hat{Q})$ 
by the diagram (\ref{tlo}) in view of (\ref{ic}) and (\ref{dgl}), 
where $A^\pm_l$ are the subspaces corresponding to $A$ in Proposition \ref{grlfb}.  
Thus $Z_t$ and $S^\pm_{l,t}$ are smooth for $t\in T-A(\hat{Q})$ and the 
structure of the surfaces $S^\pm_{l,t}$ as stated in 1) and 2) 
of the theorem is obtained from Proposition \ref{grlfb} and 2) of Remark 4.2. 
It only remains to prove the relations among $(S^\pm_{l,t})$, which is 
given in the next lemma.      \hfill q.e.d. 

\vspace{3 mm} 

\begin{Lemma}\label{wsi} 
$S^+_{l,t}$ and $S^-_{l,t}$ are isomorphic.  
$S^\pm_{i,t}$ and $S^\pm_{j,t}$ are transpositions to each other. 
\end{Lemma} 

{\em Proof}.  The weight sequence of $\tilde{S}^\pm_i$ is given by (\ref{wsl}) and similarly for 
$\tilde{S}^\pm_j$.  From the two chains between $1$ and $-1$ in (\ref{wsl}) arise 
the two cycles of $\hat{S}^\pm_l$, which in turn produces 
the two cycles  $C^{\pm ,\al}_{l,t}, \al =1,2$, on  $S^\pm _{l,t}$ via smoothing. 
The weight sequence of  $C^{\pm ,\al}_{l,t}$ are computed by (\ref{wsl}) and the formulae 
in Lemma \ref{intt}. From this already follows the first isomorphy.  For instance 
the chain $\tilde{C}^\mp_{i+1}+\cdots +\tilde{C}^\mp_j$ with weight sequence (\ref{wsl}) 
gives rise to that of $C^{\pm ,\al}_{l,t}$, 
and its toric numbering is determined by the conditions  
$\tilde{C}^\mp_{i+1}\cap H_i^\pm \neq \emptyset$ and  
$\tilde{C}^\mp_j\cap E_i^\pm \neq \emptyset$.  Since by Lemma \ref{defc} toric numbering and 
canonical numbering coincide, we see that $S^+_{l,t}$ and $S^-_{l,t}$ are isomorphic as 
we recalled in Section 3. 

The second assertion is proved  similarly as follows.  
Consider the intersections  $S^*_{i,t}\cap S^{*'}_{j,t}$, where 
$\{*,*'\}= \{+,-\}$. 
They are precisely one of the cycles contained in either of the surfaces. 
Consider for instance $S^+_{i,t}\cap S^-_{j,t}$, which is the cycle coming from the chain 
$\tilde{S}^+_i\cap\tilde{S}^-_j= \tilde{C}^-_{i+1}+\cdots +\tilde{C}^-_j$.  
Here, we have 
$\tilde{C}^-_{i+1}\cap H_i^+ \neq \emptyset$ and  
$\tilde{C}^-_j\cap E_i^+ \neq \emptyset$, while   
$\tilde{C}^-_{i+1}\cap E_j^- \neq \emptyset$ and  
$\tilde{C}^-_j\cap H_j^- \neq \emptyset$. This implies that that the toric numberings 
of  $S^+_{i,t}\cap S^-_{j,t}$ as a cycle in $S^+_{i,t}$ and in $S^-_{j,t}$ are 
reverse to each other. 
Thus again by Lemma \ref{defc} and by Section 3 we conclude that $S^+_{i,t}$ and of $S^-_{j,t}$ are 
transpositions to each other.  \hfill q.e.d. 

\vspace{3 mm} 
{\em Proof of Theorem \ref{mad}}.  
First proceed in the same way as in the proof of Theorem \ref{ma} and 
consider the universal family of log-deformations of $(\hat{Z},\hat{S})$ such that 
the given surface  $S$ is realized as a fiber $S_u, u\in I\cap T^{\sigma} -A(\hat{Q})$.  
Then we may realize  $S'$ as $S_t$ for some $t\in T -A(\hat{Q})$ which is sufficiently close to $u$. 
Note that this family is universal at any point $t$ of $T$ since 
$\dim Ext^0_{O_{\hat{Z_t}}}(\Omega_{\hat{Z_t}}(\log \hat{S_t}),O_{\hat{Z_t}}) = 0 $ 
by the upper semicontinuity of $Ext$ (cf.\ \cite{bps}).   
Now we consider this family as a germ at  $u$ and consider 
as in the proof of Theorem \ref{iit} the versal maps $\tau^\pm_i  : T \ra  \hat{T}^\pm_i$, 
but at $u$ instead of  at $o$.  
It suffices to show 
that $\tau^+_i$ maps  $T^{\sigma}$ submersively onto  $\hat{T}^+_i$ at $u$ with (smooth) 
fiber of real dimension $m$.  

$\hat{T}^+_i\times \hat{T}^-_i$ may be considered 
as the Kuranishi space of the universal deformations 
of the disjoint union $S^+_{i,u}\cup S^-_{i,u}$, where the universality is 
due to Proposition \ref{wa}.  The real structure on  $\hat{Z}$ 
interchanges $S^\pm_{i,u}$ and therefore defines a real structure on  
$S^+_{i,u}\cup S^-_{i,u}$.  Since the family over  $\hat{T}^+_i\times \hat{T}^-_i$ is universal, 
this real structure extends to the real structure on the total family 
over  $\hat{T}^+_i\times \hat{T}^-_i$.  The fixed point set $D$ of this action is clearly 
a real submanifold of dimension $m$ of  $\hat{T}^+_i\times \hat{T}^-_i$ which is mapped diffeomorphically onto each factor by the natural projections.  Moreover, $\tau^+_i \times \tau^-_i $ becomes real in the sense that it commutes with the real structure; in particular  $\tau^+_i \times \tau^-_i $ induces a smooth map $\delta : T^{\sigma}\ra D$. 

The differential of $\tau^+_i \times \tau^-_i $ at  $u$  is given by $\al$ 
in the following commutative diagram similar to (\ref{tlo}): 
\begin{equation}\label{diah} 
\begin{array}{ccccc} 
 H^1(\Theta_{Z_u}(-\log S_u)) & 
\hookrightarrow & Ext^1_{O_{Z_u}}(\Omega_{Z_u}(\log S_u),O_{Z_u})  	
&  \twoheadrightarrow & \oplus_{p\in B} \C_p  	   		\\ 
    & & \al \downarrow   & &  ||     	 	\\
 &  &  
Ext^1_{O_{S_{i,u}}}(\Omega'_{S_{i,u}}(\log C_{i,u}),O_{S_{i,u}})  	
&  \cong   & \oplus_{p\in B_i} \C_p		
\end{array}  
\end{equation}  
where  $B$ is the set of tangential points of $S_u$ and $B_i$ is the set 
of nodes of $C_{i,u}$; they are naturally identified.  
Here each term admits a natural real structure and each map is real. 
From this we immediately see that $\delta $ is a submersion and 
$T^{\sigma}$ is mapped submersively onto  $\hat{T}^+_i$ with fiber of 
real dimension equal to 
\[ \dim_{\C} H^1(\Theta_{Z_u}(-\log S_u)) 
= \dim_{\C} Ext^1_{O_{Z_u}}(\Omega_{Z_u}(\log S_u),O_{Z_u}) - 2m=m\] 
(cf.\ Proposition \ref{wa}) as desired, where we have used the constancy of 
the dimension of 
$Ext^1_{O_{Z_t}}(\Omega_{Z_t}(\log S_t),O_{Z_t})$.  
In fact, since the base space is smooth, 
it is immediate to see 
that $\Omega_{{\mathcal Z}/{\mathcal T}}(\log {\mathcal S})$ is flat ovet $T$. 
Then by the upper semiconinuity of relative $Ext$ (cf.\ \cite{bps}), we have the 
vanishing of $Ext^i_{O_{Z_t}}(\Omega_{Z_t}(\log S_t),O_{Z_t})$ for $i=0,2$, and 
then by the invariance of the alternating sum of the dimensions of $Ext^i$, we get 
the constancy of the dimension of $Ext^1$ \cite{bps}.  \hfill q.e.d.  

\vspace{3 mm} 
{\em Proof of Theorem \ref{map}}.  
We start from any of the LeBrun-Joyce twistor spaces  $Z$  
with the distinguished choice of $(i,j)$ as mentioned 
in Remark 7.3 and 
get the singular twistor space $(\hat{Z},\hat{S})$ with real structure $\sigma $. 
First of all, in order to get the universal family we fix as usual (cf.\ \cite{df}) 
any twistor line $L$  on  $Z$ other than  $L_l, 1\leq l \leq k$, and consider its proper transform $\hat{L}$ in  $\hat{Z}$.   
We then consider the Kuranishi family of log-deformations  
of the triple  $(\hat{Z},\hat{S},\hat{L})$, 
``log'' referring only to the deformations of the pair $(\hat{Z},\hat{S})$,  
which is 
universal since   Aut$_0(\hat{Z},\hat{S},\hat{L})$ now reduces to the identity.  
The Kuranishi space  $T(L)$ is a smooth fiber space over the original Kuranishi space $T$ for the deformations of the pair $(\hat{Z},\hat{S})$ with four dimensional fibers. 
We then restrict the family over the inverse image  $I(L)$  of  $I\subseteq T$.  
Then  $\sigma $ induces a canonical real structure on  $T(L)$ preserving $I(L)$.  The 
restriction to $I(L)^\sigma$ of the family of anti-self-dual bihermitian structures 
of Theorem \ref{tw} 
has as its underlying complex structures all parabolic Inoue surfaces 
with Betti number $m$. 

We consider the universal family of {\em log}-deformations of $(\hat{S}^\pm_i,\hat{C}^\pm_i)$ constructed in 2) of Proposition \ref{grlfb}.  Then the product $\hat{T}^+_i\times \hat{T}^-_i$ of the corresponding Kuranishi spaces  $\hat{T}^\pm_i$ is naturally considered as parametrizing the universal family of deformations of the disjoint union of $(\hat{S}^\pm_i,\hat{C}^\pm_i)$.   Since this family is universal and $\sigma $ interchanges both pairs 
inducing a real structure on this disjoint union, there exists a natural action of 
$\sigma $ on $\hat{T}^+_i\times \hat{T}^-_i$ interchanging the two factors.  Let  $D$  be the associated real part.  

  Now as in the proof of Theorem \ref{mad} 
we get a versal map 
$\tau : I(L) \ra  \hat{T}^+_i\times \hat{T}^-_i$  
which induces a smooth  map of the real part: $I^\sigma (L) \ra  D$. 
The map is of rank one at the origin $ o$ and the image of its differential is 
mapped surjectively onto both factors  
as in the proof of Theorem \ref{iit}. 
Therefore, $\tau$ has rank at least one also at all nearby points $u$ of $o$ 
and submersive onto the both factors. 
This implies that the image of $\tau|I(L)^{\sigma}$ contains a (local) real smooth curve  $K$ 
contained in  $D\cap \tau(I(L))$ whose images on both factors of $\hat{T}^\pm_i$  again are 
real smooth curves $K^\pm_i$.   Let  $S$  be any parabolic Inoue surface corresponding 
to a point $\kappa$ of $K^\pm_i$.  Then the family of anti-self-dual bihermitian structures on  
$M[m]$ restricted to a suitable $m$-dimensional submanifold of 
$\tau^{-1}(\kappa)$ has the desired properties.  

Finally we show that $S^+_i$ and $S^+_j$ are isomorphic. 
We may assume that the intersection $E_{t}:=S^+_{i,t}\cap S^+_{j,t}$ is 
an elliptic curve.  (Otherwise we have only to replace  $S^+_{j,t}$ by 
$S^-_{j,t}$ and define $J_{2,t}$ via $S^-_{j,t}$.) 
The twistor fibration  $Z_t \ra  M[m]$  induces the isomorphism of the 
fundamental groups of these spaces.  
Since the induced projection  $S^+_{l,t} \ra  M[m],l=i,j$, 
is diffeomorphic, the inclusion  $S^+_{l,t} \hookrightarrow Z_t$  also 
gives the isomorphism of fundamental groups.  
Thus, if $r: \tilde{Z}_t \ra  Z_t$ is the universal covering, 
the induced map $\tilde{S}^+_{l,t} = r^{-1}(S^+_{l,t}) \ra  S^+_{l,t}$ is 
the universal covering of $S^+_{l,t}$ also, 
and $\tilde{E}_t:=r^{-1}(E_t) \ra  E_t$ gives 
the common infinite cyclic unramified covering of $E_t$ 
for both of $S^+_{l,t}, l=i,j$. 
Hence by Lemma \ref{prbi} we conclude that $S^\pm_i$  are isomorphic.  
\hfill q.e.d.  

\vspace{2 mm} 
\section{Anti-self-dual hermitian structures on half Inoue surfaces}  

\vspace{3 mm} 

In the construction of the pair  $(\hat{Z},\hat{S})$ in Section 5 using 
the identification map $\varphi$ of (\ref{fai}),  
we may also use a map $\varphi$ of {\em twisted type} of three kinds 
in the sense that 
$\varphi $ maps  $(H_i^\pm, E^\pm_j)$ 
to $(E^\pm_{i}, H^\mp_j)$  
(resp.\ $(E^\mp_{i}, H^\pm_j)$, resp.\ $(E^\mp_{i}, H^\mp_j)$)
(instead of to  $(E^\pm_i, H_j^\pm)$). 
We call such a  
$\varphi$  $i$-{\em twisted} (resp.\ $j$-{\em twisted}, resp.\ {\em bi-twisted}) 
in compatible with the terminology in Section 4.  
The main geometric implication of these variations is that if e.g. $\varp$ is 
$i$-twisted, $\hat{S}_j$  becomes connected, while  $\hat{S}_i$  
consists of two connected components $\hat{S}^\pm_i$ as before.  
Similarly, if $\varp$ is bi-twisted, both $\hat{S}_i$ and $\hat{S}_j$  are connected. 

In this case we are led to anti-self-dual hermitian structures 
on half Inoue surfaces and to anti-self-dual 
bihermitian structures on hyperbolic Inoue surfaces 
on their unramified double coverings. 
(There are no `parabolic' case  
unlike in the untwisted case.) 

To explain this, 
we first note that most of the constructions and results in Case-H 
in Section \ref{ldt} are also valid for this case without any change.  
(The relevancy of Case-H comes from Proposition \ref{au}.) 
Indeed, the cohomological computations in Section 8 are either local along $\hat{Q}$ 
or those on  $Z$ or $\tilde{Z}$ for which  $\varp$ plays no role, 
and hence we get the same result also in this case.  
In particular, the obstruction for the log-deformations 
for the pair  $(\hat{Z},\hat{S})$  vanishes and we get the Kuranishi family 
\begin{equation}\label{krr'} 
 g: ({\mathcal Z} ,{\mathcal S} ) \ra T,\ (Z_o,S_o)=(\hat{Z},\hat{S}),\ o\in T 
\end{equation}  
of log-deformations of $(\hat{Z},\hat{S})$ 
with the properties in Case-H of Theorem \ref{ii}. 
In particular  $T$  is smooth of dimension $3m$. 
The main difference now lies in 
the structure of the surfaces  $S_t$ for $t\in T-A(\hat{Q})$.  
(Here and in what follows we use the notations of Section 7.) 
For  $l=i$ or $j$ denote by $l'$ the complementary index with $\{l,l'\}=\{i,j\}$ as before.  
Then also in this case, by the arguments in \cite[\S 3,\S 4]{nkm1} 
we get easily the following: 
\begin{Lemma}\label{ll}
Let $t$ be any point of $T-A(\hat{Q})$. If  $\varp$ is $l$-twisted, then 
the deformation $S_{l',t}$  of  $\hat{S}_{l'}=\hat{S}^+_{l'}\cup \hat{S}^-_{l'}$ 
is a connected smooth surface of class VII with second Betti number $2m$, while 
the deformation $S^\pm_{l,t}$  of  $\hat{S}^\pm_{l}$ are 
smooth disjoint surfaces of class VII with second Betti number $m$. 
Similarly, if $\varp$ is bi-twsited, the conclusion for $\hat{S}_{l'}$ 
above holds for both $S_{i,t}$ and $S_{j,t}$.  
\end{Lemma} 

Using this lemma and Proposition \ref{grlfb} in Case-H$'$, 
the more precise structure of the surfaces $S_{l,t}$ and $S^\pm_{l,t}$ are deduced 
as in the proof of Theorem \ref{iit}.   
Namely using the notations of Section 7 we have the following: 
\begin{Theorem}\label{iih}
1) Assume that $t \in T-A(\hat{Q})$. In the $l$-twisted case 
the fibers $Z_t$, $S^\pm_{l,t}$ and $S_{l',t}$ are all smooth. 
The minimal model $\bar{S}^\pm_{l,t}$ of $S^\pm_{l,t}$ is 
either a half Inoue surface or a diagonal Hopf surface, 
while the minimal model $\bar{S}_{l',t}$ of $S_{l',t}$ is 
either a hyperbolic Inoue surfaces or a diagonal Hopf surface. 
Diagonal Hopf case occur if and only if  $t \in T-D$. 
In the bi-twisted case the statements for 
$S_{l',t}$ above holds for both $S_{i,t}$ and $S_{j,t}$.  

2) Assume that $t \in I-A$. 
If $\varp$ is $l$-twisted, 
$S^\pm_{l,t}$ is a 
properly blown-up half Inoue surface, while 
$S_{l',t}$ is a properly blown-up hyperbolic Inoue surface 
which is 
the unique double covering of the transposition ${}^{\mbf}S^\pm_{l,t}$ of $S^\pm_{l,t}$. 
If $\varp$ is bi-twisted, 
$S_{i,t}$ and $S_{j,t}$ are transpositions of 
each other.  
Moreover, the complex surfaces above are independent of  $t$ up to isomorphisms. 

3) 
The Kuranishi family $g$ is universal.  
\end{Theorem} 
 By restricting the family obtained in the above theorem 
to the real parts  $T^\sigma$ of $T$ and  $I^\sigma$ of  $I$ respectively, 
we can immediately deduce the coclusions similar to Theorem 7.4, 7.5 and 7.6 
in the same way as we obtained these theorems from Theorem \ref{iit}. 
Here we state only an analogue of Theorem \ref{ma}, 
leaving the reader to formulate the analogues of 
Theorems \ref{tw} and \ref{mad}.  
In fact, 
considering the deformations $S^\pm_{l,t}$ of  $\hat{S}^\pm_l$ 
over $I^\sigma$ in the $l$-twisted case we now obtain: 

\begin{Theorem}\label{mah}
 Let  $S$ be an arbitrary properly blown-up half Inoue surface 
with second Betti number $m$.  Then 
there exists a real $m$-dimensional family 
of anti-self-dual hermitian structures on  $S$. 
\end{Theorem}  

For the statement of the next result 
we introduce the following terminology.  
Let  $a: M[2m] \ra  M[m]$  be the unique unramified double covering. Denote the 
covering involution by $\kappa $.   
A {\em bihermitian structure with anti-holomorphic involutions} 
on the pair $\langle M[2m],\kappa \rangle$ is by definition 
a bihermitian structure $([g], J_1, J_2)$  on  $M[2m]$ 
such that 
on  each $S_i:=(M[2m],J_i)$,  $\kappa $ is anti-holomorphic.  
If instead, $\kappa $ is holomorphic on  $S_1$ 
and anti-holomorphic on  $S_2$, 
it is called a {\em bihermitian structure 
with holomorphic and anti-holomorphic involutions}.  
Then 
from the deformations of $S_{l,t}$  of  $\hat{S}_{l}$,  
in the bi-twisted (resp.\ $l$-twisted) case 
we have the following:  

\begin{Theorem}\label{mahh}
Let $S$ be any properly blown-up half Inoue surface  
and $\tilde{S}$ 
the properly blown-up hyperbolic Inoue surface which is an unramified double covering of $S$ 
with Galois involution $\iota $.   Then 
there exists a real $m$-dimensional family $([g]_t, J_{1,t}, J_{2,t})$ 
of anti-self-dual bihermitian structures with anti-holomorphic 
(resp.\ holomorphic and anti-holomorphic) involutions 
on  $\langle M[2m],\kappa \rangle $ such that $(M[2m],J_{1,t}) \cong \tilde{S}$ 
(resp.\ $\langle (M[2m],J_{1,t}), \kappa \rangle  = \langle \tilde{S},\iota \rangle $)  
and $(M[2m],J_{2,t})\cong {}^{\mbf}\tilde{S}$, the transposition of $\tilde{S}$, 
independently of $t$. 
\end{Theorem} 
{\em Proof}.  
We consider the family obtained from $S$ as in Theorem \ref{iih}. 
Assuming that  $\varphi$ is bi-twisted, we shall show the existence of the family 
in the case of anti-holomorphic involutions, the other case being shown similarly by 
starting with $l$-twisted $\varphi$. 
Now if we restrict the obtained family to 
$U:=I^\sigma -A(\hat{Q})$ we get a smooth family  $\{[g]_t\}_{t\in U}$ of anti-self-dual structures 
on $M=M[m]$, 
and the associated family of twistor spaces $\{Z_t\}_{t\in U}$D 
By Lemma \ref{ll} 
we have surfaces $S_{l,u},\ l=i,j$, in $Z_t$ 
which are $\sigma$-invariant hyperbolic Inoue surfaces and 
are transpositions to each other  
such that the restriction of the twistor fibration $b_t: Z_t \ra M$ makes  $S_{l,u}$  
smooth unramified double coverings of  $M$.  

Take the natural fibered product 
$\tilde{b}_t:\tilde{Z}_t:=Z_t\times_M \tilde{M}\ra \tilde{M}:=M[2m]$. 
$\tilde{Z}_t$ is an unramified double covering of  $Z_t$ which is the twistor space of 
the induced anti-self-dual structure  $(\tilde{M},[\tilde{g}_t])$. 
Moreover, 
the inverse image of $S_{l,t}$ in  $\tilde{Z}_t$  is a disjoint union of two copies of $S _{l,t}$, denoted by  $\tilde{S}^\pm _{l,t}$. 
These are $\tilde{\sigma}$-conjugate 
to each other and are mapped diffeomorphically on to $\tilde{M}$, 
where $\tilde{\sigma}=\tilde{\sigma}_t$ is the real structure of $\tilde{Z}_t$. 
Moreover,  $\tilde{S}^\pm _{j,t}$ is the transposition of $\tilde{S}^\pm _{i,t}$ 
by Lemma \ref{wsi}. 
Thus in  $\tilde{Z}_t$  we get two pairs of $\tilde{\sigma}$-invariant elementary surfaces 
$\{\tilde{S}^\pm _{l,t}\},\ l=i,j$, giving rise to anti-self-dual bihermitian structures on  
$(\tilde{M},[\tilde{g}]_t)$. 

  It remains to see that this latter structures are actually those on $(\tilde{M},\kappa )$. 
In fact, since  $\kappa $  preserves the anti-self-dual structures $[\tilde{g}]_t$, 
it lifts to a biholomorphic automorphism  $\tilde{\kappa }$ of  $\tilde{Z}_t$  which interchanges  $\tilde{S}^\pm _{l,t}$. 
Then the compsition  $\tilde{\sigma }\tilde{\kappa }$  preserves $\tilde{S}^\pm _{l,t}$ and induces an anti-holomorphic 
involution on each, which is also a lift of  $\kappa $.  Since  $b_t$ induces  
an $(\tilde{\sigma }\tilde{\kappa },\kappa )$-equivariant isomorphism 
$\tilde{S}^+_{l,t}\cong (\tilde{M},\tilde{J}_{l,t})$, 
where $\tilde{J}_{l,t}$ is the pull-back of $J_{l,t}$ to $\tilde{M}$. 
Then by the definition of  $\tilde{J}_{l,t}$, we are done. 
\hfill q.e.d.

\vspace{3 mm} 
{\em Remark 9.1}.  
1) Using Proposition \ref{grlfb} in Case-H$'$, we can also obtain an analogue 
of Theorem \ref{mad} above, in which we get a family of 
anti-self-dual structures on certain blown-up half Inoue surfaces and 
blown-up diagonal Hopf surfaces.  

2) As in Case-H, for a fixed  $S$  as in Theorem \ref{mah}, 
we can construct other families with the same properties 
by starting from suitable other choice of 
$K$-actions on  $m\B^2$ and pairs $(i,j)$ 
in the construction of Section 5.  

3) The family of anti-self-dual bihermitian structures on $M[2m]$ 
in Theorem \ref{mahh} could possibly be connected by deformations to those 
obtained in 2) of Theorem \ref{iit} with $m$ replaced by  $2m$ there. This type of 
relations would deserve further study.  

4) Let $\tilde{S}\ra S$ and $\iota$ be as in the theorem. 
The theorem implies that $\tilde{S}$ always 
admits a fixed point free anti-holomorphic involution, which could be 
identified with the anti-holomorphic involution $\iota\mu$, where $\mu$ is the 
real structure of $\tilde{S}$ defined in Lemma \ref{rsi}.

\section{Differential geometric consequences} 

 Our interest in bihermitian metrics comes from anti-self-duality
and this case was first treated in \cite{pnt6} motivated by
questions of Salamon \cite{sa94} concerning existence of orthogonal
(integrable) complex structures in a given conformal class $[g]$;
here we explicitly want to exclude the well known case of
hyperhermitian structures.

More impetus came from the work of \cite{agg} who considered the
general four-dimensional case; new examples have been recently found
by Hitchin \cite{hi07} on del Pezzo surfaces in the context of
generalized K\"ahler manifolds as introduced by Gualtieri \cite{gu04}.
In dimension four this condition amounts to say that there is a
metric $g\in[g]$ which is Gauduchon for both $J_i$, $i=1,2$, and for
which the sum of the Lie forms $\theta_i$ vanishes:
$$\theta_1 +\theta_2=0 \quad\textrm{ and }\quad
\delta\theta_1=0=\delta\theta_2.$$

The above two equations always hold for a bihermitian metric on a
compact anti-self-dual four-manifold (we excluded the hyperhermitian case) 
\cite[Prop.3.5]{pnt6} and their twistor correspondence is that the
\it real \rm degree-4 divisor defined by $\{\pm J_i$, $i=1,2\}$ is
an anticanonical divisor of $Z$ \cite[Lemma 3.4.]{pnt6}.

By Lemmas \ref{aatc} and \ref{atc} an anti-self-dual bihermitian surface $S$ 
with odd first Betti number must 
be a blow-up of either a hyperbolic or parabolic Inoue or Hopf surface.  
This condition is therefore a necessary condition 
for the existence of anti-self-dual bihermitian metrics. 

However, Lemma 3.2 was known to Nakamura (unpublished), 
and the same statement can also be found in \cite[II,2]{agg} 
(in their proof Lemma (2.8) should be replaced by Lemma (2.7));  
it also appears in \cite[2.29]{dl} 
where a more general result is proved concerning 
surfaces with a global spherical shell.  We also would 
like to give here a complete proof following the work of 
Nakamura.

\vspace{3 mm} 
{\em Proof of Lemma 3.2 }
 When $b_2(S)=0$ Bombieri-Inoue surfaces have no
curves, therefore by Bogomolov theorem \cite{lyz}\cite{te94} $S$ is
a Hopf surface and is diagonal because $-K$ is disconnected. We can
therefore assume that $S\in$ VII$_0^+$ and by \cite[12.4]{na84}, $S$
contains a cycle $C$ of rational curves.
 Suppose that \cite[2.2]{nato} holds (i.e. $S$ has a branch) then 
\cite[3.1]{nato} applies with $m=1$ and $F$ the trivial line bundle
and Nakamura concludes that anti-canonical divisor itself is connected, which is 
absurd.
 The only other possibility is that \cite[2.2]{nato} does not hold, in
which case $S$ is an Enoki surface (but $-K$ has no divisor in this
case) or a half-Inoue (but $-K$ is connected in this case) or else
$S$ must be parabolic or hyperbolic Inoue.   \hfill q.e.d.

\vspace{3 mm} 
The results of Section 7 show that the above mentioned necessary condition is 
actually sufficient at least for properly blown up hyperbolic Inoue surfaces 
and also for some parabolic Inoue surfaces.  To summarize 
we have the following:  

\begin{Proposition} All metrics constructed in Section 7 on (blown-up) hyperbolic or parabolic Inoue surfaces or Hopf surfaces are (twisted) generalized K\"ahler. 
\end{Proposition} 

\vspace{1 mm} 
We also have applications to an old and basic question of Vaisman
who asked \cite{va87}: which compact complex surfaces $S$ can admit
locally conformally K\"ahler (l.c.K.) metrics ?  (cf.\ \cite{po08} for more details.) 

By a result of Tricerri \cite{tr82} one can assume that $S$ is
minimal and after the work of Belgun \cite{be00} the answer is
positive for all the locally homogeneous surfaces, i.e. all 
surfaces in class VI, Kodaira surfaces and 
all surfaces such that  $b_2=0$ and $b_1 =1$, except 
for the complement of a real line in a complex $1$-dimensional 
family of certain Bombieri-Inoue surfaces which do not admit l.c.K. 
metrics at all. 
Therefore, Vaisman question remains open only for surfaces in class-
VII$_0^+$.

Now, by a theorem of Boyer \cite{by6} anti-self-dual hermitian metrics are
automatically l.c.K. on a compact complex surface. Therefore we have 

\begin{Theorem} All the surfaces of class VII$_0^+$ in Theorems 
7.5, 7.6, 7.7, 9.3 and 9.4 have l.c.K. metrics. 
\end{Theorem}

We conclude with the following

{\em Remark 10.1}.  
Except for the anti-self-dual hermitian metrics on parabolic Inoue surfaces by
LeBrun \cite{lb}, all the examples in the proposition and the theorem above are new 
and are the only known examples in class-VII$_0^+$. 

\vspace{3 mm} 
{\em Note}. While undergoing the final draft of this work 
M.\ Brunella communicated to us that he constructed l.c.K.\ 
metrics on all Enoki surfaces \cite{br}. 

\vspace{3 mm} 
\noindent {\small \bf Acknowledgements. \rm We would like to thank 
several institutions for kind hospitality during the preparation of
this work: Instituto de Matem\'aticas, Cuernavaca; I.C.T.P. and
S.I.S.S.A., Trieste; Universit\'e de Provence, Marseille; Centro De 
Giorgi, Pisa.}

\vspace{3 mm} 

\begin{flushright}
Department of Mathematics\\
Graduate School of Science\\
Osaka University\\
Toyonaka 560-0043 Japan \\ 
{\small E-mail address: fujiki@math.sci.osaka-u.ac.jp} \\ 
Dipartimento di Matematica\\
Universit\'a  Roma Tre\\
00146 Roma, Italy\\
{\small E-mail address: max@mat.uniroma3.it}   
\end{flushright}

\end{document}